\documentclass{amsart}
\usepackage{latexsym}
\usepackage{amsmath}
\usepackage{mathrsfs}
\usepackage{amssymb}
\usepackage[all]{xy}
\usepackage{comment}

\usepackage{autobreak}
\usepackage{enumitem}
\usepackage[hypertexnames=false]{hyperref} 

\usepackage[\ifdefined\enabletodonotes\else disable\fi]{todonotes}
\setuptodonotes{fancyline}

\usepackage{autonum} 

\usepackage{tikz}
\usetikzlibrary{positioning}
\usetikzlibrary{arrows.meta}
\usepackage{tikz-cd}

\numberwithin{equation}{section}
\newtheorem{thm}{Theorem}[section]
\newtheorem{prop}[thm]{Proposition}
\newtheorem{cor}[thm]{Corollary}
\newtheorem{lem}[thm]{Lemma}

\theoremstyle{definition}
\newtheorem{defn}[thm]{Definition}

\newtheorem{assertion}[thm]{Assertion}
\theoremstyle{remark}
\newtheorem{rem}[thm]{Remark}
\newtheorem{ex}[thm]{Example}

\newcommand{\K}{{\mathbb K}}
\newcommand{\Q}{{\mathbb Q}}
\newcommand{\R}{{\mathbb R}}
\newcommand{\Z}{{\mathbb Z}}
\newcommand{\C}{{\mathcal C}}
\newcommand{\e}{\varepsilon}
\newcommand{\A}{\mathcal A}
\renewcommand{\ker}{\operatorname{Ker}}
\newcommand{\im}{\operatorname{Im}}

\newcommand{\mapright}[1]{%
 \smash{\mathop{%
  \hbox to 1cm{\rightarrowfill}}\limits_{#1}}}
\newcommand{\maprightd}[2]{%
 \smash{\mathop{%
  \hbox to 0.5cm{\rightarrowfill}}\limits^{#1}\limits_{#2}}}
\newcommand{\mapleft}[1]{%
 \smash{\mathop{%
  \hbox to 1cm{\leftarrowfill}}\limits_{#1}}}
\newcommand{\mapleftu}[1]{%
 \smash{\mathop{%
  \hbox to 0.8cm{\leftarrowfill}}\limits^{#1}}}
\newcommand{\maprightu}[1]{%
 \smash{\mathop{%
  \hbox to 1cm{\rightarrowfill}}\limits^{#1}}}
\newcommand{\maprightud}[2]{%
 \smash{\mathop{%
  \hbox to 1cm{\rightarrowfill}}\limits^{#1}_{#2}}}
\newcommand{\mapleftud}[2]{%
 \smash{\mathop{%
  \hbox to 1cm{\leftarrowfill}}\limits^{#1}_{#2}}}

\newcommand{\utot}[1]{s^{#1}}

\author[K. Kuribayashi]{Katsuhiko Kuribayashi}
\address{%
  Department of Mathematics,
  Faculty of Science,
  Shinshu University,
  Matsumoto, Nagano 390-8621, Japan
}
\email{kuri@math.shinshu-u.ac.jp}

\author[T. Naito]{Takahito Naito}
\address{%
  Nippon Institute of Technology,
  Gakuendai, Miyashiro-machi, Minamisaitama-gun, Saitama 345-8501, Japan
}
\email{naito.takahito@nit.ac.jp}

\author[S. Wakatsuki]{Shun Wakatsuki}
\address{%
  Graduate School of Mathematics, Nagoya University,
  Furo-cho, Chikusa-ku, Nagoya, Aichi 464-8601, Japan
}
\email{shun.wakatsuki@math.nagoya-u.ac.jp}

\author[T. Yamaguchi]{Toshihiro Yamaguchi}
\address{%
  Faculty of Education,
  Kochi University, Akebono-cho, Kochi 780-8520, Japan
}
\email{tyamag@kochi-u.ac.jp}

\begin{document}
\title[The cohomology interleaving distance]{The equalities of interleaving distances and cohomology interleavings of spaces over $BS^1$ 
}

\footnote[0]{{\it 2020 Mathematics Subject Classification}. Primary 55N31, 55U15; Secondary 55P62, 55U35.

{\it Key words and phrases.} Interleaving distance, barcode, persistence differential graded module, Sullivan model, cup-length.


}


\maketitle
\todo{[RR2]: The title is changed reflecting the contribution of our article.}

\begin{abstract}
The cohomology interleaving distance (CohID) is defined and considered in the category of persistent differential graded modules over a field.
As a consequence, we show that, in the category, the distance coincides with the {\it homotopy commutative interleaving distance}, the {\it homotopy interleaving distance} originally due to Blumberg and Lesnick, and
the {\it interleaving distance in the homotopy category} in the sense of Lanari and Scoccola. 
Moreover, we apply the CohID to spaces over the classifying space $BS^1$ of the circle group via the singular cochain functor.
Then, upper and lower bounds of the CohID are investigated with the cup-lengths of spaces over $BS^1$.
As a computational example, we explicitly determine the CohID for complex projective spaces by utilizing the bottleneck distance of barcodes associated with the cohomology of the spaces.
\end{abstract} \todo{[RR2] Abstract is revised with `defined' and `over a filed'.}

\tableofcontents

\section{Introduction}\label{sect:Introduction}
Over the last decades, persistence theory, in the context of topological data analysis, has developed rapidly
through the study of persistent homology and representations of algebras.
More recently, persistence objects values in a category $\mathcal {C}$, namely, objects in the functor category $\mathcal{C}^{(\R, \leq)}$, are investigated from the homological and homotopical points of view.
Indeed,  Blumberg and Lesnick \cite{Bl-L} introduce the {\it homotopy interleaving distance} $d_{\rm{HI}}$ and the {\it homotopy commutative interleaving distance} $d_{\rm{HC}}$ for persistence objects valued in the category of topological spaces.
In \cite{LS},  Lanari and Scoccola define the  {\it interleaving distance in the homotopy category}, denoted $d_{\rm{IHC}}$, in the functor category
$\mathcal{M}^{(\R, \leq)}$ for a cofibrantly generated model category $\mathcal{M}$.
Moreover, the distances $d_{\rm{HI}}$ and $d_{\rm{HC}}$ are generalized to be applicable to persistence objects valued in $\mathcal{M}$.
Here the term {\it distance} means an extended pseudo-metric on a class. As a consequence, we see that each distance $d$ for persistence objects $M$ and $N$ is a numerical two-variable homotopy invariant; that is, $d(M, N)=d(M', N')$ whenever $M$ and $N$ are weakly equivalent to $M'$ and $N'$, respectively; see  Subsection \ref{sect:ID_up_to_H} for more details. 

By a categorical consideration, for example, \cite[Proposition 3.6]{B-S}, it is readily seen that
\begin{eqnarray}\label{eq:HI}
d_{\text{HC}} \leq d_{\text{IHC}} \leq d_{\text{HI}}
\end{eqnarray}
on the class of objects in $\mathcal{M}^{(\R, \leq)}$. Furthermore, the result \cite[Theorem A]{LS} asserts that $d_{\text{HI}}\leq 2d_{\text{HC}}$.
Let $\mathsf{Top}$ be the category of topological spaces. Then, by \cite[Proposition 3.12]{LS}, we see that if there exists a positive number $c$ such that $d_{\text{HI}} \leq c \cdot d_{\text{HC}}$ on $\mathsf{Top}^{(\R, \leq)}$, then $c\geq \frac{3}{2}$. 
Hence $ d_{\text{HC}} \neq d_{\text{HI}}$ in general. \todo{[RR1, 2.4 Comment 1] cdot $\cdot$ is put. To clarify $c$, the sentence is changed.}
Remarkably, this result is proved by using the notion of the Toda bracket; see \cite[Section 3.2]{LS}.  It is worthwhile to note that a positive answer to a version of the persistent Whitehead conjecture \cite[Conjecture 8.6]{Bl-L} is given by deeply
considering interleavings in the homotopy category $\text{Ho}(\mathsf{Top}^{(\R, \leq)})$; see \cite[Theorem B, Remark 5.14]{LS}.

In this article, we introduce the {\it cohomological interleaving distance} $d_{\rm{CohI},\K}$ for persistence objects valued in $\mathsf{Ch}_\K$ the category of differential graded (dg) modules over a field $\K$. Then, we show the equalities
\[
d_{\rm{CohI},\K} = d_{\text{HC}} =d_{\text{IHC}} =d_{\text{HI}}
\]
for objects in $\mathsf{Ch}_\K^{(\R, \leq)}$
; see Theorem \ref{thm:Equalities} for more details. Thus, we may compute the distances 
 of persistence dg modules
by using the {\it bottleneck distance} of barcodes associated with the cohomology groups of the given persistence objects.

In the context of topological data analysis, 
persistence objects usually arise from sublevel sets of maps, or  certain filtrations of simplicial complexes.
In the latter half of this article, we will study persistence objects arising from spaces over $BS^1$,
specifically examining the cohomology interleaving distance $d_{\rm{CohI},\K}$ between such spaces, where $BS^1$ denotes the classifying space of the circle group $S^1$. More precisely, we bring a space over $BS^1$ to the category $\mathsf{Ch}_\K^{(\R, \leq)}$ of persistence dg modules via the singular cochain functor $C^*( \ ; \K)$ with coefficients in a field $\K$; see Section \ref{sect:BS^1}. As a consequence, we have an extended pseudodistance
\[
d_H : [X, BS^1]\times [X, BS^1]\to \R_{\geq 0}\cup\{\infty\}
\]
on the homotopy set $[X, BS^1]$ for a space $X$ by using $d_{\rm{CohI},\K}$; see Proposition \ref{prop:HomotopySet}.
We observe that, homotopically speaking, maps into $BS^1$ are the same thing as principal $S^1$-bundles; see \cite[Theorem 4.5.21]{Tamaki} or \cite[Chapter 4, Sections 12 and 13]{Husemoller}.

As Theorem \ref{thm:Equalities} aforementioned, an algebraic result (Theorem \ref{thm:IHC=CohI}) yields that the cohomology interleaving distance between spaces over $BS^1$ coincides with the homotopy interleaving distances $d_{\text{HC}}$, $d_{\text{IHC}}$ and $d_{\text{HI}}$ for the spaces.
Moreover, 
Propositions \ref{prop:cup} and \ref{prop:cup_2} allow us to obtain upper and lower bounds of the cohomology interleaving distance of spaces over $BS^1$ with the cup-lengths of the spaces. 

We have many computational examples of the distance $d_{\rm{CohI},\K}$. Some of them enable us to realize the triangle inequality of the distance. For instance, let
$S^3\to \mathbb{C}P^1$ be the Hopf bundle and $Z_j' \to Z_j$  \todo{[RR1, 2.4 Comment 2] We use the terminology `the Hopf bundle'.}  
the $S^1$-bundle described in Proposition \ref{prop:models} and Remark \ref{rem:realization} for $j = 0$ and $1$. Observe that the total space $Z_j'$ has the rational homotopy type of $S^3 \times S^3 \times S^7$ for each
$j = 0, 1$.
We regard the base spaces of the bundles as spaces over $BS^1$ with the classifying maps.
Then, by Propositions \ref{prop:M_pt}, \ref{prop:X_pt}, \ref{prop:models} and \ref{prop:M_CP}, we have the tetrahedron (\ref{eq:triangle}) below.
It also follows that the distance $d_{\rm{CohI},\Q}$ of the spaces in (\ref{eq:triangle}) and the Borel construction of the free loop space of a simply-connected space are infinite; see Example \ref{ex:ThreeClasses}.

\medskip
\begin{minipage}{0.7\hsize}
\begin{eqnarray}\label{eq:triangle}
\xymatrix@C9pt@R16pt{
\mathbb{C}P^1 \ar@{-}[rrrrrr]^-{d(\mathbb{C}P^1,  Z_1)=\frac{7}{2}}  & & & & &  &Z_1 \\
\text{pt} \ar@{-}[u]^{d({\rm pt}, \mathbb{C}P^1)=1}  \ar@{-}[rrrrrru]^(0.5){d({\rm pt},  Z_1)=\frac{7}{2} \ \ \ }  \ar@{-}[rrd]_-{d({\rm pt}, Z_0)=2}&&& & & \\
&&  Z_0 \ar@{-}[rrrruu]_-{\ d(Z_0, Z_1)=3} \ar@{-}[lluu]|(0.57)\hole_(0.45){d(\mathbb{C}P^1, Z_0)=2}  && &
}
\end{eqnarray}
\end{minipage}
\hspace{-1.5cm}
\begin{minipage}{0.35\hsize}
Here $\text{pt}$ is the space over $BS^1$ with the trivial map from the one point to a base point and $d(X,Y)$ stands for the cohomology interleaving distance $d_{\rm{CohI},\Q}(X, Y)$ between spaces $X$ and $Y$; see Remark \ref{rem:no_map} and
Proposition \ref{prop:M_CP}.
\end{minipage}

\medskip

\todo{The paragraph that begins with the sentence `The distances ...' is added.}
The distances are {\it realized} with interleavings. More precisely,
thanks to Proposition \ref{prop:D_E}, we see that for example, the singular cochain complexes $C^*(\mathbb{C}P^1; \Q)$ and $C^*(Z_1; \Q)$ are $\frac{7}{2}$-homotopy interleaved in the category of persistence dg modules via an appropriate functor $\alpha$; see Subection \ref{sect:ID_up_to_H} for interleavings up to homotopy and the paragraph after the diagram (\ref{eq:diagrams-diagram}) for the functor. 

We anticipate that the study of the cohomology interleaving of spaces developed in this article will bring new insights into persistence theory and topological comparison between spaces as the Gromov--Hausdorff distance is used in the study of Riemannian manifolds.
Indeed, when we compare two spaces $v_X : X\to BS^1$ and $v_Y: Y\to BS^1$, it is natural to rely on an appropriate morphism between the spaces over $BS^1$, namely a continuous map $f : X\to Y$ with $v_Y\circ f = v_X$.
However, we may compare spaces over $BS^1$ with the cohomology interleaving distance, which is a numerical two-variable homotopy invariant, even if there is no such morphism of the spaces; see
Remark \ref{two-variable homotopy invariants}.

An outline of the manuscript is as follows. Section \ref{sect:D} recalls the interleaving distance of persistence objects and the bottleneck distance of persistence vector spaces.
In Subection \ref{sect:ID_up_to_H}, the homotopy interleaving distances $d_{\text{HC}}$, $d_{\text{IHC}}$ and $d_{\text{HI}}$ are defined.
In Section \ref{sect:ID_up_to_Homotopy}, we show the formality of a persistence differential graded module over $(\mathbb{Z}, \leq)$; see Lemma \ref{lem:formality}. This fact allows us to prove Theorem \ref{thm:Equalities}. Section \ref{sect:ID_of_dgK[u]Mod} addresses the interleaving distance of dg $\K
[u]$-modules, where $\deg u = 2$. Moreover, we prove Theorem \ref{thm:IHC=CohI} mentioned above.
We also consider the bigraded $\K[t]$-module $E^{*,*}$ associated with a filtered $\K[t]$-module $H^*$, where $\deg t =1$.
As a consequence, Lemma \ref{lem:bigradedModules} enables us to recover the $\K[t]$-module structure of $H^*$ from that of $E^{*,*}$ with no extension problem.

In Section \ref{sect:BS^1}, by applying the results described in the previous sections, we consider the cohomology interleaving distances
for three different classes 
consisting of spaces over $BS^1$. Example \ref{ex:ThreeClasses} mentioned above indeed asserts that the distance $d_{\rm{CohI},\K}$ between spaces, which belong to the different classes, is infinite.

Appendix \ref{sect:toy_examples} explains explicit calculations of the cohomology interleaving distances of the complex projective spaces and the orbit spaces $Z_0$ and $Z_1$ in (\ref{eq:triangle}). Moreover,
Propositions \ref{prop:M_pt} and \ref{prop:M_CP1} are helpful in computing the cohomology interleaving distance between a persistence dg module with a small barcode and a general one.

 Appendix \ref{sect:RHT} deals with some rational homotopy invariants, whereby we observe a difference between spaces having the positive cohomology interleaving distance. In particular, the rational toral ranks and the cup-lengths of the orbit spaces in Appendix \ref{sect:toy_examples} are considered.
While as mentioned above, the cup-length is related to the cohomology interleaving distance, a relationship between the rational toral rank and the distance is currently unclear.

\subsection{Future work and perspective}\label{sect:perspective}
As mentioned above,  
Blumberg and Lesnick  \cite{Bl-L} have introduced the homotopy interleaving distance $d_{\text{HI}}$ for $\R$-spaces, namely objects in $\mathsf{Top}^{(\R, \leq)}$; see also \cite{LS}.  In particular, the distance satisfies the condition $d_{\text{HI}}(X, Y) =0$ whenever $X\simeq Y$; see \cite[Theorem 1.9]{Bl-L} and Subection \ref{sect:ID_up_to_H}.

By getting used to the {\it algebraic} interleaving distances in this article, we may introduce and consider the {\it rational} homotopy interleaving distance of $\R$-spaces. To this end, we deal with the interleaving distance in $(\mathsf{CDGA^{op}})^{(\R, \leq)}$ whose objects are persistent commutative differential graded algebras (cdga's) over $\mathbb{Q}$; see \cite{HLM, Z}.
We also refer the reader to \cite{CGL} for the study of {\it tame} persistence objects valued in a more general model category.
In \cite{KNSWY},
we introduce a functor $\Theta$ from the category of maps between spaces
to that of tame persistence cdga's over $\mathbb{Q}$, which is defined by constructing the minimal relative Sullivan models for the maps; see \cite[Sections 14 and 15]{FHT} for the models. We suppose that there exist maps $f$ and $g$ with the same domain and codomain such that
 $d_{\text{CohI}, \Q}(\Theta(f), \Theta(g))  \lneqq  d_{\text{IHC}}(\Theta(f), \Theta(g))$; see Remark \ref{rem:formal_HID_CDGA}.

It is worthwhile mentioning that homotopy invariants, which are obtained by applying the singular chain complex functor and the cohomology functor to a persistence space, give rise to more fruitful structures endowed with for example the cup products and Steenrod operations in persistent theory; see \cite{Be-St, CMSZ, G-L, LMT, MSZ}.
It may be possible to investigate each space but not a persistence space with the structures via the functor $C$ introduced
in Section \ref{sect:ID_of_dgK[u]Mod}.

It is also in our interest to consider multiparameter persistence theory in for example \cite{LS, Lesnick}.
In fact, spaces over the classifying space of a higher dimensional torus give rise to such objects in the theory via the singular cochain functor.
Thus we may investigate the moment-angle complexes that appear in toric topology from the viewpoint of multiparameter persistence theory; see \cite{BLPSS} for a related issue.


\section{The interleaving distance and the bottleneck distance}\label{sect:D}

\subsection{Interleavings}\label{sect:Interleavings}
We begin by reviewing the interleaving distance introduced in \cite{CCGGO} and results in \cite{B-S} related to the distance, which we use extensively in this article.
Let $(\R, \leq)$ be the poset defined with the usual order. Considering the poset as a category, we have the functor category $\C^{(\R, \leq)}$ for a category $\C$.
For a real number $\e \geq 0$, define a functor $T_\e : (\R, \leq) \to (\R, \leq)$ by $T_\e(a) = a+ \e$.  Moreover, the $\e$-{\it shift functor}
\[
( \ )^\e :  \C^{(\R, \leq)} \to \C^{(\R, \leq)}
\]
is defined by $( \ )^\e (F) = F^\e := FT_\e$.

\begin{defn}\label{defn:interleaving}
(\cite[Definition 4.2]{CCGGO}, \cite[Definition 3.1]{B-S})
Objects $F$ and $G$ in $\C^{(\R, \leq)}$ are {\it $\e$-interleaved} if there exist natural transformations $\varphi : F \to GT_\e$ and
$\psi : G \to FT_\e$
such that the diagram
\begin{eqnarray}\label{eq:interleaving_1}
\xymatrix@C35pt@R20pt{
F \ar[r]  \ar[dr]^(0.7){\varphi} &  F^{\e}  \ar[r] \ar[rd]^(0.7){\varphi^\e}|\hole & F ^{2\e}   \\
G  \ar[r]  \ar[ur]^(0.3){\psi}|\hole& G^\e \ar[ur]^(0.3){\psi^\e} \ar[r]  & G^{2\e}
}
\end{eqnarray}
is commutative, where horizontal arrows are the natural transformations defined by the structure maps of $F$ and $G$.
Such a tuple $(F, G, \varphi, \psi)$ is called an {\it $\e$-interleaving}.
\end{defn}

\begin{rem}\label{rem:Comm_diagrams}
The commutativity of the triangles in Definition \ref{defn:interleaving} implies the diagrams
\begin{eqnarray}\label{eq:interleaving}
\xymatrix@C10pt@R15pt{
&&F(i) \ar[rr]^-{F(i \to i+2\e)} \ar[dr]_-{\varphi(i)} &  & F(i+2\e)   &\text{and}& & F(i +\e) \ar[rd]^{\varphi(i+\e)}& \\
&&& G(i+\e) \ar[ur]_{\psi(i+\e)} &                                               &&G(i) \ar[rr]_-{G(i \to i+2\e)} \ar[ru]^-{\psi(i) }&& G(i+2\e)
}
\end{eqnarray}
are commutative for all $i \in \R$. We note that $F$ is isomorphic to $G$ in $\C^{(\R, \leq)}$ if and only if $F$ and $G$ are $0$-interleaved.
\end{rem}

\begin{defn}
For objects $F$ and $G$ in $\C^{(\R, \leq)}$, the interleaving distance $d_{\text{I}}(F, G)$ between $F$ and $G$ is defined by
\[
d_{\text{I}}(F, G) := \text{inf}(\{ \e \geq 0 \mid \text{$F$ and $G$ are $\e$-interleaved} \}\cup \{\infty\}).
\]
\end{defn}

\begin{rem}\label{rem:sum}
Let $\C$ be an additive category and $I$ a set. Suppose that for $i \in I$ $(F_i, G_i, \varphi_i, \psi_i)$ is an $\e$-interleaving in $\C^{(\R, \leq)}$. Then, we see that $\oplus_{i\in I}F_i$ and $\oplus_{i\in I}G_i$ are $\e$-interleaved with $\oplus_{i\in I}\varphi_i$ and $\oplus_{i\in I}\psi_i$.
\end{rem}

\begin{thm}\cite[Theorem 3.3]{B-S}
The function $d_{\text{\em I}}$ defined above is an extended pseudometric on 
the class of objects in $\C^{(\R, \leq)}$.
\end{thm}



\todo{[RR1, 2.4 Comment 3][RR2] Remark 2.6 is revised as Remark \ref{rem:revision_of_rem2.6} below.}
\subsection{Equivalence of interleaving and bottleneck distances}
In what follows, let $\K$ be a field of arbitrary characteristic and $\mathsf{Mod}_\K$ the category of vector spaces over $\K$ unless otherwise specified. We refer to an object in $\mathsf{Mod}_\K^{(\R, \leq)}$ as a {\it persistence vector space}.


\begin{ex}\label{ex:infinity} Let $F$ and $G$ be persistence vector spaces.
Suppose that there exists a real number $\delta$ such that $F(j)= 0$ for $j > \delta$. Moreover, we assume that
there exist $i \in \R$ and an element $x \in G(i)$ such that $G(i \to i+\delta')(x) \neq 0$ for every $\delta' >0$. Then, it follows that $d_{\text{I}}(F, G) = \infty$. In fact, suppose that
$F$ and $G$ are $\e$-interleaved.
We choose a positive number $\e'$ so that $\e' \geq \e$ and $i +\e'  >  \delta$.  By virtue of \cite[Lemma 3.4]{B-S}, we see that
$F$ and $G$ are $\e'$-interleaved. Then, the commutativity of the right-hand side diagram in Remark \ref{rem:Comm_diagrams} enables us to deduce that
$G(i \to i+2\e')(x) = 0$, which is a contradiction.
\end{ex}

Let $J$ be an interval, namely, a subset of $\R$ which satisfies the condition that if $r < s <t $ with $r, t \in J$, then $s \in J$.
We define an object $\chi_J$ in $\mathsf{Mod}_\K^{(\R, \leq)}$, called an {\it interval module},  by
\[
\chi_J(x) =
\left\{
\begin{array}{ll}
\K & \text{if $x \in J$,} \\
0 & \text{otherwise,}
\end{array}
\right.
 \ \ \ \ \ \chi_J(x\leq y) =
\left\{
\begin{array}{ll}
id_\K & \text{if $x, y \in J$,} \\
0 & \text{otherwise.}
\end{array}
\right.
\]
We call a persistence vector space $V$ {\it locally finite} ({\it pointwise finite dimensional}, p.f.d for short) if $\dim V(t) < \infty$ for $t\in \R$.
Observe that a locally finite persistence module can be decomposed uniquely as  a direct sum of interval modules; see \cite{B-C, CB, Z-C}.
Then, the {\it barcode} $\mathcal{B}_V$ \label{barcode}{\it associated with} a locally finite persistence vector space $V \cong \oplus_{J} \chi_{J}$ is defined by the multiset consisting of intervals  $J$.


\begin{lem}\label{lem:finiteIntervals_d}\text{\em (\cite[Proposition 4.12 (2)(3), Proposition 4.13(3)]{B-S})}
Let $J$ and $J'$ be finite intervals.
\begin{itemize}
\item[\rm{(1)}]
If $J' = \emptyset$ and $J$ has endpoints $a$ and $b$, then $d_{\rm{I}}(\chi_J, \chi_{J'})=\frac{b-a}{2}$.
\item[\rm{(2)}]
If $J$ and $J'$ have endpoints $a$, $b$ and $a'$, $b'$, respectively, then
\[
d_{\rm{I}}(\chi_J, \chi_{J'})=\text{\em min}\left\{\text{\em max}\{|a-a'|, |b-b'|\}, \text{\em max}\left\{\frac{b-a}{2}, \frac{b'-a'}{2}\right\} \right\}.
\]
\item[\rm{(3)}] If $\text{\em sup}(I) =\infty = \text{\em sup}(I')$ and $I$ and $I'$ have left end points $a$ and $a'$, then
$d_{\rm{I}}(\chi_I, \chi_{I'})=|a-a'|$.
\end{itemize}

\end{lem}

For multisets $A$ and $B$, define the multiset $A_B$ by $A_B:= A \amalg (\coprod_{|B|}\{\emptyset\})$. We write $f : A \leftrightarrow B$ for a bijection
$f : A_B \to B_A$.

\begin{defn}\label{defn:bottleneck_D} Let $S$ and $T$ be two barcodes. Define the bottleneck distance between $S$ and $T$ by
\[
d_{\rm{B}}(S, T) := \inf_{f  : S \leftrightarrow T}  \sup_{I \in \text{dom} (f)}d_{\text{I}}(\chi_I, \chi_{f(I)}),
\]
where $\chi_{\R}$ and $\chi_{\emptyset}$ denote the constant functors $\K$ and $0$, respectively.
\end{defn}

\begin{thm}\label{thm:I_and_B} \text{\em (The isometry theorem)} \text{\em (\cite[Theorem 4.16]{B-S}, \cite[Theorem 4.11]{CDGO})} For locally finite persistence vector spaces $V$ and $V'$, one has
\[
d_{\rm{I}}(V, V') = d_{\rm{B}}(\mathcal{B}_V, \mathcal{B}_{V'}).
\]
\end{thm}

Observe that the bottleneck distance with the $l^\infty$-metric introduced in \cite{CDGO} coincides with that in Definition \ref{defn:bottleneck_D}; see \cite[Section 4.3]{B-S} for details.

Let $\K[t]$ be the polynomial algebra generated by an element $t$ with degree $1$.
For a graded $\K[t]$-module
\begin{eqnarray}\label{eq:decomposition}
K =\bigoplus_{i=1}^n \Sigma^{-a_i}\K[t] \oplus \bigoplus_{j=1}^{n'}  \Sigma^{-b_j}(\K[t]/(t^{c_j})),
\end{eqnarray}
we define the {\it barcode} $\mathcal{B}_K$ {\it associated with} $K$ by the multiset consisting of intervals $[a_i, \infty)$ and  $[b_j, b_j + c_j)$.
Here, $\Sigma^l$ stands for the shift operator with degree $+l$; that is, $(\Sigma^lA)^i = A^{i+l}$.
We also deal with the case where $n$ and $n'$ are infinite.
The result \cite[Theorem 1]{W} implies that a bounded below graded $\K[t]$-module decomposes uniquely into the form such as (\ref{eq:decomposition})\footnote{The uniqueness of the decomposition follows from the Krull--Remak--Schmidt--Azumaya theorem.}.
Then, the barcode $\mathcal{B}_K$ associated with a graded $\K[t]$-module $K$ gives rise to the object $\chi(\mathcal{B}_K)$
in $\mathsf{Mod}_\K^{(\R, \leq)}$ defined by $\oplus_{J \in \mathcal{B}_K} \chi_{J}$. Observe that $\mathcal{B}_{\chi(\mathcal{B}_K)}=\mathcal{B}_K$ and then $d_{\rm{I}}(\chi(\mathcal{B}_K), \chi(\mathcal{B}_{K'})) = d_{\rm{B}}(\mathcal{B}_K, \mathcal{B}_{K'})$ by Theorem
\ref{thm:I_and_B}.

We conclude this section by recalling interleaving distances up to homotopy introduced in \cite{Bl-L, LS}.

\subsection{Interleavings up to homotopy}\label{sect:ID_up_to_H}

Let $\mathcal{M}$ be a cofibrantly generated model category and
$\mathcal{M}^{(\R, \leq)}$ the model category endowed with the {\it projective model structure} whose weak equivalences and fibrations are pointwise weak equivalences and pointwise fibrations of $\mathcal{M}$, respectively; see \cite[Theorem 11.6.1]{H}.

\medskip
\noindent
(1)
For objects $X$ and $Y$ in $\mathcal{M}^{(\R, \leq)}$, we say that $X$ and $Y$ are $\e$-{\it homotopy interleaved} if there exist $X\simeq X'$ and $Y\simeq Y'$ such that $X'$ and $Y'$ are $\e$-interleaved in $\mathcal{M}^{(\R, \leq)}$; see \cite[Section 3.3]{Bl-L}. Here
$W\simeq W'$ means that there is a zigzag of weak equivalences connecting $W$ and $W'$.

\noindent
(2) We say that objects $X$ and $Y$ in $\mathcal{M}^{(\R, \leq)}$ are $\e$-{\it interleaved in the homotopy category} if they are
$\e$-interleaved in $\text{Ho}(\mathcal{M}^{(\R, \leq)})$.
Observe that the shift functor $( \ )^\e : \mathcal{M}^{(\R, \leq)} \to \mathcal{M}^{(\R, \leq)}$ preserves weak equivalences. Thus,
we can consider the commutative diagram (\ref{eq:interleaving_1}) in $\text{Ho}(\mathcal{M}^{(\R, \leq)})$.

\noindent
(3) Let $q_* : \mathcal{M}^{(\R, \leq)} \to \text{Ho}(\mathcal{M})^{(\R, \leq)}$ be the functor induced by the localization functor $q : \mathcal{M} \to \text{Ho}(\mathcal{M})$.
We say that $X$ and $Y$ in $\mathcal{M}^{(\R, \leq)}$ are $\e$-{\it homotopy commutative interleaved} if $q_*X$ and $q_*Y$ are
$\e$-interleaved in $\text{Ho}(\mathcal{M})^{(\R, \leq)}$.

\medskip
Let  $X$ and $Y$ be objects in $\mathcal{M}^{(\R, \leq)}$. Blumberg and Lesnick \cite{Bl-L} introduce
 the {\it homotopy interleaving distance} and the
{\it  homotopy commutative interleaving distance} defined by
\[
d_{\text{HI}}(X, Y) \! := \text{inf}(\{ \e \geq 0 \mid \text{$X$, $Y$ are $\e$-homotopy interleaved} \}\cup \{\infty\} )\ \text{and}
\]
\[
d_{\text{HC}}(X, Y) \! := \text{inf}(\{ \e \geq 0 \mid \text{$X$, \!$Y$ are $\e$-homotopy commutative interleaved} \}\cup \{\infty\}),
\]
respectively. Moreover, Lanari and Scoccola \cite{LS} introduce
the {\it interleaving distance in the homotopy category}\label{index:homotopyID} define by
\[
d_{\text{IHC}}(X, Y) \! :=\! \text{inf}(\{ \e \geq 0 \mid \text{$X$, \!$Y$ are $\e$-interleaved in the homotopy category} \}\cup \{\infty\}).
\]

It is worthwhile mentioning that Berkouk has defined the distance in the derived category of multi-parameter persistence modules; see \cite[Section 3]{B}  for more detail.

We exhibit relationships among the three distances. We first observe that
the homotopy interleaving distance is also extended pseudometric on the class of objects in $\mathcal{M}^{(\R, \leq)}$; see \cite[Section 4]{Bl-L}
and \cite[Proposition 2.3]{LS}. By applying the universal property of the homotopy category of $\text{Ho}(\mathcal{M}^{(\R, \leq)})$ to the functor $q_*$ mentioned above,
we have a functor
$\theta : \text{Ho}(\mathcal{M}^{(\R, \leq)}) \to \text{Ho}(\mathcal{M})^{(\R, \leq)}$; see the diagram (\ref{eq:diagrams-diagram}) in Section \ref{sect:ID_of_dgK[u]Mod}. 
Thus, we establish the inequalities (\ref{eq:HI}).


We see that $d_{\text{HI}}$ on $\mathcal{M}^{(\R, \leq)}$ is characterized with the property that it is the largest distance which is bounded above by the interleaving distance and is invariant under weak equivalences. 
More precisely, we have

\begin{prop}\label{prop:Large D}
Suppose that a distance $d'$ satisfies
$d'(X, Y) \leq d_\text{\em I}(X, Y)$ for $X, Y \in \mathcal{M}^{(\R, \leq)}$ and $d'(X, Y) = d'(X', Y')$
whenever $X\simeq X'$ and $Y\simeq Y'$. Then $d' \leq d_{\text{HI}}$.  Moreover, the distance $d_{\text{\em HI}}$ itself satisfies the condition of $d'$ above.
 \end{prop}

 \begin{proof}
We assume that $d_{\text{HI}}(X, Y) < d'(X, Y)$ for some $X$ and $Y$. Then, there exists a positive real number $\e$ with $d_{\text{HI}}(X, Y) < \e < d'(X, Y)=:\delta$ such that $X$ and $Y$ are $\e$-homotopy interleaved.  Then, we have objects  $X'$ and $Y'$ with $X\simeq X'$ and $Y\simeq Y'$ which are $\e$-interleaved. Thus, we see that $\delta = d'(X, Y) = d'(X', Y') \leq d_\text{I}(X', Y') \leq \e$, which is a contradiction.  The latter half of the assertion follows from the definition of $d_{\text{HI}}$.
\end{proof}

In the rest of this section, we consider more general categories endowed with distances and functors between them. By applying Lemma \ref{lem:An_isometry} below, we prove Theorem \ref{thm:Equalities}.

Let $H: \C \to \mathcal{A}$ be a functor. Suppose that $\C$ and $\A$ are equipped with distances $d_{\C}$ and $d_{\A}$, respectively and $H$ is 1-Lipschitz; that is, $d_{\A}(H(c), H(c')) \leq d_{\C}(c, c')$ for objects $c$ and $c'$ in $\C$.  We assume that a distance is zero for isomorphic objects.
We call a morphism $f :  c \to c'$ an $H$-{\it quasi-isomorphism} if $H(f): H(c) \to H(c')$ is an isomorphism. Moreover,
we say that $x$ and $y$ in $\C$ are {\it weakly $H$-equivalent},    \todo{[RR1, 2.4 Comment 4] The space has been put.}
denoted $x \simeq y$, if there is a zig-zag of $H$-quasi-isomorphisms connecting $x$ and $y$.
Let $i : \mathcal{A} \to \C$ be a 1-Lipschitz functor.
We call an object $c$ to be $H$-{\it formal} if $c\simeq iH(c)$.

Let $d_{\C, H}$ be    \todo{[RR1, 2.4 Comment 5] Suppose ... is revised. }
the largest distance on $\C$ which is bounded above by $d_\C$ and is invariant under weak equivalences. Let $d_\C^\A$  be the distance on $\C$ defined by $d_\C^\A(c, c')=d_{\A}(H(c), H(c'))$. Since
$d_\C^\A$ is bounded above by $d_\C$ and is invariant under weak equivalences, it follows that $d_\C^\A \leq d_{\C, H}$.

\begin{defn}\label{defn:IHD_formalness}
An object $c$ in $\C$ is $d_{\C, H}$-{\it formal} if $d_{\C, H}(c, iH(c))=0$.
\end{defn}

\begin{rem}\label{rem:formal_v.s._HID}
It is readily seen that the $H$-formality implies the $d_{\C, H}$-formality.
\end{rem}

\begin{lem}\label{lem:An_isometry}
If  $c$ and $c'$ are $d_{\C, H}$-formal, then
$d_{\C, H}(c, c') = d_{\A}(H(c), H(c'))$. Thus,  $H : (\C, d_{\C, H}) \to (\A, d_{\A})$ is an isometry provided every object in $\C$ is  $d_{\C, H}$-formal.
\end{lem}

\begin{proof} We have
\begin{align}
d_\C^\A(c, c') &\leq d_{\C, H}(c, c') \leq d_{\C, H}(c, iH(c)) + d_{\C, H}(iH(c), iH(c')) + d_{\C, H}(iH(c'), c') \\
& =  d_{\C, H}(iH(c), iH(c')) \leq d_{\C} (iH(c), iH(c')) \\
& \leq d_{\A}(H(c), H(c')) = d_\C^\A(c, c').
\end{align}
Here, the third inequality follows from the boundedness of $d_{\C, H}$. The Lipschitz condition of $i$ yields the fourth inequality.
\end{proof}

\section{Interleavings up to homotopy between persistence dg modules}\label{sect:ID_up_to_Homotopy}

The aim of this section is to prove Theorem \ref{thm:Equalities},
namely the equalities
$d_{\text{HC}} = d_{\text{IHC}} = d_{\text{HI}}=d_{\text{CohI}}$.

\subsection{The equalities of distances} \todo{`equalities'  is used.}
Let $\mathsf{Ch}_\K$ be the category of cochain complexes whose objects are not necessarily bounded.
We may call such a cochain complex a {\it differential graded (dg) module}. 
We observe that $\mathsf{Ch}_\K$ can be equipped with the structure of cofibrantly generated model category, in which the weak equivalences are the maps that induce isomorphisms on cohomology groups, the so-called {\it quasi-isomorphisms}. 

Let $P$  be a poset. We view $P$ as a category with
the unique arrow $i \to j$ if $i \leq j$.
Then, the functor category $\mathsf{Ch}_\K^P$ is the model category endowed with the projective model structure; see \cite[Theorem 3.3]{B-M-R} and  \cite[Theorem 11.6.1]{H}. Thus, the three distances $d_{\text{HC}}$, $d_{\text{IHC}}$ and  $d_{\text{HI}}$ are defined on the class of the objects in
$\mathsf{Ch}_\K^{(\R, \leq)}$. We may call an object in $\mathsf{Ch}_\K^{(\R, \leq)}$ a {\it persistence dg module}.




For each integer $k$, we define a functor $H^k : \text{Ho}( \mathsf{Ch}_\K^{(\R, \leq)}) \to  \mathsf{Mod}_\K^{(\R, \leq)}$ by taking the $k$th cohomology.
\todo{[RR2] The first sentence in Remark 3.2 in the previous version is moved to this part.}

\begin{defn}\label{defn:H_HC} \label{index:CohI} The {\it cohomology interleaving distance} $d_{\text{CohI}}(X,Y)$ of persistence dg modules $X$ and $Y$ is defined by
\[
d_{\text{CohI}}(X,Y):= \text{sup}\{d_\text{I}(H^k(X), H^k(Y)) \mid k \in \mathbb{Z}\}.
\]
\end{defn}

\begin{rem}\label{rem:H_HC}
It is readily seen from \cite[Proposition 3.6]{B-S} that $d_{\text{CohI}}(X,Y) \leq d_{\text{HC}}(X, Y)$.
\end{rem}

The main theorem in this section is as follows.

\begin{thm}\label{thm:Equalities}
One has the equalities
$
d_{\text{\em HC}} = d_{\text{\em IHC}} = d_{\text{\em HI}}=d_{\text{\em CohI}}
$
on the class of the objects in $\mathsf{Ch}_\K^{(\R, \leq)}$. 
\end{thm}

Let $\eta^k :  \mathsf{grMod}_\K^{(\R, \leq)} \to  \mathsf{Mod}_\K^{(\R, \leq)}$\label{index:eta} be the functor defined by $(\eta^k)(V)(i) = V(i)^k$ for each integer $k$ and
$(H)_* :   \mathsf{Ch}_\K^{(\R, \leq)} \to \mathsf{grMod}_\K^{(\R, \leq)}$ the functor induced by  the homology functor, where
$\mathsf{grMod}_\K^{(\R, \leq)}$ denotes the category of persistence graded vector spaces.
Observe that the composite $\eta^k(H)_*$ is nothing but the functor $H^k$ mentioned in Remark \ref{rem:H_HC}. We may write $H( \ )$ for the functor $(H)_*( \ )$. 

\begin{rem}\label{rem:dHI-formality}
Let $\C$ be the category $\mathsf{Ch}_\K^{(\R, \leq)}$ endowed with the interleaving distance $d_\text{I}$ 
 and $\A$ the category $\mathsf{grMod}_\K^{(\R, \leq)}$.
 We define a distance $d_\A$ on $\A$ by
\[
d_\A(a, a') :=\text{sup}\{d_{\text{I}}(a^k, (a')^k) \mid k\in \mathbb{Z}\}.
\]
Let $H : \C \to \A$ be the functor mentioned above.
Since \[
d_\text{I}(X, Y)\geq d_\text{I}(\eta^k(H)_*(X), \eta^k(H)_*(Y))
\] for each $k \in \mathbb{Z}$, it follows that $H$ is $1$-Lipschitz. A persistence graded module is regarded as a persistence dg module with the trivial differential. Then, it is readily seen that the inclusion functor $i : \A \to \C$ is $1$-Lipschitz.

Moreover, Proposition \ref{prop:Large D} enables us to conclude that $d_{\text{HI}}$ is nothing but the distance $d_{\C, H}$ in Lemma \ref{lem:An_isometry} in the case we here consider. Furthermore, we see that $X$ is $H$-formal if and only if $X$ and $iH(X)=(H(X), 0)$ are $0$-homotopy interleaved. This follows from the definition of the $H$-formality.
\end{rem}

We prove Theorem \ref{thm:Equalities} by applying the following proposition.
The assertion and its proof are inspired by those of \cite[Theorem A]{LS}; see the paragraph after (\ref{eq:HI}) in the Introduction.

\begin{prop}\label{prop:Equality_R}
Let $X$ and $Y$ be objects 
in $\mathsf{Ch}_\K^{(\R, \leq)}$. 
If $(H(X), 0)$ and $(H(Y), 0)$ are $\delta$-interleaved,
then $X$ and $Y$ are $\delta'$-homotopy interleaved for each $\delta'$ greater than $\delta$. In particular, every persistence
dg-module is $d_{\text{\em HI}}$-formal. 
\end{prop}


In order to prove Proposition \ref{prop:Equality_R}, we discretize persistence dg modules and apply the following proposition whose proof is postponed to Subsection \ref{subsection:proof_prop:Equality}.
Here, discretization is necessary since this proposition does not hold in \(\mathsf{Ch}_\K^{(\R, \leq)}\); see Subsection \ref{sect:a_non-formal_ex}. 
We may make the same setup as in Remark \ref{rem:dHI-formality} by replacing $(\R, \leq)$ with the poset $(\mathbb{Z}, \leq)$.

\begin{prop}\label{prop:Equality}
Let $X$ and $Y$ be objects
in $\mathsf{Ch}_\K^{(\Z, \leq)}$.
If $(H(X), 0)$ and $(H(Y), 0)$ are $m$-interleaved in $\mathsf{Ch}_\K^{(\Z, \leq)}$,
then $X$ and $Y$ are $m$-homotopy interleaved in $\mathsf{Ch}_\K^{(\Z, \leq)}$. In particular, every object in $\mathsf{Ch}_\K^{(\Z, \leq)}$ is $H$-formal.
\end{prop}

\begin{proof}[Proof of Proposition \ref{prop:Equality_R}.]
Suppose that $\delta \neq 0$.
We recall the self-functor $(M_t)$ on $(\R, \leq)$ for each positive number $t \in \R$
defined by $M_t(r) = t \cdot r$
; see \cite[Section 3]{LS}.
Let $m$ be a positive integer. Then, by \cite[Lemma 3.1]{LS}, we see that
$(M_{\delta/m})^*H(X)$ and $(M_{\delta/m})^*H(Y)$ are $m$-interleaved and hence $\iota^*(M_{\delta/m})^*H(X)$ and
$\iota^*(M_{\delta/m})^*H(Y)$ are $m$-interleaved in $\mathsf{Ch}_\K^{(\Z, \leq)}$,
where
$\iota^* : \mathsf{Ch}_\K^{(\R, \leq)} \to \mathsf{Ch}_\K^{(\Z, \leq)}$ is the functor induced by the inclusion $\iota : \Z \to \R$.
The homology functor is compatible with the functor
$\iota^*(M_{\delta/m})^*$.
Thus, Proposition \ref{prop:Equality} enables us to deduce that
$\iota^*(M_{\delta/m})^*X$ and $\iota^*(M_{\delta/m})^*Y$ are $m$-homotopy interleaved in $\mathsf{Ch}_\K^{(\Z, \leq)}$.

It follows from \cite[Lemma 3.2]{LS} that $(M_{\delta/m})^*X$ and $(M_{\delta/m})^*Y$ are $(m+2)$-homotopy interleaved in $\mathsf{Ch}_\K^{(\R, \leq)}$. Therefore, there exist objects $X'$ and $Y'$ with $(M_{\delta/m})^*X\simeq X'$ and $(M_{\delta/m})^*Y\simeq Y'$ such that
$X'$ and $Y'$ are $(m+2)$-interleaved in $\mathsf{Ch}_\K^{(\R, \leq)}$. Then, it follows that $(M_{m/\delta})^*X'$ and $(M_{m/\delta})^*Y'$ are
$(\delta +2(\delta/m))$-interleaved. We observe that $X\simeq (M_{m/\delta})^*X'$ and
$Y\simeq (M_{m/\delta})^*Y'$. It turns out that $X$ and $Y$ are $(\delta +2(\delta/m))$-homotopy interleaved in $\mathsf{Ch}_\K^{(\R, \leq)}$.

We consider the case where $\delta = 0$. For a positive real number $\delta'$, choose $\delta''$ with $0 < \delta'' < \delta'$. By applying the same argument as above to  $\delta''$, we have the result.

Since $(H(X), 0)$ and $(H(H(X)), 0)$ are $0$-interleaved, we have the latter half of the assertion.
\end{proof}


\begin{proof}[Proof of Theorem  \ref{thm:Equalities}] Under the setup described in Remark \ref{rem:dHI-formality},
Proposition \ref{prop:Equality_R} enables us to deduce that the sufficient condition of Lemma \ref{lem:An_isometry} holds. Moreover, we see that the distance  $d_\C^\A$ coincides with the cohomology interleaving distance $d_{\text{CohI}}$. Observe that $d_\C^\A = d_\A \circ (H\times H)$ by definition. Thus, Lemma \ref{lem:An_isometry} yields Theorem \ref{thm:Equalities}.
\end{proof}

\begin{rem}\label{rem:formal_HID_CDGA} 
Under the same setting as in the proof of Theorem \ref{thm:Equalities}, we see that $X$ and $H(X)$ are $0$-homotopy interleaved if and only if $X$ is $H$-formal; see Subsection \ref{sect:ID_up_to_H}. Then, it follows  that $d_{\text{HI}}(X, H(X)) =0$ if $X$ is $H$-formal. Moreover, Proposition \ref{prop:Equality_R} implies that every object in $\mathsf{Ch}_\K^{(\R, \leq)}$ is $d_{\text{HI}}$-formal. Thus, we are interested in the problem whether there is a non-formal object in $\mathsf{Ch}_\K^{(\R, \leq)}$. This topic is discussed in Subsection \ref{sect:a_non-formal_ex}; see Theorem \ref{thm:nonformal_interleaved}.

Furthermore, in a more general setting as in Subsection \ref{sect:ID_up_to_H}, one might expect that there exists an non $d_{\C, H}$-formal object in a category $\C$. Lemma \ref{lem:An_isometry} may be helpful in finding such an object in the category of persistence cdga's over $\mathbb{Q}$. In fact, we have such an object if $d_{\A}(H(c), H(c')) \lneqq d_{\C, H}(c, c')$ for some $c$ and $c'$; see Section \ref{sect:perspective} and  \cite{KNSWY}.
\end{rem}


\begin{rem}
In \cite{B}, Berkouk has considered the homotopy interleaving distance in the derived category $\mathsf{D}^-$ of the category
$\mathsf{Mod}_\K^{(\R^d, \leq)}$ of multi-parameter persistence modules whose objects are bounded below.
In particular, the result \cite[Theorem 2]{B} asserts that the canonical functor $\iota$ from
$\mathsf{Mod}_\K^{(\R^d, \leq)}$ to $\mathsf{D}^-$ is object-wise isometric with respect to the interleaving distance and
the interleaving distance in the homotopy category described in Subsection \ref{sect:ID_up_to_H}, respectively.

Theorem \ref{thm:Equalities} implies that the functor $\iota$ factors through the category $(\mathsf{Ch}_\K^{(\R, \leq)}, d_{\rm CohI})$ with
object-wise isometric functors when $d = 1$ and  $\mathsf{Ch}_\K^{(\R, \leq)}$ is restricted to cochain complexes bounded above in objects.
Observe that objects in $\mathsf{D}^-$ are chain complexes but not cochain complexes. 
\end{rem}

\subsection{Proof of Proposition \ref{prop:Equality}}
\label{subsection:proof_prop:Equality}
In order to prove Proposition \ref{prop:Equality}, we regard an object $X$ in $\mathsf{Ch}_\K^{(\Z, \leq)}$ as a {\it differential bigraded (dbg) $\K[t]$-module} (a cochain complex of graded $\K[t]$-modules)
\[
\left( \bigoplus_{(i, n) \in \mathbb{Z}^2}X(i)^n, d \right)
\]
for which
$(\bigoplus_n X(i)^n, d)$ is a dg module for each $i$ and the module structure $\times t :  X(i)^n \to  X(i+1)^n$ is given by the structure map
$X(i \to i+1)$. Observe that $\times t \circ d = d \circ \times t$. 
We show that
\[\left(\bigoplus_{(i, n) \in \mathbb{Z}^2}X(i)^n, d \right) \simeq
\left(\bigoplus_{(i, n) \in \mathbb{Z}^2}H^n(X(i)^*), 0 \right)
\]
 as a dbg $\K[t]$-module. 


The following lemma is a generalization of the assertion of \cite[Remark 3.7]{AZ} to an unbounded case.

\begin{lem}\label{lem:formality}
For an object $(X, d)$ in $\mathsf{Ch}_\K^{(\Z, \leq)}$, there exist a dbg $\K[t]$-module $Q$ and quasi-isomorphisms
$X \stackrel{\simeq}{\leftarrow} Q \stackrel{\simeq}{\rightarrow} H(X)$
of dbg $\K[t]$-modules. As a consequence, $X \simeq H(X)$ in $\mathsf{Ch}_\K^{(\Z, \leq)}$.
\end{lem}

As seen in Remark \ref{rem:hereditaryAb} below, Lemma \ref{lem:formality} follows from a more general result.
We here prove the lemma by a constructive approach.
First we note and prove the following lemma for the reader's convenience.

Recall that a graded algebra \(R\) is a graded principal domain
if it is a domain and
every homogeneous ideal is generated by a single element
(i.e., every graded \(R\)-submodule of \(R\) is a free graded \(R\)-submodule of rank 0 or 1);
see \cite{Loeh} for details.

\begin{lem}
  \label{lem:graded_pid}
  Let \(R\) be a graded principal ideal domain 
  and \(F\) a free graded \(R\)-module.
  Then any graded \(R\)-submodule \(G\subset F\) is also free as a graded \(R\)-module.
\end{lem}
\begin{proof}
  The proof is based on \cite[Theorem IV 6.1]{Hun},
  but here we carefully choose a \emph{homogeneous} basis.

  Let \(\{x_i\}_{i\in I}\) be a basis of \(F\).
  Then we have an internal direct sum decomposition \(F=\bigoplus_{i\in I}R x_i\),
  where \(R x_i\) denotes the free \(R\)-module of rank 1 with the homogeneous basis \(x_i\).
  Now we choose a well ordering \(\le\) of the set \(I\) (by using the axiom of choice).
  Let \(J = I \coprod {\infty}\) be the well-ordered set defined by \(i < \infty\) for any \(i\in I\).
  For \(i\in I\), we define
  \(i+1 = \min\{j\in J\mid i < j\} \in J\),
  where the minimum exists since \(J\) is well-ordered.
  For \(j\in J\), we write
  \(F_j = \bigoplus_{i<j}R x_i\) and \(G_j = G\cap F_j\).

  Now we show that, for each \(i\in I\),
  there exists a \emph{homogeneous} element \(b_i\in G_{i+1}\setminus G_i\)
  such that \(G_{i+1} = G_i\oplus R b_i\)
  and \(R b_i\) is a free \(R\)-module of rank 0 or 1.
  Consider the short exact sequence
  \(0 \to G_i \to G_{i+1} \to G_{i+1}/G_i \to 0\).
  \todo{[RR1, 2.4 Comment 6] Here is the detail.}
  Note that
  \(G_{i+1}/G_i\) is a graded \(R\)-submodule of (a degree shift of) \(R\)
  since
  \begin{equation}
    G_{i+1}/G_i
    = G_{i+1}/(G_{i+1}\cap F_i)
    \cong (G_{i+1}+F_i)/F_i
    \subset F_{i+1} / F_i
    \cong R x_i.
  \end{equation}
  Now \(G_{i+1}/G_i\) is a free graded \(R\)-module of rank 0 or 1
  by the assumption that
  \(R\) is a graded principal ideal domain.
  Hence the short exact sequence splits as graded \(R\)-modules
  and we have an internal direct sum decomposition
  \(G_{i+1} = G_i \oplus R b_i\).
  Here, \(b_i = 0\) if the rank of \(G_{i+1}/G_i\) is 0,
  and \(R b_i \cong G_{i+1}/G_i\) is a free graded module of rank 1 if the rank of \(G_{i+1}/G_i\) is 1.

  Define \(B = \{b_i\mid i\in I, b_i\neq 0\}\).
  Then, by using transfinite induction,
  we can show that \(B\) is a homogeneous basis of \(G\);
  see \cite[Theorem IV 6.1]{Hun} for details.
\end{proof}

Now we give a proof of Lemma \ref{lem:formality}.

\begin{proof}[Proof of Lemma \ref{lem:formality}] 
Let $\{[b_\lambda(i)^n]\}$ be a set of generators of $H(X)$ as a bigraded $\K[t]$-module, where $b_\lambda(i)^n$ is in $X(i)^n$.
Observe that $H(X) = \bigoplus_{(i, n)} H^n(X(i))$.

Let $F_0$ be the free $\K[t]$-module generated by $\{b_\lambda(i)^n\}$.
Since $\K[t]$ is a graded principal ideal domain, 
it follows that the kernel \(\ker(p\circ\varphi)\) of the composite of
the natural map $\varphi : F_0 \to \text{Ker} \ d$ and
the projection $p : \text{Ker} \ d \to H(X)$ is a free graded $\K[t]$-module by Lemma \ref{lem:graded_pid}.
Let $B=\{f_\mu(i)^n\}$ be the basis of \(\ker(p\circ\varphi)\),
where  $f_\mu(i)^n$ is of bidegree $(i,n)$.

Let $F_1$ be the $\K[t]$-module obtained from \(\ker(p\circ\varphi)\) by shifting the second degree by one.
Namely, \(F_1\) is the free \(\K[t]\)-module with the basis $\{ \alpha_\mu(i)^{n} \}$,
where the indices run in the same set as in \(\{f_\mu(i)^n\}\)
and \(\alpha_\mu(i)^n\) is of bidegree \((i,n-1)\).
We define the differential $D :F_1 \to F_0$ by  $D(\alpha_\mu(i)^{n}) =  f_\mu(i)^n$.
Since $(p\circ \varphi)(f_\mu(i)^n)=0$ for each element $f_\mu(i)^n$ in the basis $B$, there exists an element $z_\mu(i)^n$ in $X$ such that
$\varphi(f_\mu(i)^n)=d(z_\mu(i)^n)$. We define a morphism
$
\psi : Q:=(F_0\oplus F_1, D) \to X
$
of dg $\K[t]$-modules by $\psi(b_\lambda(i)^n) = b_\lambda(i)^n$ and $\psi(\alpha_\mu(i)^{n}) = z_\mu(i)^n$.
Moreover, we have a quasi-isomorphism
$g: Q \stackrel{\simeq}{\to} H(X)$ defined by $g(x +y) = [x]$, where $x \in F_0$ and $y \in F_1$.
\end{proof}

\begin{rem}\label{rem:hereditaryAb}
An object $M$ in $\mathsf{Ch}_\K^{(\Z, \leq)}$ is regarded as a chain complex of graded
$\K[t]$-module and then an object in 
the derived category of graded $\K[t]$-modules. Thus, Lemma \ref{lem:formality} follows from the more general result
\cite[Proposition 4.4.15]{Krause} which asserts that every object in the derived category of an abelian category $\mathcal{A}$ is formal ({\it quasi-isomorphic to its cohomology} in the sense in \cite{Krause}) if and only if $\mathcal{A}$ is {\it hereditary}; that is,
the functor $\text{Ext}_{\mathcal{A}}^2(\text{-}, \text{-})$ vanishes.
Therefore, the result on the derived category implies that Lemma \ref{lem:formality} cannot be generalized to an assertion for multi-parameter persistence dg-modules.  In fact, the second Ext functor does not vanish in the category of graded $\K[t]^{\otimes n}$-modules for $n\geq 2$.
\end{rem}

\begin{proof}[Proof of Proposition \ref{prop:Equality}]  Let $X$ and $Y$ be objects in $\mathsf{Ch}_\K^{(\Z, \leq)}$.
Then, the assumption and Lemma \ref{lem:formality} imply that $X$ and $Y$ are $m$-homotopy interleaved.
\end{proof}

\subsection{A non-formal example in \(\mathsf{Ch}_\K^{(\R, \leq)}\)}\label{sect:a_non-formal_ex}
\newcommand{\nf}{P}
\newcommand{\fm}{Q}
\newcommand{\fmd}{\overline{Q}}
\newcommand{\bo}{\alpha}
\newcommand{\bz}{\beta}
\newcommand{\ktr}{\K[t^{\R_+}]}
In this subsection, we give a persistence dg module
\(X\in \mathsf{Ch}_\K^{(\R, \leq)}\)
which illustrates differences
between \(H\)-formality and \(d_{\text{HI}}\)-formality,
and between \(\mathsf{Ch}_\K^{(\Z, \leq)}\) and \(\mathsf{Ch}_\K^{(\R, \leq)}\).
More precisely, we prove

\begin{thm}
  \label{thm:nonformal_interleaved}
  Let \(\nf\in\mathsf{Ch}_\K^{(\R, \leq)}\) be the persistence dg module
  defined in Definition \ref{defn:nonformal}.
  Then we have
  \begin{enumerate}
    \item \label{item:nonformal}%
      \((\nf, d)\) is not
      \(H\)-formal; that is,
      \((\nf, d)\) and \((H(\nf), 0)\) are not \(0\)-homotopy interleaved.
    \item \label{item:interleaved}
    \((\nf, d)\) is \(d_{\text{\em HI}}\)-formal; that is,
    \((\nf, d)\) and \((H(\nf), 0)\) are \(\varepsilon\)-homotopy interleaved
    for any \(\varepsilon>0\).    \todo{[RR1, 2.4 Comment 7] See the paragraph below and our report `Responses'. }
  \end{enumerate}
  See Remark \ref{rem:dHI-formality} for the terminology on formality.
\end{thm}

Hence (\ref{item:nonformal}) implies that \(\nf\) is a counterexample for
Proposition \ref{prop:Equality} and Lemma \ref{lem:formality}
in \(\mathsf{Ch}_\K^{(\R, \leq)}\),
and shows that Proposition \ref{prop:Equality_R} does not hold when \(\delta = \delta'\).
Although (\ref{item:interleaved}) follows from Proposition \ref{prop:Equality_R},
here we give an \textit{explicit} proof
by giving an \(\varepsilon\)-homotopy interleaving
(an \(\varepsilon\)-interleaving together with quasi-isomorphisms)
to illustrate the situation cleanly. 
We give constructive proofs for lemmas below
although some of them follow from general theorems.
\todo{[RR1, 2.4 Comment 9] This sentence explains why we give constructive proof for Lemma \ref{lem:lifting}}

Similarly to the case of \(\mathsf{Ch}_\K^{(\Z, \leq)}\),
we regard a persistence dg module \(X\in\mathsf{Ch}_\K^{(\R, \leq)}\)     \todo{[RR1, 2.4 Comment  8] This concerns Notation. See our report `Responses'.}
as a complex of \(\R\)-graded \(\ktr\)-module
\((\bigoplus_{(i,n)\in\R\times\Z}X(i)^n, d)\).
Following the notation of \cite[Subsection 5.2]{HLM},
here \(\R_+\) denotes the set of non-negative real numbers and
\(\ktr\) denotes the \(\R\)-graded \(\K\)-algebra
with the basis \(\{t^i\mid i\in \R_+\}\).
By an abuse of notation, we call it a \textit{dbg \(\ktr\)-module}.

\begin{defn}
  \label{defn:nonformal}
 A dbg \(\ktr\)-module \(\nf\) is defined by
  \begin{itemize}
    \item \(\nf^1 = \ktr\bo\), i.e.\ the free \(\ktr\)-module
      generated by \(\bo\in \nf(0)^1\)
    \item \(\nf^0 = \bigoplus_{n\in\Z_{>0}}\ktr\bz_n\), i.e.\ the free \(\ktr\)-module
      generated by \(\bz_n\in\nf(\frac{1}{n})^0\) for \(n\in\Z_{>0}\)
    \item \(d\colon \nf(i)^0 \to \nf(i)^1\) is the map of \(\ktr\)-modules
      defined by \(d(\bz_n)=t^{\frac{1}{n}}\bo\).
  \end{itemize}
  See Figure \ref{fig:nonformal}.
\end{defn}

\begin{figure}[htbp]
  \begin{tikzpicture}
    \newcommand{\len}{5}
    \newcommand{\maxN}{4}
    \newcommand{\fy}{1}
    \newcommand{\drawLine}[4][black]{%
      \draw (#2,#3) -- (\len, #3);
      \draw[fill=#1] (#2,#3) circle (2pt);
      \node[anchor=east] at (#2,#3) {#4};
    }
    \node (X1) at (-1,\fy) {\(\nf^1\)};
    \node (X0) at (-1,-0.5) {\(\nf^0\)};
    \draw[->] (X0) -- (X1) node [pos=0.5, left] {\(d\)};
    \drawLine{0}{\fy}{\(\bo\)}
    \foreach \n in {1,...,\maxN} {
      \pgfmathsetmacro\x{3/\n}
      \pgfmathsetmacro\y{0.4*(\n-\maxN)}
      \drawLine{\x}{\y}{\(\bz_\n\)}
      \draw[->] (\x,\y+0.1) -- (\x,\fy-0.05);
      \ifnum\n<3
        \node at (\x,\fy+0.3) {\(t^{\frac{1}{\n}}\bo\)};
      \fi
    }
    \node[rotate=20] at (3/\maxN-0.2,0.4+0.1) {\(\vdots\)};
    \node at (4,0.4) {\(\vdots\)};
  \end{tikzpicture}
  \caption{Diagram of \(\nf\)}
  \label{fig:nonformal}
\end{figure}

The following proposition is the key to proving Theorem \ref{thm:nonformal_interleaved} (\ref{item:nonformal}).

\begin{prop}
  \label{prop:nonformal}
  Let \(\varphi\colon (\nf, d)\to (H(\nf), 0)\) be a morphism of dbg \(\ktr\)-modules.
  Then \(\varphi\) is not a quasi-isomorphism.
\end{prop}
\begin{proof}
  Assume that \(\varphi\) is a quasi-isomorphism.
  First we show that \(\ker(\varphi(\frac{1}{n})^0) \cong \K\) for any \(n\ge 1\).
  Since \(H^0(\nf) = \ker d\) and \(\varphi\) is a quasi-isomorphism,
  the map \(\varphi(\frac{1}{n})^0\) restricts to the isomorphism
  \(\ker d \xrightarrow{\cong}H^0(\nf(\frac{1}{n}))\).
  Hence we have an internal direct sum decomposition
  \(\nf(\frac{1}{n}) = \ker d \oplus \ker(\varphi(\frac{1}{n})^0)\).
  This implies that
  \begin{equation}
    \ker(\varphi(\tfrac{1}{n})^0)
    \cong \nf(\tfrac{1}{n}) / \ker d
    \xrightarrow[\cong]{d} \im d
    = \nf(\tfrac{1}{n})^1
    \cong \K.
  \end{equation}

  Now consider the map \(\nf(\frac{1}{n+1} \to \frac{1}{n})^0\),
  which is injective by definition.
  This induces an injection
  \(\ker(\varphi(\frac{1}{n+1})^0) \to \ker(\varphi(\frac{1}{n})^0)\),
  which is also surjective since both modules are isomorphic to \(\K\).
  Hence we have
  \begin{equation}
    \ker(\varphi(\tfrac{1}{n})^0)
    \subset
    \bigcap_{m\ge n}\im(\nf(\tfrac{1}{m}\to\tfrac{1}{n}))
    = 0,
  \end{equation}
  which contradicts the above isomorphism \(\ker(\varphi(\frac{1}{n})^0) \cong \K\).
\end{proof}

Since \(\nf^n\) is a free \(\ktr\)-module for each \(n\),
it satisfies the following lifting property for quasi-isomorphisms.

\begin{lem}
  \label{lem:lifting}
  \todo{[RR1, 2.4 Comment 9] We prefer a concrete proof.}
  Let \(X\) and \(Y\) be dbg \(\ktr\)-modules and
  consider morphisms \(f\colon \nf \to Y\) and \(g\colon X\to Y\) of dbg \(\ktr\)-modules.
  Assume that \(g\) is a quasi-isomorphism.
  Then, there exists a morphism \(\tilde{f}\colon \nf\to X\) of dbg \(\ktr\)-modules
  such that \(g\circ\tilde{f}\) and \(f\) are homotopic.
  \begin{equation}
    \begin{tikzcd}
      & X \ar[d,"g","\simeq"'] \\
      \nf \ar[r,"f"] \ar[ru,"\tilde{f}",dashed] & Y
    \end{tikzcd}
  \end{equation}
\end{lem}
The above lemma essentially follows from \cite[Proposition 6.4 (ii)]{FHT},
since \(\nf\) is semifree in the sense of \cite{FHT}
if we ignore the difference between the \(\R\)- and \(\Z\)-gradings.
However, since our object \(P\) is \(\R\)-graded
whereas only the \(\Z\)-graded case is dealt with in \cite{FHT},
the result cannot be applied directly.
For the reader's convenience, we therefore provide a complete and constructive proof here.
\begin{proof}
  Since \(H(g)\colon H(X)\to H(Y)\) is surjective,
  there exists \([x]\in H^1(X(0))\) such that \(H(g)[x]=[f \bo]\); that is,
  there exist \(x\in X(0)^1\) and \(y\in Y(0)^0\)
  such that \(dx=0\) and \(dy=g x - f \bo\).

  Since \(g(t^{\frac{1}{n}}x)=d(t^{\frac{1}{n}}y+f \bz_n)\),
  we have \(H(g)[t^{\frac{1}{n}}x]=0 \in H^1(Y(\frac{1}{n}))\) for all \(n>0\).
  Hence injectivity of \(H(g)\colon H(X)\to H(Y)\) implies that,
  for each \(n>0\), there exists \(x'_n\in X(\frac{1}{n})^0\)
  such that \(dx'_n=t^{\frac{1}{n}}x\).

  Note that the above definitions imply \(d(g x'_n-f \bz_n-t^{\frac{1}{n}}y)=0\).
  Hence, by using surjectivity of \(H(g)\colon H(X)\to H(Y)\) again,
  we have \(H(g)[x''_n] = [g x'_n-f \bz_n-t^{\frac{1}{n}}y] \in H^0(Y(\frac{1}{n}))\)
  for some \([x''_n]\in H^0(X(\frac{1}{n}))\);
  that is, there are \(x''_n\in X(\frac{1}{n})^0\) and \(y'_n\in Y(\frac{1}{n})^{-1}\)
  such that \(dx''_n=0\) and \(dy'_n=(g x'_n-f \bz_n-t^{\frac{1}{n}}y)-g x''_n\).

  Now we define maps
  \(\tilde f\colon \nf\to X\) and \(h\colon \nf\to Y\)
  of bigraded \(\ktr\)-modules by
  \(\tilde f(\bo)=x\), \(\tilde f(\bz_n)=x'_n-x''_n\),
  \(h(\bo)=y\), and \(h(\bz_n)=y'_n\).
  Then, it is readily seen that \(\tilde f\) is a map of dbg \(\ktr\)-modules,
  namely,  \(\tilde f d = d \tilde f\),
  and  \(hd+dh=g\tilde f-f\); that is,
  \(h\) is a homotopy from \(g\tilde f\) to \(f\).
  This completes the proof of the lemma.
\end{proof}

Now we are ready to prove Theorem \ref{thm:nonformal_interleaved} (\ref{item:nonformal}).
\begin{proof}[Proof of Theorem \ref{thm:nonformal_interleaved} (\ref{item:nonformal})]
  Assume that \((\nf, d)\) is \(H\)-formal.
  By definition, there is a zig-zag of quasi-isomorphisms between \((\nf, d)\) and \((H(\nf), 0)\).
  By Lemma \ref{lem:lifting}, we have a quasi-isomorphism
  \(\varphi\colon (\nf, d)\xrightarrow{\simeq} (H(\nf), 0)\).
  This contradicts Proposition \ref{prop:nonformal}.
\end{proof}

Next we introduce two dbg \(\ktr\)-modules \(\fm\) and \(\fmd\)
to prove Theorem \ref{thm:nonformal_interleaved} (\ref{item:interleaved}).

\begin{lem}
  \label{lem:removed_formal}
  Let \(\fm\) be the dbg \(\ktr\)-submodule of \(\nf\) defined by \(\fm(i) = \nf(i)\) for \(i>0\) and \(\fm(i) = 0\) for \(i\le 0\);
  see Figure \ref{fig:formal}.
  Then, there is a quasi-isomorphism
  \((H(\fm), 0) \xrightarrow{\simeq} (\fm, d)\).
\end{lem}
\begin{proof}
  Since \(d\colon \fm(i)^0 \to \fm(i)^1\) is surjective
  for all \(i\in\R\) (including \(i=0\)),
  we have \(H(\fm) = H^0(\fm) = \ker(d^0\colon \fm^0\to\fm^1)\).
  Hence, the inclusion map
  \(H(\fm) = \ker d^0 \hookrightarrow \fm\)
  is a map of dbg \(\ktr\)-modules
  and is a quasi-isomorphism.
\end{proof}


\begin{figure}[htbp]
  \begin{minipage}{0.48\textwidth}
    \begin{tikzpicture}
      \newcommand{\len}{4.5}
      \newcommand{\maxN}{4}
      \newcommand{\fy}{1}
      \newcommand{\drawLine}[4][black]{%
        \draw (#2,#3) -- (\len, #3);
        \draw[fill=#1] (#2,#3) circle (2pt);
        \node[anchor=east] at (#2,#3) {#4};
      }
      \node (X1) at (-0.5,\fy) {\(\fm^1\)};
      \node (X0) at (-0.5,-0.5) {\(\fm^0\)};
      \draw[->] (X0) -- (X1) node [pos=0.5, left] {\(d\)};
      \drawLine[white]{0}{\fy}{}
      \foreach \n in {1,...,\maxN} {
        \pgfmathsetmacro\x{3/\n}
        \pgfmathsetmacro\y{0.4*(\n-\maxN)}
        \drawLine{\x}{\y}{\(\bz_\n\)}
        \draw[->] (\x,\y+0.1) -- (\x,\fy-0.05);
        \ifnum\n<3
          \node at (\x,\fy+0.3) {\(t^{\frac{1}{\n}}\bo\)};
        \fi
      }
      \node[rotate=20] at (3/\maxN-0.2,0.4+0.1) {\(\vdots\)};
      \node at (4,0.4) {\(\vdots\)};
    \end{tikzpicture}
  \end{minipage}
  \begin{minipage}{0.48\textwidth}
    \begin{tikzpicture}
      \newcommand{\len}{4.5}
      \newcommand{\maxN}{4}
      \newcommand{\fy}{1}
      \newcommand{\drawLine}[4][black]{%
        \draw (#2,#3) -- (\len, #3);
        \draw[fill=#1] (#2,#3) circle (2pt);
        \node[anchor=east] at (#2,#3) {#4};
      }
      \node (X1) at (-0.5,\fy) {\(\fmd^1\)};
      \node (X0) at (-0.5,-0.5) {\(\fmd^0\)};
      \draw[->] (X0) -- (X1) node [pos=0.5, left] {\(d\)};
      \drawLine[white]{0}{\fy}{}
      \node (f1) at (0,\fy+0.4) {};
      \fill (f1) circle (2pt);
      \node[anchor=south west] at (f1) {\(\bar\bo\)};
      \foreach \n in {1,...,\maxN} {
        \pgfmathsetmacro\x{3/\n}
        \pgfmathsetmacro\y{0.4*(\n-\maxN)}
        \drawLine{\x}{\y}{\(\bz_\n\)}
        \draw[->] (\x,\y+0.1) -- (\x,\fy-0.05);
        \ifnum\n<3
          \node at (\x,\fy+0.3) {\(t^{\frac{1}{\n}}\bo\)};
        \fi
      }
      \node[rotate=20] at (3/\maxN-0.2,0.4+0.1) {\(\vdots\)};
      \node at (4,0.4) {\(\vdots\)};
    \end{tikzpicture}
  \end{minipage}

  \caption{Diagrams of \(\fm\) and \(\fmd\)}
  \label{fig:formal}
\end{figure}

By using \(\fm\), we define a dbg \(\ktr\)-module \(\fmd \) by \(\fmd = \K \bar\bo \oplus \fm\),
where \(\K\bar\bo\) denotes a dbg \(\ktr\)-module
concentrated in \((\K\bar\bo)(0)^1 = \K\bar\bo\); see Figure \ref{fig:formal}.

\begin{lem}
\label{lem:dotted}
  \begin{enumerate}
    \item There is a quasi-isomorphism
      \((H(\fmd), 0) \xrightarrow{\simeq} (\fmd, d)\).
    \item \(H(\nf)\) is isomorphic to \(H(\fmd)\).
  \end{enumerate}
\end{lem} \begin{proof}
  By \(\fmd = \K \bar\bo \oplus \fm\),  \todo{[RR2] Lemma 3.16 can be applied here through this direct sum}
  Lemma \ref{lem:removed_formal} immediately implies that
  there is a quasi-isomorphism \((H(\fmd), 0) \xrightarrow{\simeq} (\fmd, d)\).
  Since
  \(d\colon \nf(i)^0 \to \nf(i)^1\) and
  \(d\colon \fmd(i)^0 \to \fmd(i)^1\)
  are surjective for \(i>0\),
  we have an isomorphism \(H(\nf) \cong H(\fmd)\).
\end{proof}

\begin{prop}
  \label{prop:interleaved}
  For any \(\varepsilon>0\),
  \((\nf, d)\) and \((\fmd, d)\) are \(\varepsilon\)-interleaved.
\end{prop}
\begin{proof}
  Define morphisms \(\varphi\colon \nf \to \fmd^\varepsilon\) and \(\psi\colon \fmd \to \nf^\varepsilon\)
  of dbg \(\ktr\)-modules by
  \(\varphi(\bo)=t^\varepsilon\bo\),
  \(\varphi(\bz_n)=t^\varepsilon\bz_n\),
  \(\psi(t^i\bo)=t^{i+\varepsilon}\bo\) (for \(i>0\)),
  \(\psi(\bo')=0\), and
  \(\varphi(\bz_n)=t^\varepsilon\bz_n\).
  It is readily seen that these satisfy (\ref{eq:interleaving_1}).
\end{proof}

Now we give a proof for Theorem \ref{thm:nonformal_interleaved} (\ref{item:interleaved}).

\begin{proof}[Proof of Theorem \ref{thm:nonformal_interleaved} (\ref{item:interleaved})]
  By Lemma \ref{lem:dotted},
  we have  quasi-isomorphisms
  \begin{equation}
    (\fmd, d) \xleftarrow{\simeq} (H(\fmd), 0) \cong (H(\nf), 0).
  \end{equation}
  Now Proposition \ref{prop:interleaved} completes the proof.
\end{proof}

\section{The cohomology interleaving of dg $\K[u]$-modules}\label{sect:ID_of_dgK[u]Mod}

\subsection{Interleavings of dg $\K[u]$-modules}\label{sect:I_M}
Let $\K[u]$ be the polynomial algebra generated by an element $u$ with degree $2$.
\todo{[RR1, 2.1 i)] Some functors are summarized in the paragraph after the diagram (\ref{eq:diagrams-diagram}).  We think that List of Symbols before the References is useful.}

\begin{defn} \todo{[RR1, 2.4 Comment 10] `if' is changed to `together with'.}
A {\it differential graded (dg) $\K[u]$-module}
is
a differential graded module  $\{M^l, \partial \}$ together with
a $\K$-linear map $u : M^l\to M^{l+2}$ which satisfies the condition that
$u\circ \partial = \partial \circ u$.
\end{defn}



Let $\K[u]\text{-}\mathsf{Ch}$ denote the category of dg $\K[u]$-modules.
In order to develop persistence theory for dg $\K[u]$-modules,
we assign a persistence dg module to each dg $\K[u]$-module via a functor.
For a dg $\K[u]$-module
$M=\{M^l, \partial\}$,
we define a functor $C: \K[u]\text{-}\mathsf{Ch} \to  \mathsf{Ch}_\K^{(\Z, \leq)}$ by
$C(\{M^l, \partial \})(i) = \Sigma^{2i}M$ and
\begin{eqnarray*}
C(\{M^l, \partial \})(i \to i+1) :  C(\{M^l, \partial \})(i)  \stackrel{u \cdot }{\longrightarrow} C(\{M^l, \partial \})(i+1)
\end{eqnarray*}
with the multiplication by $u$.

As it will 
be seen in Section \ref{sect:BS^1}, the functor $C$ allows us to bring topological spaces over $BS^1$ to persistence theory.

Let $\K[t]\text{-}\mathsf{grMod}$ stand for the category of graded $\K[t]$-modules. We denote by $\mathsf{grMod}_\K$ the category of graded vector spaces.
An object $K$ in $\mathsf{Mod}_\K^{(\Z, \leq)}$ gives the graded $\K[t]$-module $\gamma(K):=\oplus_{i\in \mathbb{Z}} K(i)$ with the module structure defined by $t\cdot x = K(i\to i+1)(x)$ for $x \in K(i)$. 
It is readily seen that $\gamma$ gives rise to an equivalence 
$\gamma : \mathsf{Mod}_\K^{(\Z, \leq)} \stackrel{\simeq}{\to} \K[t]\text{-}\mathsf{grMod}$  of the categories.

\begin{minipage}{0.65\hsize}     \begin{eqnarray}\label{eq:complexes}
\xymatrix@C7pt@R7pt{
  &\cdots  \Sigma^{-2}M \ar[r]^-{u \cdot}& M \ar[r]^-{u \cdot}& \Sigma^{2}M \ar[r]^-{u \cdot} &  \Sigma^{4}M &\hspace{-0.6cm} \cdots \\
&\vdots &\vdots & \vdots & &  &  \\
&M^0 \ar[u]& M^2 \ar[u]& M^4 \ar[u] & \vdots\\
&M^{-1} \ar[u]& M^1 \ar[u]&  M^3  \ar[u]& M^5  \ar[u]\\
\cdots &M^{-2}\ar[u] & M^0 \ar[u]& M^2 \ar[u]& M^4 \ar[u]& \cdots \\
&\vdots \ar@{->}[u] & M^{-1} \ar[u]& M^1 \ar[u]& M^3 \ar[u] 
}
\end{eqnarray}
\end{minipage}
\hspace{-1.2cm}
\begin{minipage}{0.35\hsize}
The diagram (\ref{eq:complexes}) explains the persistence module $C(M)$.
The vertical arrows are differentials in the dg $\K[u]$-module $M$.

\ \ As for the functor $h^k$ defined below, $h^0M$  and $h^1M$ are direct sums of the cohomology groups of the third and second rows, respectively.
\end{minipage}
\todo{[RR1, 2.4 Comment 11] $\times$ to $\cdot$.}

\medskip
We recall the functors
$
\mathsf{Ch}_\K^{(\R, \leq)} \stackrel{(H)_*}{\to} \mathsf{grMod}_\K^{(\R, \leq)} \stackrel{\eta^k}{\to} \mathsf{Mod}_\K^{(\R, \leq)}
$
 mentioned before Proposition \ref{prop:Equality}.
Moreover, we have maps
\[
(\Z, \leq) \stackrel{\iota}{\to} (\R, \leq) \stackrel{\lfloor \ \rfloor}{\to} (\Z, \leq)
\]
of posets defined by the usual inclusion $\iota$ and the floor function $\lfloor \ \rfloor$\label{index:floor}, respectively.
Observe that with $\lfloor \ \rfloor \circ \iota = 1$.
Thus, we obtain a commutative diagram consisting of categories and functors

\begin{eqnarray}\label{eq:diagrams-diagram}
\xymatrix@C7pt@R20pt{
& & \K[u]\text{-}\mathsf{Ch} \ar[r]^{C} \ar[rd]^-{h^k:=\utot{k}hq} \ar[d]_{q}&  \mathsf{Ch}_\K^{(\Z, \leq)} \ar[r]_-{(H)_*} \ar@/^2.0pc/[rr]_-{(\lfloor \ \rfloor)^*} & \mathsf{grMod}_\K^{(\Z, \leq)} \ar[d]^-{\eta^k}& \mathsf{Ch}_\K^{(\R, \leq)} \ar[d]^-{\eta^k  (H)_*} \ar[rd]^\pi \ar[r]^-{q_*}&  \text{Ho}(\mathsf{Ch}_\K)^{(\R, \leq)} \\
&&\text{D}(\K[u]) \ar[d]_{h}  \ar@{.>}@/_2.0pc/[rrrr]_-{\mu} & \K[t]\text{-}\mathsf{grMod}  \ar@{~>}[d]_-{\mathcal{B}_{(\ )}} & \mathsf{Mod}_\K^{(\Z, \leq)}  \ar[r]^{(\lfloor \ \rfloor)^*}_{\text{embedding}} \ar[l]^-{\simeq}_-\gamma &
 \mathsf{Mod}_\K^{(\R, \leq)} \ar@{~>}[d] & \text{Ho}(\mathsf{Ch}_\K^{(\R, \leq)}) \ar@{.>}[l]^{\xi^k}  \ar[u]_-{\theta}\\
&&\K[u]\text{-}\mathsf{grMod}  \ar[ru]_{\utot{k}}&  (\mathcal{BAR}, d_B) \ar[rr]_\chi & & (\mathsf{Mod}_\K^{(\R, \leq)}, d_{\text{I}}). &
}
\end{eqnarray}
Here $\text{D}(\K[u])$
\label{index:derivedcat}
denotes the derived category of dg $\K[u]$-modules; see \cite{Keller, KM}, $q$ is the localization, $h$ is the homology functor, $\utot{k}$\label{index:S} is the functor defined by 
$\utot{0}(M)=\oplus_{i}M^{2i}$ and $\utot{1}(M)=\oplus_{i}M^{2i+1}$.
We remark that
$(\utot{0}(M))^i= M^{2i}$ and $(\utot{1}(M))^i= M^{2i+1}$.


\todo{[RR1, 2.2] We put this paragraph to summarize our functors.}
A functor $h^k :  \K[u]\text{-}\mathsf{Ch} \to \K[t]\text{-}\mathsf{grMod}$ is defined by $h^k =s^khq$.  
We have the natural functor
$\theta:=\overline{q^{\mathbb R}} : \text{Ho}(\mathsf{Ch}_\K^{(\R, \leq)}) \to  \text{Ho}(\mathsf{Ch}_\K)^{(\R, \leq)}$,
which is induced by the localization functor $q$.
Let $\nu^k:  \K[u]\text{-}\mathsf{Ch} \to \mathsf{Mod}_\K^{(\R, \leq)}$ be the composite $(\lfloor \ \rfloor)^*\circ\eta^k \circ (H)_*\circ C$.\label{index:nu}
Moreover,
let $\alpha : \K[u]\text{-}\mathsf{Ch} \to \mathsf{Ch}_\K^{(\R, \leq)}$ \label{index:alpha} be the functor $(\lfloor \ \rfloor)^* \circ C$.
Then, we observe that $d_{\text{HC}}(\alpha M, \alpha N)= d_{\text{I}}((\theta\circ \mu\circ q)M, (\theta\circ \mu\circ q)N)$
and
$d_{\text{IHC}}(\alpha M, \alpha N) = d_{\text{I}}((\mu\circ q)M, (\mu\circ q)N)$; see Subsection \ref{sect:ID_up_to_H} for the distance
$d_{\text{HC}}$ and $d_{\text{IHC}}$.

The pair $(\mathcal{BAR}, d_B)$ stands for the set of barcodes (multisets of intervals) endowed with the bottleneck distance and
$(\mathsf{Mod}_\K^{(\R, \leq)}, d_{\text{I}})$ is the {\it class} of the diagrams endowed with the interleaving distance.
The wave arrows denote the assignments of the objects, where
$\mathcal{B}_{(\ )}$ is defined in the class of locally finite $\K[t]$-modules. Here, by definition,  a $\K[t]$-module $M$ is locally finite if the vector space
$M^i$ is of finite dimension for each $i$.

The result \cite[Theorem 4.16]{B-S} asserts that $\chi$ gives an isometric embedding if the domain is restricted to the set of finite barcodes.
%
%
Moreover, the functors $\mu$ and $\xi^k$ in (\ref{eq:diagrams-diagram})
are induced by the universality of the homotopy categories $\text{D}(\K[u])$ and
$\text{Ho}(\mathsf{Ch}_\K^{(\R, \leq)})$, respectively.

\begin{rem}\label{rem:Ho} \todo{[RR1, 2.4 Comment 12] We do not take the cofibrant replacement for $G$.}
We may regard the set of morphisms from $F$ to $G$ in $\text{Ho}(\mathsf{Ch}_\K^{(\R, \leq)})$ as the homotopy set of maps from $\widetilde{F}$ to
$G$, where $\widetilde{F}$ denotes the cofibrant replacements of $F$. In this manuscript, we do not use an explicit form of the cofibrant replacement; see \cite[Section 11.6]{H} for the form.
Observe that all objects in $\text{Ho}(\mathsf{Ch}_\K^{(\R, \leq)})$ are fibrant; see \cite[Theorem 1.4]{B-M-R}.
\end{rem}

\begin{rem}\label{rem:h} (i)
It follows from the definition of the functor $h^k$ that $H(M) \cong H(N)$ for $M$ and $N$ in $\K[u]\text{-}\mathsf{Ch}$ if
$h^kM\cong h^kN$ for $k = 0$ and $1$.

(ii) Every dg $\K[u]$-module $M$ is formal in the sense that $M \cong (H(M), 0)$ in
$\text{D}(\K[u])$. This fact follows from the proof of Lemma \ref{lem:formality}.

\end{rem}

\begin{defn}\label{defn:CohI} \label{index:eo_CohI} \todo{[RR1, 2.2 ii) The definition is revised with the interleaving distance.}
Let $M$ and $N$ be dg $\K[u]$-modules. Then, the {\it even cohomology interleaving distance}
$d_{\text{CohI}}^0(M, N)$ and the {\it odd cohomology interleaving distance} $d_{\text{CohI}}^1(M, N)$ are defined
by $d_{\text{I}}(\nu^0(M), \nu^0(N))$ and $d_{\text{I}}(\nu^1(M), \nu^1(N))$, respectively.
\end{defn}

We observe that
\begin{eqnarray}\label{eq:CohI--IHC}
d_{\text{CohI}}^k(M, N)&=& d_{\text{I}}(\nu^kM, \nu^kN) \\
& \leq & d_{\text{I}}((\theta\circ \mu\circ q)M, (\theta\circ \mu\circ q)N) \nonumber \\
&\leq& d_{\text{I}}((\mu\circ q)M, (\mu\circ q)N) \nonumber
\end{eqnarray}
for $k= 0$, $1$ and dg $\K[u]$-modules $M$ and $N$.
The inequalities follow from \cite[Proposition 3.6]{B-S}.

By Theorem \ref{thm:I_and_B} and the commutativity of the diagram (\ref{eq:diagrams-diagram}), we establish
\begin{prop}\label{prop:dI-dB} Let $M$ and $N$ be dg $\K[u]$-modules whose homologies are of finite dimension for each degree. Then,
$d_{\text{\em CohI}}^k(M, N) = d_\text{\em I}(\chi(\mathcal{B}_{\utot{k}hqM}), \chi(\mathcal{B}_{\utot{k}hqN})) =
d_{\rm B}(\mathcal{B}_{\utot{k}hq(M)}, \mathcal{B}_{\utot{k}hq(N)})$ for $k = 0$ and $1$.
\end{prop}

The following proposition shows the reason why we consider the interleaving distances $d_{\text{CohI}}^k(M, N)$ for $k =0$ and $1$ only.
\begin{prop}\label{prop:nu}
For each $l \in {\mathbb Z}$, it holds that $d_{\text{\em CohI}}^0(M, N) = d_{\text{\em I}}(\nu^{2l}M, \nu^{2l}N)$ and
$d_{\text{\em CohI}}^1(M, N) = d_{\text{\em I}}(\nu^{2l+1}M, \nu^{2l+1}N)$.
\end{prop}

\begin{proof}
We recall a translation functor $(l) : (\R, \leq) \to (\R, \leq)$ defined by $(l)t = t+ l$ and the functor
$(l)^* : \mathsf{Mod}_\K^{(\R, \leq)} \to \mathsf{Mod}_\K^{(\R, \leq)}$ induced by $(l)$. We see that
$(\nu^{2l})M = (l)^*(\nu^0)M$ and $(\nu^{2l+1})M = (l)^*(\nu^1)M$. It follows that
\begin{eqnarray*}
d_{\text{CohI}}^0(M, N) = d_{\text{I}}((-l)^*(\nu^{2l})M, (-l)^*(\nu^{2l})N)&\leq & d_{\text{I}}((\nu^{2l})M, (\nu^{2l})N) \\
&\leq& d_{\text{CohI}}^0(M, N).
\end{eqnarray*}
By the same argument as above, we have the second equality.
\end{proof}

We now state our main theorem of this section. 

\begin{thm}\label{thm:IHC=CohI} \todo{[RR1, 2.2 iii)] The assertion is revised with the final equality only.}
The equality
\[
d_{\text{\em HI}}(\alpha M, \alpha N) = \text{\em max}\{d_{\text{\em CohI}}^k(M, N) \mid k = 0, 1 \}
\]
holds for dg $\K[u]$-modules $M$ and $N$. 
\end{thm}


In what follows, we may write $d_{\text{CohI}}(M, N)$\label{index:CohI_total} for $\text{max}\{d_{\text{CohI}}^k(M, N) \mid k = 0, 1 \}$ and call it the {\it cohomology interleaving distance} of dg $\K[u]$-modules $M$ and $N$.

\begin{rem} (i)
By the commutativity of the diagram (\ref{eq:diagrams-diagram}) and Proposition \ref{prop:nu},      \todo{[RR1, 2.4 Comment 13] The sentence `By the commutativity ...' is left as a remark. }
we see that
$d_{\text{CohI}}(M, N) = d_{\text{CohI}}(\alpha M, \alpha N)$ for dg $\K[u]$-modules $M$ and $N$, where the right-hand side stands for the distance of persistence modules described in Definition \ref{defn:H_HC}. \\
(ii) Combining Theorems \ref{thm:IHC=CohI} and \ref{thm:Equalities}, we have the equalities
\[
d_{\text{HC}}(\alpha M, \alpha N) =
d_{\text{IHC}}(\alpha M, \alpha N) = d_{\text{HI}}(\alpha M, \alpha N) = d_{\text{CohI}}(M, N).
\]
\end{rem}





\begin{proof}[Proof of Theorem \ref{thm:IHC=CohI}] We recall the inequalities (\ref{eq:CohI--IHC}). In order to prove the assertion, it suffices to show that
$d_{\text{HI}}(\alpha M, \alpha N) \leq d_{\text{CohI}}(M, N)=:\e$.
We observe that the funcor $\alpha$ commutes with taking homology. 
Then,
Lemma \ref{lem:formality} allows us to deduce that $\alpha L \simeq H(\alpha L)=\alpha H(L)$ for a dg $\K[u]$-module $L$, 
where $\alpha = (\lfloor \ \rfloor)^* \circ C$ by definition.
Therefore, it follows that $d_{\text{HI}}(\alpha M, \alpha N) \leq  d_{\text{I}}(\alpha H(M), \alpha H(N))$.
Moreover, with the same notation as in the proof of Proposition \ref{prop:nu}, we see that for each dg $\K[u]$-module $L$,
\[
\alpha H(L) = \oplus_{l \in {\mathbb Z}}((l)^*\nu^0H(L)\oplus (l)^*\nu^1H(L)),
\]
where $(l)^*\nu^kH(L)$ is regarded as a persistent dg module concentrated at degree $2l+k$ for $k = 0$ and $1$.
We have $\e\geq  d_{\text{CohI}}^k(M, N)= d_{\text{CohI}}^k(H(M), H(N)) = d_{\text{I}}(\nu^kH(M), \nu^kH(N))$. 
\end{proof}


\begin{prop}\label{prop:Ametric} Let $M$ and $N$ be dg $\K[u]$-modules. Suppose that $d_{\text{\em CohI}}^k(M, N) < \frac{1}{2}$. Then $h^kM \cong h^kN$ as a $\K[t]$-module. Moreover, if $d_{\text{\em CohI}}(M, N) < \frac{1}{2}$, then $M\cong N$ in $\text{\em D}(\K[u])$.
Thus, the distance $d_{\text{\em CohI}}$ is an extended metric on $\text{\em D}(\K[u])$. 
\end{prop}

\begin{proof} 
Let $F$ and $G$ be the persistence modules $\nu^kM$ and $\nu^kN$, respectively. By the assumption, there exists a positive real number $\e$ less than
$\frac{1}{2}$ such that $F$ and $G$ are $\e$-interleaved. For each integer $i$, we consider the commutative triangles in the diagram (\ref{eq:interleaving}).  In view of a property of the floor function,  we see that $F(i \to i+2\e)$ and $G(i \to i+2\e)$ are isomorphism for each integer $i$. Therefore, the maps
$\varphi(i)$ and $\varphi(i+\e)$ are injective and surjective, respectively.
Moreover, we have a commutative diagram 
\begin{eqnarray*}
\xymatrix@C45pt@R20pt{
F(i) \ar[r]^-{F(i\to i+\e)}_-{\cong}  \ar[dr]^(0.7){\varphi(i)} &  F(i+\e)  \ar[r] \ar[rd]^(0.7){\varphi(i+\e)}|\hole & F(i+2\e)   \\
G(i)  \ar[r]  \ar[ur]^(0.3){\psi(i)}|\hole& G(i+\e) \ar[ur]^(0.3){\psi(i+\e)} \ar[r]^-{\cong}_-{G(i+\e \to i+2\e)}  & G(i+2\e)
}
\end{eqnarray*}
in which horizontal arrows $F(i\to i+\e)$ and $G(i+\e \to i+2\e)$ are isomorphisms. Thus, it follows that $\varphi(i)$ is an isomorphism for each integer $i$. We observe that $G(i+\e)=G(i)$. It turns out that $h^0 M = \oplus_i H^{2i}M = \oplus_i\nu^0M (i) \cong \oplus_i\nu^0N (i) =\oplus_i H^{2i}N = h^0 N$ and
$h^1 M = \oplus_i H^{2i+1}M  =  \oplus_i\nu^1M (i)\cong  \oplus_i\nu^1N (i) = \oplus_i H^{2i+1}N = h^1 N$ as
$\K[t]$-modules.

Suppose that  $d_{\text{CohI}}(M, N) < \frac{1}{2}$. Remark \ref{rem:h} (i) yields that $H(M)\cong H(N)$ as $\K[u]$-module. Moreover,
it follows from Remark \ref{rem:h} (ii) that $M\cong N$ in $\text{D}(\K[u])$. We have the result.
\end{proof}

\begin{rem} 
We regard the category $\K[u]\text{-}\mathsf{Ch}$ as a category with a flow in the sense of de Silva, Munchi and Stefanou \cite{SMS18}. In fact, the flow $\mathcal{T} : ([0,  \infty), \leq), +, 0) \to \text{End}(\K[u]\text{-}\mathsf{Ch})$ is defined by
$T_\e(x) := \lfloor \e \rfloor u \cdot x$ for $x \in M$. Since $\lfloor \e \rfloor + \lfloor \delta \rfloor \leq \lfloor \e+\delta \rfloor$,
we have a natural transformation $\mu_{\e, \delta} : \mathcal{T}_{\e}\mathcal{T}_{\delta} \Longrightarrow \mathcal{T}_{\e+\delta}$.
In general, we have the {\it interleaving distance $d_{(\mathcal{C}, \mathcal{T})}$ with respect to $\mathcal{T}$} for objects of a category $(\mathcal{C}, \mathcal{T})$ with a flow $\mathcal{T}$; see \cite[Section 2.2]{SMS18}.
One can see that $d_{(\K[u]\text{-}\mathsf{Ch}, \mathcal{T})} = d_{\text{CohI}}$. \todo{[RR1, 2.2 iv] The last sentence is revised.}
\end{rem}

 The value of the cohomology interleaving distance between dg $\K[u]$-modules
 is {\it realized} with the homotopy interleaving via the functor $\alpha : \K[u]\text{-}\mathsf{Ch} \to \mathsf{Ch}_\K^{(\R, \leq)}$
 under  some  finiteness condition on the given dg $\K[u]$-modules.
 \todo{Remark\ref{rem:revision_of_rem2.6} and Proposition \ref{prop:D_E} are new faces.}

To describe the result more precisely, 
we recall a finitely presented persistence module in the sense of Lesnick \cite{Lesnick}.
While a representation for a multiparameter persistence module is introduced in \cite[Section 4]{Lesnick}, we here define that for one parameter.
Let $\mathcal{W}$ be an $\R$-graded set and $\ktr[\mathcal{W}]$ the free $\ktr$-module generated by $\mathcal{W}$. Observe that $\ktr$ is nothing but the ring $P_1$ in \cite{Lesnick}. 
By definition, 
a {\it presentation} of a $\ktr$-graded module $M$ is a pair $(\mathcal{W}, \mathcal{Y})$ for which $M \cong \ktr[\mathcal{W}]/\langle \mathcal{Y} \rangle$ as a $\ktr$-graded module,
where $\mathcal{Y}$ is a set of homogeneous elements of $\ktr[\mathcal{W}]$ and $\langle \mathcal{Y} \rangle$ is a submodule generated by
$\mathcal{Y}$. We say that $M$ is {\it finitely presented} if $M$ admits a presentation $(\mathcal{W}, \mathcal{Y})$ with
 $\mathcal{W}$ and  $\mathcal{Y}$ finite.

\begin{rem}\label{rem:revision_of_rem2.6}
Let $G$ and $F$ be finitely presented persistence modules.
Then, the {\it closure theorem} \cite[Theorem 6.1]{Lesnick} says that $F$ and $G$ are $\delta$-interleaved for $F,G\in \mathsf{Mod}_\K^{(\R, \leq)}$ if $d_{\text{I}}(F,G)=\delta$.
 Then, we see that the interleaving distance is an extended metric on isomorphism classes of finitely presented ({\it multidimensional}) persistence modules. 

We consider the full subcategory $\mathcal{F}$ of $\mathsf{Mod}_\K^{(\R, \leq)}$ consisting of the image of finitely generated graded $\K[t]$-modules by the functor $(\lfloor \ \rfloor)^*\circ \gamma^{-1}$; see the diagram (\ref{eq:diagrams-diagram}).
Then, it follows that
each object $M$ in $\mathcal{F}$ is  finitely presented. To see this,
define the functor $\Theta$ by the composite
\[
\xymatrix@C25pt@R20pt{
\K[t]\text{-}\mathsf{grMod}  & \mathsf{Mod}_\K^{(\Z, \leq)}  \ar[r]^{(\lfloor \ \rfloor)^*}_{\text{embedding}} \ar[l]^-{\simeq}_-\gamma &
 \mathsf{Mod}_\K^{(\R, \leq)} \ar[r]^-{\gamma'}_-{\simeq}& \ktr\text{-}\mathsf{grMod}.
}
\]
Here,  $\ktr\text{-}\mathsf{grMod}$ is the category of $\R$-graded $\ktr$-modules and the equivalence $\gamma'$ is defined by the same fashion as that of $\gamma$; see the paragraph before the diagram (\ref{eq:complexes}).
Then, we see that the functor $\Theta$ assigns the $\ktr$-module $\ktr$ to the $\K[t]$-module $\K[t]$.  Observe that $\Theta$ is an exact functor.
Thus, it follows that a free $\K[t]$-module is mapped to a free $\ktr$-module by $\Theta$.
Moreover, by virtue of Lemma \ref{lem:graded_pid},  we see that for a finitely generated $\K[t]$-module $M$, the $\ktr$-module $\Theta(M)$ is finitely represented.
Thus, objects $K$ and $L$ in $\mathcal{F}$ are $\e$-interleaved whenever $d_{\rm I}(K, L) = \e$. 
In particular, the interleaving distance is an extended metric on $\mathcal{F}$.
\end{rem}

\begin{prop}\label{prop:D_E}
Let $M$ and $N$ be dg $\K[u]$-modules whose homologies are finitely generated $\K[u]$-modules.  
Suppose that $d_{\text{\em CohI}}(M, N) =\e$. Then $\alpha M$ and $\alpha N$ are $\e$-homotopy interleaved and not $\delta$-homotopy interleaved for any $\delta$ less than $\e$ in the category $\mathsf{Ch}_\K^{(\R, \leq)}$.
\end{prop}

As seen in Section \ref{sect:BS^1}, the singular cochain functor assigns persistence modules to spaces over $BS^1$.
Thus, the cohomology interleaving distance for spaces over $BS^1$ is defined; see the beginning of Section \ref{sect:BS^1} for more detail.
In Appendix \ref{sect:toy_examples}, we compute explicitly the cohomology interleaving distances for spaces over $BS^1$.
All of the cohomology groups with coefficients in $\K$ of spaces we deal with in the appendix are finitely generated modules over $H^*(BS^1; \K) \cong \K[u]$.
While the computations use barcodes via Lemma \ref{lem:finiteIntervals_d} and Theorem \ref{thm:I_and_B} (the isometry theorem), Proposition \ref{prop:D_E} allows us to conclude that there are indeed homotopy interlevings between the singular cochain complexes of the spaces in the category $\mathsf{Ch}_\K^{(\R, \leq)}$ with the functor $\alpha$.


\begin{proof}[Proof of Proposition \ref{prop:D_E}] 
By assumption, we see that $(\lfloor \ \rfloor)^*\circ \gamma^{-1}\circ h^k)(M)$ and $(\lfloor \ \rfloor)^*\circ \gamma^{-1}\circ h^k)(N)$ are finitely generated $\K[t]$-modules for $k = 0$ and $1$.
The argument in Remark \ref{rem:revision_of_rem2.6} enables us to deduce that the persistence modules are $\e$-interleaved for $k = 0$ and $1$.
We observe that
\begin{align}
d_{\text{CohI}}(M, N) &= \text{max}\{d_{\text{CohI}}^k(M, N) \mid k = 0, 1 \} \\
&=  \text{max}\{d_{\text{I}}(\nu^kH(M), \nu^kH(N)) \mid k = 0, 1 \} \\
 &= \text{max}\{d_{\text{I}}((\lfloor \ \rfloor)^*\circ \gamma^{-1}\circ h^k)(M), (\lfloor \ \rfloor)^*\circ \gamma^{-1}\circ h^k)(N)) \mid k = 0, 1 \} .
 \end{align}
 The first and second equalities follow from the definition. The commutativity of the diagram (\ref{eq:diagrams-diagram}) yields the third one.
Then $(H)_*\alpha M$ and  $(H)_*\alpha N$ are $\e$-interleaved. Proposition \ref{prop:Equality} enables us deduce that $(H)_*\alpha M \simeq \alpha M$ and $(H)_*\alpha N \simeq \alpha N$. 
The latter half follows from Theorem \ref{thm:IHC=CohI}. This completes the proof.
\end{proof}

\begin{rem}
We consider a persistence dg module $M$, namely an object in
$\mathsf{Ch}_\K^{(\R, \leq)}$. Suppose that $M$ is $H$-formal whenever
$(\gamma' \circ \eta^k(H)_*)(M)$ is finitely presented for each $k$; see (\ref{eq:diagrams-diagram}) for the functors.
Then, we may generalize Proposition \ref{prop:D_E}] to an assertion for objects $\mathsf{Ch}_\K^{(\R, \leq)}$. We do not pursue this issue in this manuscript.
\end{rem}

\subsection{A filtered $\K[t]$-module}  
As mentioned in the Introduction, we consider a filtered $\K[t]$-module, where $\deg t = 1$.

Our motivated example of such a module comes from a map $X \to BS^1$.
For example, if the map is a fibration, we may use the Leray--Serre spectral sequence when computing the cohomology of $X$. Then, we have to consider the extension problems which appear in the $E_\infty$-term even if all terms are computed. The lemma below asserts that there is no extension problem as $H^*(BS^1; \K)$-modules; see Remark \ref{rem:even-odd}.

Let $H^*$ be a non-negatively graded $\K[t]$-module with a filtration
\[
H^p= F^0 H^p\supset F^1H^p \supset \cdots \supset F^iH^p \supset \cdots \supset F^{p +1}H^p = 0
\]
of $\K[t]$-submodules for $p\geq 0$. Suppose that $t F^iH^p \subset F^{i + 1}H^{p+1}$. Then, we have a bigraded $\K[t]$-module $E^{*,*}$ defined by $E^{p, q} := F^pH^{p+q}/ F^{p+1}H^{p+q}$. Observe that $t \cdot E^{p, q} \subset E^{p+1, q}$. 

For a bigraded module $E^{*,*}$, we define a graded module $\text{Tot}\  \!E^{*,*}$, which is called the {\it total complex of} $E^{*,*}$, by $(\text{Tot}\  \!E^{*,*})^i := \oplus_{p+q = i}E^{p,q}$.

\begin{lem}\label{lem:bigradedModules}
As a graded $\K[t]$-module, $\text{\em Tot} \ \!E^{*,*} \cong H^*$ provided $\dim H^i < \infty$ for each $i$.
\end{lem}

\todo{[RR1, 2.4 Comment 14] The explanation on the equivalence classes is revised.}
\begin{proof}
%
%
We say that an element $a\in H^*$ has the {\it filtration degree} $p$, denoted $\text{fil-deg}\ a = p$, if $a \in F^p$ and
$a \notin F^{p+1}$.
We prove the lemma by the induction on the degrees and filtration degrees of elements $a_\lambda^k$ of $H^*$
which desend to a basis $\{a_\lambda^k \}_{k\geq 0, \lambda}$ of
$H^* / (t) H^*$, where $\deg a_\lambda^k = k$.

Let $S_k$ be the subset $\{a_{\lambda_1}^k, \ldots,  a_{\lambda_{s_k}}^k\}$ of $H^*$ which gives rise to linearly independent elements
of $E^{*,*} / (t) E^{*,*}$ with degree $k$. We may view $S_k$ as a subset of $H^*$ with
$\text{fil-deg} \ \!a_{\lambda_i}^k \geq \text{fil-deg} \ \!a_{\lambda_j}^k$ for $i < j$. Let $[S_k]$ be the subset  $\{[a_{\lambda_1}^k], \ldots, [a_{\lambda_{s_k}}^k]\}$ of $E^{*,*}$,
where $[a]$ denotes the image of $a$ by the projection of $F^{\text{fil-deg} a} \to E^{\text{fil-deg} a, *}$.

Since $F^iH^0 =0$ for $i> 0$, it follows that
the $\K[t]$-submodule of $H^*$ generated by $S_0$ coincides with
that of $E^{*,*}$ generated by $[S_0]$.
Assume that the map $\varphi_k$ defined by $\varphi_k ([a_\mu])= a_\mu$ is an isomorphism from the $\K[t]$-submodule $E_k$ of $\text{Tot} \ \!E^{*,*}$ generated by $[S_0]\cup\cdots \cup [S_k]$ to the $\K[t]$-submodule $H_k$ of $H^*$ generated by $S_0\cup\cdots \cup S_k$.
We may replace elements in $S_i$ and $[S_i]$ to construct the isomorphism preserving the linear independence of the elements in each set if necessary\footnote{This procedure is clarified by considering the induction step described below.}.

Suppose that $[a_{\lambda_1}^{k+1}]t^n \neq 0$ in $(E_k+ \K[t]\cdot[a_{\lambda_1}^{k+1}])/E_k$ for each $n \geq 0$.
Then, it follows that  $E_k+ \K[t]\cdot[a_{\lambda_1}^{k+1}] = E_k\oplus \K[t]\cdot [a_{\lambda_1}^{k+1}]$. Thus, we have an isomorphism
$E_k\oplus \K[t]\cdot [a_{\lambda_1}^{k+1}] \stackrel{\cong}{\to} H_k\oplus \K[t]\cdot a_{\lambda_1}^{k+1}$ extending $\varphi_k$.

Suppose that $[a_{\lambda_1}^{k+1}]t^m = 0$ and $[a_{\lambda_1}^{k+1}]t^{m-1} \neq 0$ in $Z_k := (E_k+ \K[t]\cdot [a_{\lambda_1}^{k+1}])/E_k$ for some $m \geq 1$. Then, we have
\begin{eqnarray}\label{eq:filtrations}
a_{\lambda_1}^{k+1}t^m - \sum_{i\leq k, j} \beta_{j, i} a_{\lambda_j}^{i} t^{l_{j, i}}= 0
\end{eqnarray}
for some nonzero elements $\beta_{j, i} \in \K$. By degree reasons, we see that ${l_{j, i}} \geq m$. We rewrite $a_{\lambda_1}^{k+1}$ for
$a_{\lambda_1}^{k+1} - \sum_{i\leq k, j} \beta_{j, i} a_{\lambda_j}^{i} t^{l_{j, i}-m}$. Then,
we have an isomorphism
\begin{eqnarray*}
\varphi_{(1)} : E_k\oplus \Sigma^{\deg [a_{\lambda_1}^{k+1}]}(\K[t]/(t^m)) \cong E_k +  \K[t]\cdot [a_{\lambda_1}^{k+1}] & \\
& \hspace*{-3.5cm} \longrightarrow H_k +  \K[t]\cdot a_{\lambda_1}^{k+1} \cong H_k\oplus\Sigma^{\deg a_{\lambda_1}^{k+1}}(\K[t]/(t^m) )
\end{eqnarray*}
defined by $\varphi_{(1)}([a_{\lambda_1}^{k+1}])= a_{\lambda_1}^{k+1}$ and $\varphi_{(1)} |_{E_k} = \varphi_{k}$.
The isomorphism in the codomain of the map above follows from the fact that $t^l\cdot a_{\lambda_1}^{k+1}$ is not in $H_k$ for $l < m$.
Observe that if $t^l\cdot a_{\lambda_1}^{k+1}$ is in $H_k$ for some $l < m$, then $t^l\cdot [a_{\lambda_1}^{k+1}]=0$ in $Z_k$.
\todo{[RR1, 2.4 Comment 15] Explain why the element is not in $H_k$ more precisely. }
This contradictions the choice of the positive integer $m$.
Moreover, we see that $a_{\lambda_1}^{k+1}\cdot t^{m-1}\neq 0$. In fact, if $a_{\lambda_1}^{k+1}\cdot t^{m-1}=0$, then for the original
$a_{\lambda_1}^{k+1}$,
$[a_{\lambda_1}^{k+1}]\cdot t^{m-1}=0$ in $Z_k$,
which is a contradiction.

Moreover, by applying the same construction of the isomorphism as above to elements $[a_{\lambda_2}^{k+1}], \ldots , [a_{\lambda_{s_{k+1}}}^{k+1}]$,
we obtain an isomorphism $\varphi_{k+1} : E_{k+1} \stackrel{\cong}{\to} H_{k+1}$. As a consequence, we have an isomorphism
$\varphi : \text{Tot} \ \!E^{*,*} \to H^*$ of $\K[t]$-modules.
\end{proof}


\begin{rem} \label{rem:even-odd}
Let $M^*$ be a non-negatively graded $\K[u]$-module with a filtration
$M^p= F^0M^p \supset F^1M^p \supset \cdots \supset F^i M^p\supset \cdots \supset F^{p +1}M^p = 0$ of $\K[u]$-submodules for $p\geq 0$.
Then, we may  apply the same argument as in the proof of Lemma \ref{lem:bigradedModules} to $M^*$ provided
$uF^iH^p \subset F^{i+1}H^{p+2}$ for $i\geq 0$.
As a consequence, we see that $\text{Tot}\ \!E^{*,*} \cong M^*$ as a graded $\K[u]$-module.
\end{rem}

\begin{rem}
For $q\in \mathbb{Z}$, define a new graded
$\K[t]$-module $H^{*,q}$ by setting $H^{p,q}:= F^{p-q}H^p/F^{p-q+1}H^p$ and letting $t$ act by the multiplication
$H^{p,q} \to H^{p+1, q}$. Then Lemma \ref{lem:bigradedModules} implies that there is an isomorphism
$H^*\cong \oplus_{q\in \mathbb{Z}}H^{*, q}$ of graded $\K[t]$-modules.
\end{rem}

So far we consider dg $\K[v]$-modules with $\deg v =1$ or $2$. Even if the degree of $u$ is a positive integer, the same arguments as in this section are applicable to the case and the results remain true with an appropriate degree shift.  For example, the functor $C$ in (\ref{eq:diagrams-diagram}) is replaced with one defined by
$C(\{M^l, \partial \})(i) = \Sigma^{(\deg u)i}M$ for a dg $\K[u]$-module $\{M^l, \partial \}$.


\section{The cohomology interleaving of spaces over $BS^1$}\label{sect:BS^1}

Let $\K$ be a field. Unless otherwise explicitly stated, it is assumed that a space $X$
is connected and the singular cohomology of $X$ with coefficients $\K$ is locally finite;  that is, the $i$th cohomology of $X$ is of finite dimension for $i \geq 0$.

Let $f : X \to BS^1$ be a space over $BS^1$.
We have a quasi-isomorphism $\kappa: \K[u] \to C^*(BS^1; \K)$ and the morphism $f^* : C^*(BS^1; \K)\to C^*(X; \K)$ of differential graded algebras (DGAs). Then, the singular cochain complex $C^*(X; \K)$ is regarded as a $\K[u]$-module via the maps $f^*\circ \kappa$.
We may write $C^*(X; \K)_f$ for the $ \K[u]$-module.

For example, the set of the isomorphism classes of principal $S^1$-bundles over a space $X$
is isomorphic to the homotopy set $[X, BS^1]$. As mentioned below, an $S^1$-space $Y$ gives rise to an fibration over $BS^1$ with $Y$ the fibre via the Borel construction. As $S^1$-spaces, we may deal with the free loop spaces of a space and a manifold, which are main objects in string topology \cite{C-S, CHV}.
Thus, we have a lot of interesting examples of dg $\K[u]$-modules from such spaces.

The cohomology interleaving distances in Definition \ref{defn:CohI} give
the {\it even and odd cohomology interleaving distances} $d_{\text{CohI}, \K}^k((X, f), (Y, g))$, for brevity  $d_{\text{CohI}, \K}^k(X, Y)$, \label{index:CohI_space}
between the spaces $f\colon X\to BS^1$ and
$g\colon Y\to BS^1$ over $BS^1$ defined
by the distance $d_{\text{CohI}}^k(C^*(X; \K)_f, C^*(Y; \K)_g)$ for $k = 0$ and $1$, respectively.
Moreover, we define the  {\it cohomology interleaving distance} by
\begin{align}
d_{\text{CohI}, \K}((X, f), (Y, g))&=d_{\text{CohI}, \K}(X, Y) \\
&:=\text{max}\{d_{\text{CohI}}^k(C^*(X; \K)_f, C^*(Y; \K)_g) \mid k = 0, 1 \}.
\end{align}
Theorem \ref{thm:IHC=CohI} implies that the distance
$d_{\text{CohI}, \K}(X, Y)$ determines $d_{\text{HC}}$, $d_{\text{IHC}}$ and $d_{\text{HI}}$ between $\alpha C^*(X; \K)$ and $\alpha C^*(Y; \K)$.
Our first result in this section is as follows.
\begin{prop}\label{prop:HomotopySet} 
The function $d_H : [X, BS^1]\times [X, BS^1] \to \R_{\geq 0}\cup \{\infty\}$
defined by $d_H([f], [g]) = d_{\text{\em CohI}, \K}((X, f), (X, g))$ is a well-defined extended pseudodistance on the homotopy set $[X, BS^1]$.
\end{prop}

\begin{proof} Suppose that maps $f$ and $g$ from $X$ to $BS^1$ are homotopic with a homotopy $K : X\times I \to BS^1$. Then $C^*(X; \K)_f$  and $C^*(X; \K)_g$ are quasi-isomorphic to $C^*(X\times I; \K)_K$ in $\K[u]\text{-}\mathsf{Ch}$. The result follows from Proposition \ref{prop:nu}.
\end{proof}

\begin{rem}\label{two-variable homotopy invariants} Let $\mathsf{Top}/BS^1$ be the category of spaces over $BS^1$ and $\zeta : \mathsf{Top}/BS^1 \to \K[u]\text{-}\mathsf{Ch}$ the functor defined by the composite $\alpha \circ C^*(\ ; \K)$. Theorem \ref{thm:IHC=CohI} allows us to obtain a numerical
two-variable homotopy invariants $d_{\text{CohI}, \K}= d_{\text{HC}}\circ (\zeta\times \zeta) = d_{\text{IHC}}\circ (\zeta\times \zeta) = d_{\text{HI}}\circ (\zeta\times \zeta)$ on $\mathsf{Top}/BS^1$, which are determined explicitly for some cases with the results in Appendix \ref{sect:toy_examples} and are evaluated with the results in this section. 
\end{rem}

Let $Y$ be an $S^1$-space. We consider the Borel construction $Y_{hS^1}: = ES^1\times_{S^1}Y$
which fits into the Borel fibration $Y \to Y_{hS^1} \stackrel{p}{\to} BS^1$. 
Observe that $\text{\rm pt}_{hS^1} \simeq BS^1$.
Let $LX$ be the free loop space of $X$, namely, the space of continuous maps from $S^1$ to $X$ endowed with the compact-open topology.
The rotation on $S^1$ induces the action $\rho : S^1\times LX \to LX$ on the free loop space.
Thus, we have the Borel fibration $p :(LX)_{hS^1} \to BS^1$.
For a space $X$, we denote by $l(X)_{\K}$ the integer $\text{max}\{i \mid H^{i}(X; \K)\neq 0, i\geq 0\}$.
We investigate the cohomology interleaving distance between spaces, which are in the classes defined below.

\begin{itemize}
\item
Class (I) consists of the Borel constructions $(LX)_{hS^1}$ of the free loop spaces $LX$ of simply-connected spaces $X$.
\item
Class (II) consists  of the spaces $X$ for each of which $X$ fits in a fibration $\mathcal{F} :  F \to X \to BS^1$ with $l(F)_{\K}< \infty$.
\item
Class (III) consists of the spaces $X \to BS^1$ over $BS^1$ with $l(X)_{\K}< \infty$. As a consequence, the local finiteness condition of the cohomology implies that $H^*(X; \K)$ is of finite dimension.
\end{itemize}

In order to exhibit our result on the cohomology interleaving distance between spaces in Class (I),
we here introduce the {\it BV-exactness} of a simply-connected space $X$; see \cite[Definition 2.9]{KNWY2023}. By definition, the BV-operator
$\Delta$ on $H^*(LX; \Q)$ is the composite
\[
\xymatrix@C20pt@R20pt{
\Delta: H^*(LX; \Q) \ar[r]^{\rho^*} & H^*(S^1\times LX; \Q) \ar[r]^-{\int_{S^1}} & H^{*-1}(LX; \Q),
}
\]
where $\int_{S^1}$ stands for the integration along the fundamental class $t$ of $S^1$; that is, the map is defined by
$\int_{S^1} (t\otimes x)= x$ for $t\otimes x \in H^*(S^1; \Q)\otimes H^*(LX; \Q)\cong H^*(S^1\times LX; \Q)$.

Originally, the BV-operator is defined on the {\it homology} of the free loop space of a manifold $M$ and gives the loop homology $\mathbb{H}_*(LM): = H_{*+\dim M}(LM)$ endowed with 
a Batalin--Vilkovisky algebra structure; see, for example, \cite[Chapter 1, 1.3]{CHV}. Our definition above is regarded as the dual version.

\begin{defn} \label{defn:BV-exact} A simply-connected space $X$ is {\it Batalin-Vilkovisky exact (BV-exact)}
if
$\im \widetilde{\Delta} = \ker \widetilde{\Delta}$,
where
$\widetilde{\Delta}\colon \widetilde{H}^*(LX; \Q)  \to \widetilde{H}^*(LX; \Q)$ is the restriction of the $BV$-operator to the reduced cohomology groups.
\end{defn}

We also recall the {\it $S$-action} on $H^*_{S^1}(LX; \Q):=H^*((LX)_{hS^1}; \Q)$ which is the multiplication $S :=\times u : H^*_{S^1}(LX; \Q) \to H^*_{S^1}(LX; \Q)$ defined by $S(x) := p^*(u)x$ for $x \in H^*_{S^1}(LX; \Q)$, where
$p : (LX)_{hS^1} \to BS^1$ is the projection. We view the one-point space $\text{pt}$ as the $S^1$-space with the trivial action.


\begin{thm}\cite[Theorem 2.11]{KNWY2023}\label{thm:BV-S} A simply-connected space $X$ is BV-exact if and only if the reduced $S$-action on
$\widetilde{H}^*_{S^1}(LX; \Q)$ is trivial, where $\widetilde{H}^*_{S^1}(LX; \Q)$ denotes the cokernel of the map
$\Q[u]\cong H^*_{S^1}(\text{\rm pt}; \Q) \to {H}^*_{S^1}(LX; \Q)$ induced by the trivial map.
\end{thm}

We call a simply-connected space $X$ {\it formal} if there exists a zig-zag of quasi-isomorphisms of differential graded algebras between the singular cochain algebra $C^*(X; \mathbb{Q})$ and the cohomology algebra $H^*(X; \mathbb{Q})$ with the trivial differential.

\begin{cor}\text{\em (\cite[Corollary 2.13]{KNWY2023})}
  \label{cor:BV}
  If a simply-connected space $X$ is formal,
  then it is BV-exact. 
\end{cor}

The cohomology interleaving distance between BV-exact spaces in Class (I) is determined explicitly.

\begin{prop}\label{prop:LX} Let $X$ and $Y$ be formal spaces, more generally, BV-exact spaces. 
Then, it holds that for $k=0$ and $1$,
\begin{eqnarray*}
d^k_{\text{\em CohI}, \Q}((LX)_{hS^1},  (LY)_{hS^1}) =
\left \{
\begin{array}{ll}
0 & \text{if $h^kC^*((LX)_{hS^1}; \Q) \cong h^kC^*((LY)_{hS^1}; \Q)$} \\
   & \text{as a $\Q[t]$-module,} \\
\frac{1}{2} & \text{otherwise}.
\end{array}
\right.
\end{eqnarray*}
In particular, $d_{\text{\em CohI}, \Q}((LX)_{hS^1},  (LY)_{hS^1}) = 0$ if and only if $C^*((LX)_{hS^1}; \Q) \cong C^*((LY)_{hS^1}; \Q)$ in
$\rm{D}(\Q[u])$.
\end{prop}

\begin{proof}
We first prove that the cohomology interleaving distance
$d^k_{\text{CohI}, \Q}$ is less than or equal to $\frac{1}{2}$.
For a simply-connected space $X$, the Sullivan minimal model for $(LX)_{hS^1}$ in \cite[Theorem A]{V-B} enables us to conclude that
the sequence
\[
0 \to H^*_{S^1}(\text{\rm pt}; \Q) \to {H}^*_{S^1}(LX; \Q) \to \widetilde{H}^*_{S^1}(LX; \Q) \to 0
\]
is a split exact one of $\Q[u]$-modules. Then, we see that  $1\cdot u^s \neq 0$ for each $s\geq 0 $ and  the unit $1 \in H^0_{S^1}(LX; \Q)$.
Moreover, it follows from the BV-exactness that $x\cdot u = 0$ for each element $x \in \widetilde{H}^i_{S^1}(LX; \Q)$ with $i > 0$.
Then, the barcode associated with $\utot{k}C^*((LX)_{hS^1}; \Q)$ for each $k = 0$ and $1$ consists of one interval $[0, \infty)$ and intervals of the form $[i, i+1)$. Observe that the interval $[0, \infty)$ appears in the barcode only if $k =0$. By \cite[Propositions 4.13]{B-S} and
Lemma \ref{lem:finiteIntervals_d}, we see that
$d_{\text{I}}(\chi_{[0, \infty)}, \chi_{[0, \infty)}) = 0$, $d_{\text{I}}(\chi_{[0, \infty)}, \chi_{[i, i+1)}) = \infty$, 
$d_{\text{I}}(\chi_{\emptyset}, \chi_{[0, \infty)}) = \infty$,
$d_{\text{I}}(\chi_{\emptyset}, \chi_{[i, i+1)}) = \frac{1}{2}$ and $d_{\text{I}}(\chi_{[j, j+1)}, \chi_{[i, i+1)}) \leq \frac{1}{2}$.
We consider the bottleneck distance between barcodes $\mathcal{B}_{\utot{k} H^*_{S^1}(LX; \Q)}$ and $\mathcal{B}_{\utot{k} H^*_{S^1}(LY; \Q)}$.
If a bijection $f$ in Definition \ref{defn:bottleneck_D} assigns
$[0, \infty)$ to $[0, \infty)$, then the supremum
$\sup_{I \in \text{dom} (f)}d_{\text{I}}(\chi_I, \chi_{f(I)})$
is less than or equal to $\frac{1}{2}$. On the otherwise, the supremum is infinite.
Therefore, Proposition \ref{prop:dI-dB} enables us to deduce that $d^k_{\text{CohI}, \Q}((LX)_{hS^1},  (LY)_{hS^1}) \leq \frac{1}{2}$ for $k =0$ and $1$.

Assume further that $M:=h^kC^*((LX)_{hS^1}; \Q)$ and $N:=h^kC^*((LY)_{hS^1}; \Q)$ are $\e$-interleaved for some $\e < \frac{1}{2}$.  By Proposition \ref{prop:Ametric}, we see that
$M\cong N$ as a $\Q[t]$-module and $C^*((LX)_{hS^1}; \Q)\cong C^*((LY)_{hS^1}; \Q)$ in $\text{D}(\Q[u])$.
\end{proof}

Before describing upper and lower bounds of the cohomology interleaving distance of spaces, we recall the {\it cup-length} $\text{cup}(f)_R$\label{index:cup1} of a map $f :X \to Y$ with coefficients in a commutative ring $R$.    \todo{[RR1, 2.4 Comment 16] `the coefficient in' to `coefficients in'.  }
By definition,  the integer $\text{cup}(f)_R$ is the length of the longest non-zero product in the image
of the homomorphism $f^* : \widetilde{H}^*(Y; R) \to \widetilde{H}^*(X; R)$ between the reduced cohomology groups.
We observe that $\text{cup}(f)_R \leq \text{cat}(f)$,
where  $\text{cat}(f)$ denotes the category of the map $f$, namely the least integer $n$ such that $X$ can be covered by $n+1$ open subsets $U_i$, for which the restriction of $f$ to each $U_i$ is nullhomotopic, see \cite[Proposition 1.10]{B-G}.

The following proposition gives a rough evaluation of the interleaving distance between spaces over $BS^1$.

\begin{prop}\label{prop:cup}
Let $v_X : X \to BS^1$ and  $v_Y : Y \to BS^1$ be spaces over $BS^1$.
 Then, it holds that for $k = 0$ and $1$,
\[
d^k_{\text{\em CohI}, \K}(X, Y) \leq  \frac{1}{2}\text{\em max}\{\text{\em cup}(v_X)_\K+1, \text{\em cup}(v_Y)_\K+1\}.
\]
In particular, the cohomology interleaving distances between spaces in \text{\em Class (III)} are finite.
\end{prop}

\begin{lem}\label{lem:cup} Let $v : X \to BS^1$ be a space over $BS^1$.
Then, the length of the longest bar $J$ in
$\mathcal{B}_{\utot{0}H^*(X; \K)}$ and $\mathcal{B}_{\utot{1}H^*(X; \K)}$ is less than or equal to $\text{\em cup}(v)_\K +1$.
\end{lem}

\begin{proof}
Let $n$ be the integer $\text{cup}(v)_\K +1$. Then, it follows from the definition of the cup-length that
$v^*(u)^n =0$ in $H^*(X; \K)$. Therefore, we see that $m_i v^*(u)^n =0$ for each element $m_i$ of a basis
$\{m_i\}_{i\in \Lambda}$ of $H^*(X; \K)/ (v^*(u))H^*(X; \K)$. This fact enables us to deduce that
the length of $J$ is less than or equal to $n$.
\end{proof}

\begin{proof}[Proof of Proposition \ref{prop:cup}]
The result follows from Lemmas \ref{lem:finiteIntervals_d} and \ref{lem:cup}.
\end{proof}

\begin{ex}\label{ex:ThreeClasses}

Let $(LM)_{hS^1}$ and $Y$ be in Classes (I) and (III), respectively.

(1) It follows that
$d^0_{\text{CohI}, \Q}((LM)_{hS^1},  Y)=\infty$. In fact, we see that $1\cdot t^l\neq 0$ for each $l \geq 1$ and the unit $1 \in H^0_{S^1}(LM; \Q)$.
The argument in Example \ref{ex:infinity} allows us to obtain the result.

(2) Let $F$ be the fiber of a fibration ${\mathcal F} : X \to BS^1$ in Class (II). Assume that the dimension of the cohomology $H^*(F; \K)$ is greater than or equal to $2$ and the Leray--Serre spectral sequence for $\mathcal{F}$ with coefficients in $\K$ collapses at the $E_2$-term.
Then, we see that
$d_{\text{CohI}, \K}(X,  Y)=\infty$ and $d_{\text{CohI}, \Q}(X,  (LM)_{hS^1})=\infty$
if $M$ is BV-exact.
These facts follow from Example \ref{ex:infinity}, Remark \ref{rem:even-odd} and Theorem \ref{thm:I_and_B}.
\end{ex}

Let $v:X\to BS^1$ be a space over $BS^1$.
We will denote by ${\rm cup}^k(v)_{{\mathbb K}}$\label{index:cup2} the largest positive integer $n$ such that the action of $t^n$ on $\utot{k} H^{*}(X;{\mathbb K})$ is nontrivial; see the diagram \eqref{eq:diagrams-diagram} for the functor $\utot{k}$.
\todo{[RR1, 2.4 Comment 17] The space has been put.}
Observe that the integer ${\rm cup}^0(v)_{{\mathbb K}}$ coincides with the cup-length of $f$ mentioned above: ${\rm cup}^0(v)_{{\mathbb K}} = {\rm cup}(v)_{{\mathbb K}}$.
Refer to the notation ${\rm cup}^k (C^*(X;{\mathbb K}))$ introduced later in Proposition \ref{prop:M_pt}.
It follows from the definition that
\begin{eqnarray}\label{eq:length}
{\rm cup}^k (C^*(X;{\mathbb K})) = {\rm cup}^k(v)_{{\mathbb K}}.
\end{eqnarray}

\begin{prop}\label{prop:X_pt}
Let $v:X\to BS^1$ a space over $BS^1$ in Class (III).
Then, the cohomology interleaving distance between $v:X\to BS^1$ and ${\rm pt} \to BS^1$ is computed as follows.
\[
d^k_{{\rm CohI},{\mathbb K}}(X,{\rm pt})
= \left\{
\begin{array}{cl}
 0 & (H^*(X;{\mathbb K})\cong {\mathbb K})\\
\frac{1}{2}({\rm cup}^k (v)_{{\mathbb K}} + 1) & (\text{otherwise}).
\end{array}
\right.
\]
\end{prop}

\begin{proof}
Proposition \ref{prop:M_pt} and \eqref{eq:length} yield the result.
\end{proof}

\begin{prop}\label{prop:cup_2}
Let $v_X : X\to BS^1$ and $v_Y : Y\to BS^1$ be spaces over $BS^1$ in Class (III).
Assume further that  $H^*(X; {\mathbb K}) \not\cong {\mathbb K}$ and $H^*(Y; {\mathbb K})\not\cong {\mathbb K}$.
Then, it holds that
\[
d^k_{{\rm CohI},{\mathbb K}}(X,Y) \geq \frac{1}{2}|{\rm cup}^k(v_X)_{{\mathbb K}} - {\rm cup}^k(v_Y)_{{\mathbb K}}|.
\]
\end{prop}
\begin{proof}
By the triangle inequality, we have
\begin{eqnarray*}
d^k_{{\rm CohI},{\mathbb K}}(X,Y) &\geq& |d^k_{{\rm CohI},{\mathbb K}}(X,{\rm pt}) - d^k_{{\rm CohI},{\mathbb K}}(Y,{\rm pt})|.
\end{eqnarray*}
Thus, Proposition \ref{prop:X_pt} allows us to deduce the result.
\end{proof}

An argument on a spectral sequence is helpful to consider the cohomology interleaving distance between given spaces over $BS^1$.

\begin{prop}\label{thm:fibrations} Let $\mathcal{F}_i :  F_i \to X_i \to BS^1$ be a fibration for $i = 1$ and $2$. Assume that $F_i$ is a connected and $H^*(F_i; \K)$ is locally finite for each $i$.
Let $\{{}_iE_r^{*,*}, d_r\}$ be the Leray--Serre spectral sequence for $\mathcal{F}_i$ with coefficients in $\K$.
Suppose that the spectral sequences collapse at the $E_{r+1}$-term.
Then,
\[
d_{\text{\em CohI}, \K}(X, Y)  
= d_{\text{\em IHC}}(\alpha(\text{\em Tot}({}_1E_r^{*,*}, d_r)), \alpha(\text{\em Tot}({}_2E_r^{*,*}, d_r))).
\]
In particular,
$
d_{\text{\em CohI}, \K}(X, Y) = d_{\text{\em IHC}}(\alpha(\text{\em Tot}({}_1E_2^{*,*}, 0)), \alpha(\text{\em Tot}({}_2E_2^{*,*}, 0)))
$
if the spectral sequences collapse at the $E_2$-term.
\end{prop}

\begin{proof}
Since the spectral sequences collapse at the $E_{r+1}$-term, it follows that
\begin{eqnarray*}
d_{\text{CohI}, \K}(\text{Tot}\ \! {}_1E_\infty^{*,*}, \text{Tot}\ \! {}_2E_\infty^{*,*})& = &
d_{\text{CohI}, \K}((\text{Tot}\  \!{}_1E_{r+1}^{*,*}, 0), (\text{Tot}\  \!{}_2E_{r+1}^{*,*}, 0)) \\
&=& d_{\text{IHC}}(\alpha(\text{Tot}({}_1E_r^{*,*}, d_r)), \alpha(\text{Tot}({}_2E_r^{*,*}, d_r))).
\end{eqnarray*}
Observe that Theorem \ref{thm:IHC=CohI} gives the second equality. The result follows from Lemma \ref{lem:bigradedModules}; see Remark \ref{rem:even-odd}.
\end{proof}

\begin{rem}
The same result as above holds for the cobar type Eilenberg--Moore spectral sequence converging to $H^*(X_{hS^1}; \K)$ for an $S^1$-space $X$; see, for example,
\cite{E-M}, \cite[Theorem 2.2, ii)]{K2024} for the spectral sequence.  In fact, let $\{E_r^{*,*}, d_r\}$ and $\{ '\!E_r^{*,*}, '\!d_r\}$ be the Eilenberg--Moore spectral sequences
converging to $H^*(X_{hS^1}; \K)$ and  $H^*(\text{pt}_{hS^1}; \K)$, respectively. We have the $S^1$-equivariant map $f : X \to \text{pt}$.
Then, the naturality of the multiplicative spectral sequence gives a morphism $\{f_r \} : \{'\!E_r^{*,*}, '\!d_r\} \to \{E_r^{*,*}, d_r\}$
of spectral sequences with
\[
f_2 :  \K[u]\cong \!
 \ '\!E_2^{p,q}\cong \text{Cotor}^{*,*}_{H^*(S^1)}(\K, \K) \to E_2^{*,*}\cong \text{Cotor}^{p,q}_{H^*(S^1)}(\K, H^*(X)),
\]
where $\text{bideg} \ \! u = (1, 1)$. Thus, the spectral sequence $\{E_r^{*,*}, d_r\}$ has a dg $\K[u]$-module structure which is compatible with the $\K[u]$-module structure on
$H^*(X_{hS^1}; \K)$.
\end{rem}

The following corollary provides an approach for computing the interleaving distance between spaces in Class (II).

\begin{cor}\label{prop:E_2} Let $\mathcal{F}_i :  F_i \to X_i \to BS^1$ be a fibration with connected fiber for $i = 1$ and $2$.
Suppose further that for each $i$,  the spectral sequence for $\mathcal{F}_i$ collapses at the $E_2$-term and $l(F_i)_\K <\infty$. Then, the equality
\[
d^k_{\text{\em CohI}, \K}(X_1, X_2) = \inf_{f : J^k_{F_1} \leftrightarrow J^k_{F_2}}  \sup_{j \in \text{\em dom}(f)}\{ | j -f(j) | \}
\]
holds for $k =0$ and $1$, where $J^k_{F_i}$ denotes the multiset defined by \[\displaystyle{\coprod_{l=2m+k \ \text{with} \  H^l(F_i; \K) \neq 0} \left(\coprod_{\dim H^l(F_i; \K)} \{ \lfloor \frac{l}{2}  \rfloor \}\right)}.
\]
\end{cor}

\begin{proof}
The collapsing of the spectral sequence for $\mathcal{F}_i$ yields that the barcode $B_i$ associated with $H^*(X_i; \K)$ consists
of infinite intervals $[\lfloor \frac{l}{2} \rfloor, \infty)$ with $\dim H^l(F_i; \K) \neq 0$. We observe that each barcode $B_i$ is finite.
Then, the result follows from Theorem \ref{thm:I_and_B} and Lemma \ref{lem:finiteIntervals_d} (3).
\end{proof}

We conclude this section with a result which describes an upper bound of the cohomology interleaving distance between manifolds.
It is worthwhile mentioning that a map between the manifolds gives rise to one of the interleavings which induce the upper bound.

\begin{prop}\label{prop:shriek_map} \todo{[RR1, 2.3] We put a remark on the inequality together with [BL14, Corollary 6.6] after the proof of the proposition.}
Let $v_X : X\to BS^1$ and $v_Y : Y\to BS^1$ be connected closed oriented manifolds over $BS^1$.
Suppose that there exists a continuous map $f : X\to Y$ with $v_Y\circ f = v_X$.
Then
\begin{itemize}
\item[{\em (i)}] $d_{\text{\em CohI}, \K} (X, Y) \leq \frac{1}{2}(\dim Y - \dim X)$ if $\dim X$ and $\dim Y$ are even and $\dim Y \geq 2 \dim X$, and
\item[{\em (ii)}]  $d_{\text{\em CohI}, \K} (X, Y) < \frac{1}{2}(\dim Y - \dim X)$ if $\dim X$ and $\dim Y$ are odd and $\dim Y > 2 \dim X$.
\end{itemize}
\end{prop}

Before proving the result, we recall a {\it $\delta$-trivial} persistence module $M$ which satisfies the condition that
$M(i \to i+ \delta): M(i) \to M(i+\delta)$ is trivial for any $i$. Moreover, in the proof of Proposition \ref{prop:shriek_map},
we use the following lemma on a pointwise finite dimensional (p.f.d) persistence module $M$; that is, each of the vector space $M(i)$ is of finite dimension.

\begin{lem}\label{lem:n--s_condition}{\rm (\cite[Corollary 6.6]{Ba-L})} Two p.f.d persistence modules $M$ and $N$ are $\delta$-interleaved if and only if there exists a morphism $g : M \to N^\delta$ with ${\rm Ker} \ \!g$ and ${\rm Coker} \  \!g$ both $2\delta$-trivial.  Here $( \ )^{\delta}$ denotes the shift functor defined in Section \ref{sect:Interleavings}.
\end{lem}
\todo{[RR1, 2.3] We recall a result due to Bauer and Lesnick.}

\begin{proof}[Proof of Proposition \ref{prop:shriek_map}]
Let $m$ be the non-negative integer $\dim Y - \dim X$. The shriek map
$f^!$ is an element of $\text{Ext}_{C^*(Y;  \K)}^m(C^*(X; \K), C^*(Y; \K))$
which assigns the volume form of $Y$ to that of $X$, where the $\text{Ext}$ group is defined in the derived category of
$C^*(Y; \K)$-modules; see, for example, \cite{F-T}.
We have the composite map $\K[u] \stackrel{\kappa}{\to} C^*(BS^1; \K) \stackrel{v_Y^*}{\to} C^*(Y; \K)$ of morphisms of dg algebras, where $H(\kappa)(u)$ is the generator of $H^*(BS^1; \K)$.
Then, the map $H(f^!) : H^*(X; \K) \to \Sigma^{m}H^*(Y ; \K)$ induced by shriek map $f^!$ is a morphism of $\K[u]$-modules.
Observe that the map $H(f^!)$ gives rise to map $H(f^!) : h^kC^*(X; \K) \to (h^kC^*(Y; \K))^{\frac{m}{2}}$ for each $k = 0, 1$
because $m$ is even.  Here, we regard the codomain of the functor $h^k$ as $\mathsf{Mod}_\K^{(\R, \leq)}$ suppressing the isomorphism $\gamma$ and the embedding $(\lfloor \ \rfloor)^*$ in the diagram (\ref{eq:diagrams-diagram}).

(i) In view of Lemma \ref{lem:n--s_condition}, in order to prove that $h^kC^*(X; \K)$ and $h^kC^*(Y; \K)$ for $k = 0$ and $1$ are $\frac{m}{2}$-interleaved,
it suffices to show that the kernel and the cokernel of $H(f^!)$ are $2(\frac{m}{2})$-trivial.
The shriek map $f^!$ preserve the volume forms.  \todo{[RR1, 2.3] We revised the proof adding a comment on the support of  $\text{Ker} \ \! H(f^!)$.}
Then, since  the module $\text{Ker} \ \! H(f^!)$ is supported on $[0, \frac{\dim X}{2})$, it follows that
$\text{Ker} \ \! H(f^!)$ is $2(\frac{m}{2})$-trivial if $2(\frac{m}{2}) \geq  \frac{\dim X}{2}$.  Moreover, we see that $\text{Coker} \ \! H(f^!)$ is $2(\frac{m}{2})$-trivial if $2(\frac{m}{2}) \geq  \frac{\dim Y}{2}$. Thus, the result follows from the assumption that $\dim Y \geq 2 \dim X$.

(ii) The same argument as in the proof of (i) enables us to obtain the result (ii). We observe that the maps $H(f^!)|_{h^kC^*(X; \K)}$ for $k= 0$ and $1$ are $2(\frac{m}{2})$-trivial if $\dim Y > 2 \dim X$.
\end{proof}

The following remark is based on a comment due to a referee.
\begin{rem}
Let $X$ and $Y$ be manifolds  in Proposition \ref{prop:shriek_map} (i), where we do not assume the existence of a map between them.
Let $F^k(X)$ and $F^k(Y)$ be the persistence modules $h^kC^*(X; \K)$ and  $h^kC^*(Y; \K)$, respectively.
Then, the persistence module  $F^k(X)$ is supported on the interval $[0 , \frac{1}{2}\dim X + 1)$ and
$F^k(Y)$ is supported on the interval $[0 , \frac{1}{2}\dim Y + 1)$. Thus, it follows that
the zero morphism $F^k(X) \to F^k(Y)^{\frac{1}{4}\dim Y + \frac{1}{2}}$ is part of a $(\frac{1}{4}\dim Y + \frac{1}{2})$-interleaving. Then, we have
\[
d_{\rm I}(F^k(X), F^k(Y)) \leq \frac{1}{4}\dim Y + \frac{1}{2}
\left\{
\def\arraystretch{1.1}
\begin{array}{ll}
\! < \frac{1}{2}(\dim Y - \dim X) \  \  \text{if}\  \dim Y > 2\dim X +2 \\
\! = \frac{1}{2}(\dim Y - \dim X)  \  \  \text{if}\  \dim Y =  2\dim X +2\\
\! = \frac{1}{2}(\dim Y - \dim X) + \frac{1}{2}  \ \  \text{if}\  \dim Y =  2\dim X.
\end{array}
\right.
\]
In the case where $\dim Y = 2\dim X$, our upper bound  is better than the trivial one although our bound is worse than the trivial one
in other cases.
\end{rem}

\section*{Acknowledgments}
The authors thank Grigory Solomadin for discussions with the first author by which relationships between the cup-length and the cohomology interleaving distance of spaces are realized. They are also grateful to Nicolas Berkouk and Ling Zhou for their valuable comments on the first version of this manuscript and for providing additional references on variants of the interleaving distance and homotopy invariants appearing in persistence theory.
The authors are grateful to anonymous two referees for careful reading of the first version of the article and many valuable comments. One of the referees has suggested a general lemma on distances. The suggestion enables authors to revise the proof of Theorem \ref{thm:Equalities} with Lemma \ref{lem:An_isometry}, which is obtained by modifying the idea of the referee. 

The first author was partially supported by JSPS KAKENHI
Grant Numbers JP23K20795, JP26K00604. The second author was partially supported by JSPS
KAKENHI Grant Numbers JP23K03097, JP26K06813. The fourth author was partially supported
by JSPS KAKENHI Grant Number JP25K17250.

\appendix

\section{Computational examples of the CohID}\label{sect:toy_examples}\label{sect:appA}
This section is devoted to computing the cohomology interleaving distance explicitly using interval decompositions.

First, we compute the cohomology interleaving distances $d^k_{{\rm CohI},{\mathbb K}}(M,{\mathbb K})$ and $d^k_{{\rm CohI},{\mathbb K}}(M,{\mathbb K}[u]/(u^2))$ for a dg ${\mathbb K}[u]$-module $M$ whose cohomology is of finite dimension.
Here, ${\mathbb K}$ and ${\mathbb K}[u]/(u^2)$ are regarded as dg ${\mathbb K}[u]$-modules with zero differentials.

Let ${\mathcal B}_1^k$ and ${\mathcal B}_2^k$ be the barcodes associated with $\utot{k} ({\mathbb K})$ and $\utot{k} ({\mathbb K}[u]/(u^2))$, respectively; see Section \ref{sect:ID_of_dgK[u]Mod} for the functor $\utot{k}$.
In particular, ${\mathcal B}_1^0 = \{ [0,1)\}$, ${\mathcal B}_2^0 = \{ [0,2)\}$ and ${\mathcal B}_1^1 ={\mathcal B}_2^1 = \emptyset$.
According to the interval decomposition, the two distances mentioned above coincide respectively with the bottleneck distances $d_{{\rm B}}({\mathcal B}^k_M , {\mathcal B}^k_1)$ and $d_{{\rm B}}({\mathcal B}^k_M , {\mathcal B}^k_2)$ for the barcode ${\mathcal B}^k_M$ associated with $\utot{k} H^*(M)$.

As can be seen from the calculations below, these bottleneck distances are related to the {\it cup length}.
For the graded ${\mathbb K}[t]$-module $\utot{k} H^*(M)$, 
we denote ${\rm cup}^k(M)$\label{index:cup_M} by the greatest non-negative integer $n$ such that $(t \times)^n \neq 0$ on $\utot{k} H^*(M)$. 
For convenience, we set ${\rm cup}^k(M)=-1$ if $\utot{k} H^*(M)$ is the zero module.
Observe that ${\rm cup}^0(M)$ concerns the cup-length of spaces over $BS^1$; see the paragraph before Proposition \ref{prop:X_pt} for details.
In the following propositions, we write $l^k := {\rm cup}^k(M)$ for short.


Let ${\mathcal B}^k_M$  be the barcode
$\{ [b_{\lambda} , b_{\lambda} + c_{\lambda}) \mid \lambda \in \Lambda \}$
associated with $\utot{k} H^*(M)$, where the index set $\Lambda$ is finite.
Then, it is readily seen that
${\rm cup}^k(M) + 1 = \max \{ c_{\lambda} \mid \lambda \in \Lambda \}$.
We set $\Lambda_{i} :=\{ \lambda \in \Lambda \mid c_{\lambda} = i \}$ and $\Lambda_{0,i} :=\{ \lambda \in \Lambda_{i} \mid b_{\lambda} = 0 \}$.

\begin{prop}\label{prop:M_pt}
Let $M$ be a dg ${\mathbb K}[u]$-module concentrated in non-negative degrees whose cohomology is of finite dimension. Then, for $k=0,1$ it holds that
\[
d^k_{{\rm CohI},{\mathbb K}}(M,{\mathbb K})
= \left\{
\begin{array}{cl}
0 & (H^*(M)\cong {\mathbb K})\\
\frac{1}{2} & (k=0 \ \ \text{and} \ \  \utot{0}H^*(M) = \{ 0 \})\\
\frac{1}{2}(l^k + 1) & (\text{otherwise}).
\end{array}
\right.
\]
\end{prop}

\begin{proof}[Proof of Proposition \ref{prop:M_pt}]
If $H^*(M)\cong {\mathbb K}$ , then it is immediate that $d^k_{{\rm CohI},{\mathbb K}}(M,{\mathbb K}) = 0$.
In what follows, we assume that $H^*(M) \not\cong {\mathbb K}$.
By virtue of Theorem \ref{thm:I_and_B}, it suffices to compute the bottleneck distance $d_{{\rm B}}({\mathcal B}^k_M , {\mathcal B}^k_1)$ instead of $d^k_{{\rm CohI},{\mathbb K}}(M,{\mathbb K})$.


For $k=1$, it is readily seen that
\[
d_{{\rm B}}({\mathcal B}^1_M , {\mathcal B}^1_1) = \inf_{h : {\mathcal B}_M^1 \leftrightarrow {\mathcal B}_1^1}\sup_{J\in {\rm dom}(h)} d_{{\rm I}}(\chi_{J}, \chi_{\emptyset}) =\inf_{h : {\mathcal B}_M^1 \leftrightarrow {\mathcal B}_1^1} \left\{ \frac{1}{2}(l^1 + 1) \right\} = \frac{1}{2}(l^1 + 1).
\]

Next, we consider the case $k=0$.
If $\utot{0}H^*(M) = \{ 0 \}$, then by Lemma \ref{lem:finiteIntervals_d}(1), the bottleneck distance $d_{{\rm B}}({\mathcal B}^0_M , {\mathcal B}^0_1)$ is $1/2$.
We now assume that $\utot{0}H^*(M) \neq \{ 0 \}$, that is, ${\mathcal B}_M^0 \neq \emptyset$.
Let $I_{\lambda} := [b_{\lambda}, b_{\lambda}+c_{\lambda})$ in ${\mathcal B}^0_M$.
Then, Lemma \ref{lem:finiteIntervals_d}(2) enables us to deduce that
\begin{align}
d_{{\rm I}}(\chi_{I_{\lambda}}, \chi_{[0,1)})
= \left\{
\begin{array}{cl}
c_{\lambda}-1 & (b_{\lambda}=0, c_{\lambda}=1,2)\\
\frac{c_{\lambda}}{2} & (\text{otherwise}).
\end{array}
\right.
\end{align}
Given a bijection $h : {\mathcal B}_M^0 \leftrightarrow {\mathcal B}_1^0$ with $h(I_{\lambda})=[0,1)$.
First, consider the case $\lambda \in \Lambda_{0,1}$. Since $H^*(M) \not\cong {\mathbb K}$, it follows that $\Lambda \setminus \{ \lambda \} \neq \emptyset$. Thus, we have
\begin{align}
\sup_{J\in {\rm dom}(h)} d_{{\rm I}}(\chi_{J}, \chi_{h(J)})
&= \max \left\{  d_{{\rm I}}(\chi_{I_{\lambda}}, \chi_{[0,1)}), \ d_{{\rm I}}(\chi_{I_{\mu}}, \chi_{\emptyset}) \mid \mu \in \Lambda \setminus \{ \lambda \} \right\}
\\
&= \max \left\{ 0,  \frac{c_{\mu}}{2} \ \middle|  \ \mu \in \Lambda \setminus \{ \lambda \} \right\} = \frac{1}{2}(l^0 + 1).
\end{align}
If $\lambda \in \Lambda_{0,2}$, then
\begin{align}
\sup_{J\in {\rm dom}(h)} d_{{\rm I}}(\chi_{J}, \chi_{h(J)})
&= \max \left\{  d_{{\rm I}}(\chi_{I_{\lambda}}, \chi_{[0,1)}), \ d_{{\rm I}}(\chi_{I_{\mu}}, \chi_{\emptyset}) \ \middle|  \ \mu \in \Lambda \setminus \{ \lambda \} \right\}
\\
&= \max \left\{ 1,  \frac{c_{\mu}}{2} \ \middle|  \ \mu \in \Lambda \setminus \{ \lambda \} \right\}
= \frac{1}{2}(l^0 + 1).
\end{align}
Observe that $l^0 \geq 1$ in the case where $\Lambda_{0,2} \neq \emptyset$.
Furthermore, if $\lambda \in \Lambda \setminus \left( \Lambda_{0,1} \cup \Lambda_{0,2} \right)$,
\begin{align}
\sup_{J\in {\rm dom}(h)} d_{{\rm I}}(\chi_{J}, \chi_{h(J)})
= \max \left\{ \frac{c_{\lambda}}{2},  \frac{c_{\mu}}{2} \ \middle|  \ \mu \in \Lambda \setminus \{ \lambda \} \right\}
=\frac{1}{2}(l^0 + 1).
\end{align}
The computations of suprema enables us to obtain the equality
\[
d_{{\rm B}}({\mathcal B}^0_M , {\mathcal B}^0_1) = \inf_{h : {\mathcal B}_M^0 \leftrightarrow {\mathcal B}_1^0}\sup_{J\in {\rm dom}(h)} d_{{\rm I}}(\chi_{J}, \chi_{h(J)}) = \frac{1}{2}(l^0 + 1).
\]
We have the result.
\end{proof}

With the same notation as above,  
we have the following result.

\begin{prop}\label{prop:M_CP1}
Let $M$ be a dg ${\mathbb K}[u]$-module concentrated in non-negative degrees whose cohomology is of finite dimension.
Then, one gets

\begin{itemize}
\item[\text{\rm (1)}]
$d^0_{\text{\em CohI},{\mathbb K}}(M , {\mathbb K}[u]/(u^2)) \! = \! \left\{
\def\arraystretch{1.1}
\begin{array}{cl}
1 & (l^0 =-1, 0)\\
l^0 - 1 & (\# \Lambda = 1, \ l^0 =1,2) \\
\frac{1}{2}l^0 & (\# \Lambda \geq 2, \ l^0 = i, \ \# \Lambda_{i}=1, \ i=1,2)\\
\frac{1}{2}(l^0 + 1) & (\text{otherwise}),
\end{array}
\right.$
\item[\text{\rm (2)}]
$d^1_{\text{\em CohI},{\mathbb K}}(M , {\mathbb K}[u]/(u^2))\! = \!\frac{1}{2}(l^1 + 1)$.
\end{itemize}
Here $\#S$ denotes the cardinal of a set $S$.
\end{prop}


\begin{proof}[Proof of Proposition \ref{prop:M_CP1}]
Since ${\mathcal B}^1_2 = \emptyset$, the assertion (2) follows from Theorem \ref{thm:I_and_B} and Lemma \ref{lem:finiteIntervals_d} (1).

Let $\Pi$ be the set of all bijections between ${\mathcal B}^0_M$ and ${\mathcal B}_2^0$, and $\Pi_{0,i}$ the subset of $\Pi$ consisting of bijections $h:{\mathcal B}_M^0 \leftrightarrow {\mathcal B}_2^0$ such that $h([b_{\lambda}, b_{\lambda}+c_{\lambda}))=[0,2)$ with $\lambda \in \Lambda_{0,i}$.
We write $\Pi_+$ for the complement of the union $\cup_i \Pi_{0,i}$ in $\Pi$.
By applying Theorem \ref{thm:I_and_B}, we see that the right-hand side of the equality in (1) coincides with
the bottleneck distance $d_{{\rm B}}({\mathcal B}^0_M,{\mathcal B}^0_2)$, which
is the smallest value of the infima
\[
{\mathcal I_+} := \inf_{h\in \Pi_+}  \sup_{J\in {\rm dom}(h)}d_{{\rm I}}(\chi_{J}, \chi_{h(J)}) \hspace{1em} \text{and} \hspace{1em}
 {\mathcal I_{0,i}} := \inf_{h\in \Pi_{0,i}}  \sup_{J\in {\rm dom}(h)}d_{{\rm I}}(\chi_{J}, \chi_{h(J)})
\]
for $i=1,2,\ldots , l^0 + 1$. To obtain the equality in (1), we determine the values of these infima.

For any barcode $I_{\lambda} := [b_{\lambda} , b_{\lambda} + c_{\lambda} )$ in ${\mathcal B}_M^0$, the assumption implies $b_{\lambda}\geq 0$. It follows from Lemma \ref{lem:finiteIntervals_d}  that
\begin{eqnarray}\label{prop:ID+}
d_{{\rm I}}(\chi_{I_{\lambda}}, \chi_{[0,2)})
=
\left\{
\begin{array}{cl}
1 & (c_{\lambda}=1)\\
\frac{1}{2}c_{\lambda} & (c_{\lambda} \geq 2)
\end{array}
\right.
\end{eqnarray}
for the case $b_{\lambda} \geq 1$, and
\begin{eqnarray}\label{prop:ID0i}
d_{{\rm I}}(\chi_{I_{\lambda}}, \chi_{[0,2)})
=
\left\{
\begin{array}{cl}
1 & (c_{\lambda}=1)\\
c_{\lambda}-2 & (c_{\lambda}=2,3)\\
\frac{1}{2}c_{\lambda} & (c_{\lambda} \geq 4)
\end{array}
\right.
\end{eqnarray}
for the case $b_{\lambda}=0$.
Given a bijection $h : {\mathcal B}_M^0 \leftrightarrow {\mathcal B}^0_2$.
If $h\in \Pi_+$, then \eqref{prop:ID+} yields
\begin{equation}
\sup_{J\in {\rm dom}(h)}d_{{\rm I}}(\chi_{J}, \chi_{h(J)})
= \left\{
\begin{array}{cl}
1 & (l^0 = 0)\\
\frac{1}{2}(l^0 + 1) & ( l^0 \geq 1).
\end{array}
\right.
\end{equation}
Hence, ${\mathcal I_+} = 1$ if $l^0 = 0$, and ${\mathcal I_+} = (l^0 + 1)/2$ if $l^0 \geq 1$.
We next consider the case $h\in \Pi_{0,i}$ with $h(I_{\lambda})=[0,2)$ for some $\lambda \in \Lambda_{0,i}$.
Then, the equality \eqref{prop:ID0i} enables us to deduce that
\begin{eqnarray}
\sup_{J\in {\rm dom}(h)}d_{{\rm I}}(\chi_{J}, \chi_{h(J)})
&=&
\max\left\{  d_{{\rm I}}(\chi_{I_{\lambda}}, \chi_{[0,2)}), \ d_{{\rm I}}(\chi_{I_{\mu}}, \chi_{\emptyset}) \ \middle|  \ \mu \in \Lambda \setminus \{ \lambda \} \right\}
\notag\\
&=& \left\{
\begin{array}{cl}
\max \left\{ 1, \frac{1}{2}c_{\mu} \ \middle| \ \mu \in \Lambda \setminus \{ \lambda \}  \right\} & (c_{\lambda}=1)\\
\max \left\{ c_{\lambda} - 2, \frac{1}{2}c_{\mu} \ \middle| \ \mu \in \Lambda \setminus \{ \lambda \}  \right\} & (c_{\lambda}=2,3)\\
\frac{1}{2}(l^0 + 1) & (c_{\lambda} \geq 4).\\
\end{array}
\right.
\label{prop:sup_ID0i}
\end{eqnarray}

Consider the case $\# \Lambda = 1$ with $\Lambda = \{ \lambda \}$.
We observe that $\Pi_{0,i} = \emptyset$ for $i=1,2,\ldots , l^0$.
Since $c_{\lambda}=l^0 + 1$ and $\Lambda \setminus \{ \lambda \} = \emptyset$, it follows from \eqref{prop:sup_ID0i} that
\[
{\mathcal I}_{0,l^0 + 1} = \left\{
\begin{array}{cl}
1 & (l^0 = 0)\\
l^0 - 1 & (l^0 = 1,2)\\
\frac{1}{2}(l^0 + 1) & ( l^0 \geq 3).
\end{array}
\right.
\]

In the case $\# \Lambda \geq 2$, we see that ${\mathcal I}_{0,i}=(l^0 + 1)/2$ for $l^0 \geq 1$ and $i =1,2,\ldots , l^0 $. Indeed, for any $h\in \Pi_{0,i}$ associated with
$I_{\lambda}$; that is, $h(I_\lambda) = [0, 2)$ and $I_\lambda = [0, i )$, the equality \eqref{prop:sup_ID0i} gives
\[
\sup_{J\in {\rm dom}(h)}d_{{\rm I}}(\chi_{J}, \chi_{h(J)}) = \frac{1}{2}(l^0 + 1)
\]
since there exists $\mu \in \Lambda \setminus \{ \lambda \}$ such that $c_{\mu}=l^0 + 1$.
Furthermore, it follows from \eqref{prop:sup_ID0i} that
\begin{eqnarray*}
{\mathcal I}_{0,l^0+1}
&=&\left\{
\begin{array}{cl}
\min \left[ \max\left\{ 1, \   \frac{1}{2}c_{\mu} \ \middle| \ \mu \in \Lambda \setminus \{ \lambda \}      \right\}  \  \middle|  \ \lambda \in \Lambda_{0,1} \right] & (l^0 = 0)\\
\min \left[ \max\left\{ l^0 -1, \   \frac{1}{2}c_{\mu} \ \middle| \ \mu \in \Lambda \setminus \{ \lambda \}      \right\}  \  \middle|  \ \lambda \in \Lambda_{0,l^0 + 1} \right] & (l^0 = 1,2 )\\
\frac{1}{2}(l^0 + 1) & (l^0 \geq 3)
\end{array}
\right.
\notag\\
&=& \left\{
\def\arraystretch{1.1}
\begin{array}{cl}
1 & (l^0 = 0)\\
\frac{1}{2}l^0 & (l^0 = 1,2, \ \# \Lambda_{l^0 + 1} = 1)\\
\frac{1}{2}(l^0 + 1) & (\text{otherwise}).
\end{array}
\right.
\end{eqnarray*}
We remark that $c_{\lambda}=1$ for any $\lambda \in \Lambda$ in the case $l^0 =0$. The condition $\# \Lambda_{l^0 + 1} = 1$ implies $\Lambda_{l^0 + 1} = \Lambda_{0,l^0 + 1}$ and the inequality $c_{\mu}<l^0 + 1$ for any $\mu \in \Lambda \setminus \Lambda_{0,l^0 + 1}$.
On the other hand, if $\# \Lambda_{l^0 + 1} \geq 2$, then for any $\lambda \in  \Lambda_{0,l^0 + 1}$, there exists $\mu \in \Lambda $ such that $c_{\mu}= l^0 + 1$ and $\mu \neq \lambda$.
Therefore, by taking the smallest value among ${\mathcal I}_+$ and ${\mathcal I}_{0,i}$ for $i=1,2,\ldots , l^0 + 1$ computed above, we have the assertion (1).
\end{proof}

By applying Proposition \ref{prop:dI-dB}, we give computational examples of the cohomology interleaving distances between complex projective spaces.



\begin{prop}\label{example:CP}
Let $f_{n,j} : {\mathbb C} P^n \to BS^1$ be a map which represents an integer $j$ under the identifications $[{\mathbb C} P^n , BS^1] \cong H^2 ({\mathbb C} P^n ; \Z) \cong \Z$.
Then, it holds that
\begin{enumerate}
\item[\text{\em (1)}] $d^0_{\text{\em CohI},\Q}(({\mathbb C} P^n , f_{n,1}) , ({\mathbb C} P^m , f_{m,1})) = \min \left\{ |n-m|, \max\left\{ \frac{m+1}{2}, \ \frac{n+1}{2}  \right\}    \right\}$,
\item[\text{\em (2)}]  $d^0_{\text{\em CohI},\Q}(({\mathbb C} P^n , f_{n,0}) , ({\mathbb C} P^n , f_{n,1})) = \left\lceil \frac{n}{2} \right\rceil$,
\item[\text{\em (3)}]  $d^0_{\text{\em CohI},\Q}(({\mathbb C} P^n , f_{n,0}) , ({\mathbb C} P^m , f_{m,0})) =
\left\{
\begin{array}{ll}
0 & (n=m), \\
\frac{1}{2} & (n \neq m).
\end{array}
\right.$
\end{enumerate}
Here $\lceil \ \rceil$\label{index:ceiling} denotes the ceiling function.
\end{prop}

\begin{rem}
Since the cohomology of ${\mathbb C}P^n$ is concentrated in even degrees, it follows that $d^1_{\text{CohI},\Q}(({\mathbb C} P^n , f_{n,j}) , ({\mathbb C} P^m , f_{m,j'}))=0$.
\end{rem}

To prove Proposition \ref{example:CP}, we now set up some notations.
Observe that the algebra map $f_{n,j}^* : \Q[u]\cong H^*(BS^1 ; \Q) \to H^*({\mathbb C} P^n ; \Q) \cong \Q[x]/(x^{n+1})$ induced by $f_{n,j}$ in rational cohomology satisfies the condition that
$f_{n,0}^* (u) = 0$ and $f_{n,1}^* (u) = x$, where $\deg x = 2$.
These $\Q[u]$-module structures give the $\Q[t]$-module structures on $\utot{0}H^*({\mathbb C} P^n ; \Q)$ and then the barcodes associated with the modules as in Section \ref{sect:BS^1}.
Let ${\mathcal B}_{n,j}$ denote the barcode obtained by $f_{n,j}^*$.
Then, it is readily seen that
\[
{\mathcal B}_{n,j} = \left\{
\begin{array}{ll}
\{ [0,1), [1,2), \ldots , [n,n+1) \} & (j=0), \\
\{ [0,n+1) \} & (j=1).
\end{array}
\right.
\]
For simplicity, we put $\chi_{n,j} = \chi ({\mathcal B}_{n,j})$.

\begin{proof}[Proof of Proposition \ref{example:CP}]
The assertion (1) follows immediately from Lemma \ref{lem:finiteIntervals_d} (2).
In view of Proposition \ref{prop:dI-dB},
in order to show (2), it suffices to determine the bottleneck distance $d_B ({\mathcal B}_{n,0}, {\mathcal B}_{n,1})$.
Given a bijection $h:{\mathcal B}_{n,0} \leftrightarrow {\mathcal B}_{n,1}$,
if $h^{-1}([0,n+1)) = [i,i+1)$ for some $i=1,2,\ldots , n$, then we have
\begin{align}
\sup_{I\in {\rm dom}(h)} d_{{\rm I}} (\chi_I , \chi_{h (I)}) &= d_{{\rm I}}(\chi_{[i,i+1)}, \chi_{[0, n+1)})
=\min \left\{ \max \{ i , n-i \} ,  \frac{n+1}{2}  \right\}
\end{align}
by Lemma \ref{lem:finiteIntervals_d} (1) and (2).
If $h^{-1}([0,n+1)) = \emptyset$, then Lemma \ref{lem:finiteIntervals_d} (1) shows that 
\[
\sup_{I\in {\rm dom}(h)}d_{{\rm I}} (\chi_I , \chi_{h (I)}) = d_{{\rm I}}(\chi_{\emptyset}, \chi_{[0, n+1)}) = \frac{n+1}{2}.
\]
Hence, we have 
\[
d_B ({\mathcal B}_{n,0}, {\mathcal B}_{n,1}) = \inf_{h:{\mathcal B}_{n,0} \leftrightarrow {\mathcal B}_{n,1}} \sup_{I\in {\rm dom}(h)} d_I (\chi_I , \chi_{h (I)})
=
\min_{1\leq i \leq n}\left\{ \max \{ i , n-i \} ,  \frac{n+1}{2}  \right\}.
\]
Observe that the right-hand side integer coincides with $\lceil n/2 \rceil$, which completes the proof for (2).

The assertion (3) for $n=m$ is trivial.
We consider the case where $n \neq m$.
Since $d_{{\rm I}} (\chi_{[i,i+1)}, \chi_{[j,j+1)})\leq 1/2$ and $d_I (\chi_{\emptyset}, \chi_{[j,j+1)}) = 1/2$ from Lemma \ref{lem:finiteIntervals_d}, we see that
\[
\sup_{I\in {\rm dom}(h)}d_{{\rm I}} (\chi_I , \chi_{h (I)}) = \frac{1}{2}
\]
for every bijection $h:{\mathcal B}_{n,0} \leftrightarrow {\mathcal B}_{m,0}$.
Therefore, we have $d_B ({\mathcal B}_{n,0}, {\mathcal B}_{m,0})=1/2$. Theorem \ref{thm:I_and_B} yields the result (3).
\end{proof}

In the rest of this section, we use terminology in rational homotopy theory; see Appendix \ref{sect:RHT}
for (relative) Sullivan models for spaces.

\begin{prop}\label{prop:models}
For each $j = 0, 1$,
let $v_j : Z_j \to BS^1$ be a space over $BS^1$ whose relative Sullivan model has the form $(\wedge u, 0) \to (\wedge (x, y, z, u), d)$ with
$dz = jxyu + u^4$ and $dx = 0 =dy$, where $\deg x = \deg y = 3$, $\deg z=7$ and $\deg u =2$. 
Then, one has
\[
d_{\text{\em CohI}, \Q}^0(Z_0, Z_1) =3 \  \ \ \text{and}  \ \  \ d_{\text{\em CohI}, \Q}^1(Z_0, Z_1) =0.
\]
\end{prop}

%
%
In order to prove Proposition \ref{prop:models}, we first determine the $\Q$-cohomology of $Z_j$ as a $\Q[u]$-module.
It is readily seen that $H^*(Z_0; \Q) \cong \wedge (x,  y) \otimes \Q[u]/(u^4)$
as an algebra. In order to compute the cohomology of $Z_1$ by using a spectral sequence, we filter the model $\mathscr{M}$ for $Z_1$
with the {\it weights} of elements in $\mathscr{M}$.
Define the weighs 
of elements $x$, $y$, $z$ and $u$ by
$\text{weight} (x) = \text{weight} (y) = \text{weight} (z) = 0$ and $\text{weight} (u) =2$. The weight of a monomial is defined by the sum of the weights of elements constructing the monomial. We define a decreasing filtration $F^*$ of the model $\mathscr{M}$ by
\[
F^i : = \{ w \in \mathscr{M} \mid \text{$w$ is a linear combination of monomials with $\text{weight} \geq i$}\}.
\]
Then, the filtration gives rise to the first quadrant multiplicative spectral sequence
$\{E_r^{*,*}, d_r\}$ conversing to $H^*(\mathscr{M})= H^*(Z_1; \Q)$.
We see that
\[
E_2^{*,*} \cong \wedge (x, y, z)\otimes \Q[u]
\]
and $d_2(z) = xyu$, $d_2(x) = 0= d_2(y)$. It follows that as a $\Q[u]$-module,
\[
E_3^{*,*}\cong \Q[u]\{1, x, y, xz, yz, xyz\}\oplus  (\Q[u]/(u))\{xy\}.
\]
The next nontrivial differentials $d_r$ are given by $d_8(xz) = xu^4$ and $d_8(yz) = yu^4$.
The element $xyz$ in the $E_8$-term represents the element $xyz - u^3z$ in $\mathscr{M}$. Therefore, we have $d_{8}(xyz) = d_{8}(xyz -u^3z)= 0$.
Thus, we see that as a $\Q[u]$-module,
\[
E_9^{*,*}\cong \Q[u]/(u^4)\{x, y\} \oplus  \Q[u]\{1, xyz\}\oplus  (\Q[u]/(u))\{xy\}.
\]
Since
$d_{14}(xyz) = d_{14}(xyz -u^3z)= u^7$, it follows that as $\Q[u]$-modules,
\[
E_\infty\cong E_{15}^{*,*}\cong \Q[u]/(u^4)\{x, y\} \oplus  \Q[u]/(u^7)\{1\}\oplus  (\Q[u]/(u))\{xy\}.
\]
Thus, Lemma \ref{lem:bigradedModules} implies that $H^*(Z_1; \Q)\cong \text{Tot} E_\infty^{*,*}$ as a $\Q[u]$-module.

\begin{proof}[Proof of Proposition \ref{prop:models}]
By applying the functors $h^0$ and $h^1$ in the diagram (\ref{eq:diagrams-diagram}), we see that
\begin{eqnarray*}
h^1(C^*(Z_0; \Q)) &\!\!\!\cong \!\!\!& \Sigma^{-1}(\Q[t]/(t^4))^{\oplus 2} \ \cong \ h^1(C^*(Z_1; \Q)),  \\
C_0:=h^0(C^*(Z_0; \Q))&\!\!\cong \!\!&\Sigma^0(\Q[t]/(t^4)) \oplus  \Sigma^{-3}(\Q[t]/(t^4))  \  \  \text{and} \\
C_1:=h^0(C^*(Z_1; \Q)) &\!\!\cong \!\!& \Sigma^0(\Q[t]/(t^7)) \oplus \Sigma^{-3}(\Q[t]/(t)).
\end{eqnarray*}
The results follow from the computation of the cohomology mentioned above.  Thus, we have the assertion on $d^1_{\text{CohI}, \Q}$

We prove the first equality. Observe that $d^0_{\text{CohI}, \Q}(Z_0, Z_1)= d_{\text{I}}(\chi\mathcal{B}_{C_0}, \chi\mathcal{B}_{C_1}) = d_{\text{B}}(\mathcal{B}_{C_0}, \mathcal{B}_{C_1})$; see Theorem \ref{thm:I_and_B}.
Let $I_1$, $I_2$, $I'_1$ and $I'_2$ be the interval modules in $C_0$ and $C_1$ corresponding the intervals $[0, 4)$, $[3, 7)$, $[0, 7)$ and $[3, 4)$, respectively.
It follows from Lemma \ref{lem:finiteIntervals_d} and Remark \ref{rem:sum} that
$d^0_{\text{CohI}, \Q}(Z_0, Z_1)= d_{\text{I}}(\chi\mathcal{B}_{C_0}, \chi\mathcal{B}_{C_1}) = d_{\text{B}}(\mathcal{B}_{C_0}, \mathcal{B}_{C_1})  \leq 3$.

Suppose that $C_0$ and $C_1$ are $\delta$-interleaved,
where $\delta < 3$. Then, there exist natural transformations $\varphi : C_0 \longleftrightarrow C_1 : \psi$ which give the $\delta$-interleaving. Since $I_2(i) = 0$ for $i < 3$, it follows that the nontrivial image of restriction $\psi : I'_1 \longrightarrow I_1 \oplus I_2$ is in $I_1$.
By the same reason for $I_2'$ as that for $I_2$, we see that
the nontrivial image of the restriction $\varphi : I_1 \longrightarrow I'_1 \oplus I'_2$ is in $I'_1$.
Thus, the restrictions of $\varphi$ and $\psi$ induce a $\delta$-interleaving $\varphi : I_1 \longleftrightarrow I_1' : \psi$. Therefore, we have a commutative diagram
\[
\xymatrix@C40pt@R15pt{
   & I_1(\delta) \ar[rd]^(0.6){\varphi(\delta)}&  I_1(0) \ar[l]_\cong \ar[rd]^(0.6){\varphi(0)} \ar[r]^{I_1(0\to 4)} & I_1(4)  \ar[rd]^(0.6){\varphi(4)} \\
I'_1(0) \ar[rr]_-{I_1'(0 \to 2\delta)}^\cong \ar[ru]^-{\psi(0)}&& I'_1(2\delta)  & I'_1(\delta) \ar[l]_\cong \ar[r]_-{I'_1(\delta \to 4+\delta)}^-\cong & I'_1(4+\delta).
}
\]
Observe that the horizontal arrows are isomorphisms  except for $I_1(0\to 4)$ and $I_1(0\to 4)=0$.
Therefore, the map $\varphi(\delta)$ is nontrivial and hence $\varphi(0)$ is. This yields that $I_1'(\delta \to 4+\delta)\circ \varphi(0)$ is nontrivial, which is a contradiction.
We have
$d^0_{\text{CohI}, \Q}(Z_0, Z_1) =3$.
\end{proof}

One might be interested in a relationship between $Z_j$ in Proposition \ref{prop:models} and an $S^1$-action and a higher dimensional torus action on a space.  The issue is dealt with in the following remark.

\begin{rem}\label{rem:realization}
In general, for a given relative Sullivan algebra of the form
$\iota : (\wedge (u), 0) \to (\wedge W\otimes \wedge(u), d)$, there exists a fibration $M \to X \to BS^1$ whose model is the given Sullivan algebra. In fact,
by \cite[Proposition 17.9]{FHT}, we have a fibration $|\iota|: |(\wedge W\otimes \wedge(u), d)| \to |(\wedge (u), 0)|$. The pullback of the fibration along the rationalization map
$l : BS^1 \to |(\wedge (u), 0)|$ gives rise to a commutative diagram
\[
\xymatrix@C15pt@R12pt{
M \ar[r]^\simeq  \ar@{=}[d]  & X' \ar[r] \ar[d]^q & ES^1 \ar[d]^p \\
M  \ar[r] \ar@{=}[d] & X \ar[r] \ar[d] & BS^1 \ar[d]^l \\
|(\wedge W, \overline{d})| \ar[r] &|(\wedge W\otimes \wedge(u), d)| \ar[r]_-{|\iota|} & |(\wedge (u), 0)|
}
\]
in which $p$ is the universal $S^1$-bundle and the right-hand upper squares is also  pullback.
The result \cite[Proposition 15.6]{FHT} yields that the map $\iota : (\wedge (u), 0) \to (\wedge W\otimes \wedge(u), d)$ is  the relative Sullivan model for $X \to BS^1$.
Since $ES^1$ is contractible, it follows that
$X'$ is weak homotopy equivalent to the fiber $M$.
Moreover, we see that $X$ is the orbit space of the $S^1$-space $X'$ with a free action.

For example, it follows that each space $Z_j$ in Proposition \ref{prop:models} is the orbit space of an $S^1$-space $Z'_j$ which is rationally homotopy equivalent to $S^3\times S^3\times S^7$.
In particular, the bundle $p=1\times \pi : Z_0'=(S^3 \times S^3) \times S^7 \to Z_0=(S^3\times S^3) \times \mathbb{C}P^3$ is given by the usual principal $S^1$-bundle
$\pi : S^7 \to \mathbb{C}P^3$.
Moreover, we see that
the free $S^1$-action on $Z'_0$ does not extend to any free $S^1\times S^1$-action.
This follows from Proposition \ref{prop:ranks} which computes the rational toral ranks of $Z_0$ and $Z_1$.
\end{rem}

\begin{rem}\label{rem:differences}
While the computation before the proof of Proposition \ref{prop:models} yields that $H^*(Z_0; \Q)\cong H^*(Z_1; \Q)$ as a graded vector space,
Proposition \ref{prop:ranks} in particular implies that the rational homotopy types of $Z_0$ and $Z_1$ are different from each other.
Moreover, Theorem \ref{thm:IHC=CohI} and Proposition  \ref{prop:models} enable us to deduce that $\alpha C^*(Z_0; \Q)$ is not isomorphic to $\alpha C^*(Z_1;  \Q)$
in the category $\text{Ho}(\mathsf{Ch}_\Q)^{(\R, \leq)}$.
\end{rem}

It may hold that $d_{\rm{CohI}}(X,Y)=0$ for spaces $X$ and $Y$ over $BS^1$ even
if $H^*(X; \Q)$ is not isomorphic to $H^*(Y;\Q)$ as an algebra. We describe such an example.

\begin{rem}
For $a\in \Q\backslash\{0\}$, 
let $p_a : X_a \to BS^1$ be a space over $BS^1$ whose relative Sullivan minimal model is given by
\[
\iota : (\wedge (u), 0) \to \mathscr{M}(X_a):=(\wedge (u,x,y,z),d_a)
\]
with $|x|=|u|=2$, $|y|=|z|=3$, $d_au=d_ax=0$, $d_ay=ux$,
$d_az=x^2+au^2$ and $\iota(u)  = u$.
We observe that
$A_a:=H^*(X_a;\Q)\cong \Q[u,x]/(ux,x^2+au^2)$ as an algebra. Moreover, it follows that $A_a \cong A_b$ as an algebra if and only if $ab^{-1}$ is in $\Q^2$; see
\cite[Proposition 3.2]{M-S}.
On the other hand, it is readily seen that $A_a \cong A_b$ as a $\Q[u]$-module for $a ,b \in \Q\backslash\{0\}$ and hence $C^*(X_a; \Q) \cong C^*(X_b; \Q)$ in
$\rm{D}(\Q[u])$.
Thus, there exist spaces over $BS^1$ with infinitely many different rational homotopy types one another such that their cohomology interleaving distances are zero.

The spaces $X_{-1}$ and $X_{1}$ are realized as spaces $\mathbb{C}P^2\sharp \mathbb{C}P^2 \to BS^1$ and
$\mathbb{C}P^2\sharp \overline{\mathbb{C}P^2}\to BS^1$ over $BS^1$, respectively, for each which the map from the connected sum is defined by the composite of the pinching map, the projection in the first factor and the map
$f_{2,1}$ in Proposition \ref{example:CP};
see \cite[Example 3.7]{FOT} for the Sullivan model of such a connected sum.
\end{rem}

\begin{rem}\label{rem:no_map} We consider a map between $(\mathbb{C} P^n,f_{n,1})$ and the space $Z_j$ over $BS^1$ in Proposition \ref{prop:models}.
The minimal model of $\mathbb{C} P^n$ is given by $\mathscr{M}(\mathbb{C} P^n)=(\wedge (u,w),d)$ where $dw=u^{n+1}$.
Therefore, if there is a map between the two spaces, it is one of the cases.
\begin{enumerate}
\item[\text{(1)}] $f:\mathbb{C}P^n\to Z_j$ $(j=0,1)$ whose Sullivan representative is given by $\mathscr{M}(f)(x)=\mathscr{M}(f)(y)=0$ and $\mathscr{M}(f)(z)=u^{3-n}w$ for $1\leq n\leq 3$.
\item[\text{(2)}]
$f:Z_0\to \mathbb{C}P^n$  whose Sullivan representative is given by  $\mathscr{M}(f)(w)=u^{n-3}z+au^{n-1}x+bu^{n-1}y$ ($a,b\in \Q$) for $n\geq3$.
\end{enumerate}
We refer the reader to \cite[Section 12 (c)]{FHT} for a Sullivan representative for a map.
\end{rem}

\begin{assertion} There is no morphism between $\mathbb{C}P^n \ (n>3)$ and $Z_1$ over $BS^1$.
\end{assertion}
\begin{proof}
Suppose that there is a morphism  $f:\mathbb{C}P^n\to Z_1$ of spaces over $BS^1$ for $n>3$.
Then, since $|w|=2n+1>7$, it follows that $\mathscr{M}(f)(z)=0$. However, $\mathscr{M}(f)$ is a morphism of DGAs with $\mathscr{M}(f)(u) = u$, which is a contradiction.

If  there is a morphism $f:Z_1\to \mathbb{C}P^n$ of spaces over $BS^1$,
then we have  $\mathscr{M}(f)(w)=u^{n-3}z+g(u,x,y)$ for some $g\in \Q [u] \otimes \wedge^+ (x,y)$.
It follows that $d(g)=0$. This contradicts that $\mathscr{M}(f)$ is a morphism of DGAs.
\end{proof}


\begin{prop}\label{prop:M_CP}     Let $f_{n,1}:{\mathbb C}P^n \to BS^1$ and $v_j : Z_j \to BS^1$ be the spaces over $BS^1$ described in Proposition \ref{example:CP} and \ref{prop:models}, respectively.
    Then, 
    \[
    d^0_{\text{\em CohI},\Q}((Z_0 , v_0) , ({\mathbb C} P^n , f_{n,1})) = \left\{
\begin{array}{cl}
2 & (1 \leq n \leq 5)\\
3 & (6 \leq n \leq 9)\\
n-6 & (10 \leq n \leq 13)\\
\frac{n+1}{2} & (14 \leq n ),
\end{array}
\right.
\]
\[
d^0_{\text{\em CohI},\Q}((Z_1 , v_1) , ({\mathbb C} P^n , f_{n,1})) = \left\{
\begin{array}{cl}
\frac{7}{2} & (1 \leq n \leq 2)\\
-n+6 & (3 \leq n \leq 5)\\
\frac{1}{2} & (n =6)\\
n-6 & (7 \leq n \leq 13)\\
\frac{n+1}{2} & (14 \leq n )
\end{array}
\right.
\]
\[
\text{and} \hspace{0.5cm} d^1_{\text{\em CohI},\Q}((Z_j , v_j) , ({\mathbb C} P^n , f_{n,1})) = 2.
\]
\end{prop}


\begin{proof}
First, we prove the first two equalities by computing the bottleneck distances.
Recall the barcode ${\mathcal B}_{n,1}=\{ [0,n+1) \}$ associated with $\utot{0}H^*({\mathbb C}P^n ; {\mathbb Q})$ described above. We also recall the barcodes associated with $C_j = \utot{0}H^*(Z_j;{\mathbb Q})$ in the proof of Proposition \ref{prop:models} which are given by
\[
{\mathcal B}_{C_0}=\{[0,4),[3,7) \} \hspace{1em} \text{and} \hspace{1em} {\mathcal B}_{C_1}=\{[0,7),[3,4) \},
\]
respectively.
Given a bijection $h:{\mathcal B}_{C_0} \leftrightarrow {\mathcal B}_{n,1}$,
if $h([0,4))=[0,n+1)$, then Lemma \ref{lem:finiteIntervals_d} (1) and (2) allow us to deduce that
\begin{eqnarray}\label{ID1:M0_CP}
\sup_{J\in {\rm dom}(h)}d_{{\rm I}}(\chi_{J}, \chi_{h(J)}) &=& \max \{  d_{{\rm I}}(\chi_{[0,4)},\chi_{[0,n+1)}), d_{{\rm I}}(\chi_{[3,7)},\chi_{\emptyset} ) \}
\\
&=& \left\{
\begin{array}{cl}
2 & (1 \leq n \leq 5)\\
n-3 & (6 \leq n \leq 7)\\
\frac{n+1}{2} & (8 \leq n ).
\end{array}
\right.
\notag
\end{eqnarray}
Similarly, it is readily seen that
\begin{eqnarray}\label{ID2:M0_CP}
\sup_{J\in {\rm dom}(h)}d_{{\rm I}}(\chi_{J}, \chi_{h(J)}) = \left\{
\begin{array}{cl}
2 & (1 \leq n \leq 3)\\
\frac{n+1}{2} & (4 \leq n \leq 5)\\
3 & (6 \leq n \leq 9)\\
n-6 & (10 \leq n \leq 13)\\
\frac{n+1}{2} & (14 \leq n )
\end{array}
\right.
\end{eqnarray}
in the case where $h([3,7))=[0,n+1)$, and
\begin{eqnarray}\label{ID3:M0_CP}
\sup_{J\in {\rm dom}(h)}d_{{\rm I}}(\chi_{J}, \chi_{h(J)}) = \left\{
\begin{array}{cl}
2 & (1 \leq n \leq 3)\\
\frac{n+1}{2} & (4 \leq n )
\end{array}
\right.
\end{eqnarray}
in the case where $h(\emptyset)=[0,n+1)$.
Since the distance $d_{{\rm B}}({\mathcal B}_{C_0}, {\mathcal B}_{n,1})$ is the smaller value of \eqref{ID1:M0_CP}, \eqref{ID2:M0_CP} and \eqref{ID3:M0_CP}, the result for $d^0_{\text{\rm CohI},\Q}((Z_0 , v_0) , ({\mathbb C} P^n , f_{n,1}))$ is shown from Theorem \ref{thm:I_and_B}.
By the same argument above, we compute the bottleneck distance between ${\mathcal B}_{C_1}$ and ${\mathcal B}_{n,1}$, which completes the proof of (1).

More precisely, let $h:{\mathcal B}_{C_1} \leftrightarrow {\mathcal B}_{n,1}$ be a bijection satisfying $h([0,7))=[0, n+1)$. Then, we have
\begin{eqnarray}\label{ID1:M1_CP}
\sup_{J\in {\rm dom}(h)}d_{{\rm I}}(\chi_{J}, \chi_{h(J)}) = \left\{
\begin{array}{cl}
\frac{7}{2} & (1 \leq n \leq 2)\\
-n+6 & (3 \leq n \leq 5)\\
\frac{1}{2} & (n=6)\\
n-6 & (7 \leq n \leq 13)\\
\frac{n+1}{2} & (14 \leq n ).
\end{array}
\right.
\end{eqnarray}
Similarly, we have
\begin{eqnarray}\label{ID2:M1_CP}
\sup_{J\in {\rm dom}(h')}d_{{\rm I}}(\chi_{J}, \chi_{h'(J)}) = \left\{
\begin{array}{cl}
\frac{7}{2} & (1 \leq n \leq 6)\\
\frac{n+1}{2} & (7 \leq n )
\end{array}
\right.
\end{eqnarray}
for a bijection $h':{\mathcal B}_{C_1} \leftrightarrow {\mathcal B}_{n,1}$ satisfying $h'([3,4))=[0,n+1)$, and
\begin{eqnarray}\label{ID3:M1_CP}
\sup_{J\in {\rm dom}(h'')}d_{{\rm I}}(\chi_{J}, \chi_{h''(J)}) = \left\{
\begin{array}{cl}
\frac{7}{2} & (1 \leq n \leq 6)\\
\frac{n+1}{2} & (7 \leq n )
\end{array}
\right.
\end{eqnarray}
for a bijection $h'':{\mathcal B}_{C_1} \leftrightarrow {\mathcal B}_{n,1}$ satisfying $h''(\emptyset)=[0,n+1)$.
Since the bottleneck distance $d_{{\rm B}}({\mathcal B}_{C_1}, {\mathcal B}_{n,1})$ is the smaller value of \eqref{ID1:M1_CP} and \eqref{ID2:M1_CP}, Theorem \ref{thm:I_and_B} shows the assertion on  $d^0_{\text{\rm CohI},\Q}((Z_1 , v_1) , ({\mathbb C} P^n , f_{n,1}))$.

It follows from the ${\mathbb Q}[t]$-module structure of $h^1(C^*(Z_j;{\mathbb Q}))=\utot{1}H^*(Z_j ; \Q)$ described in the proof of Proposition \ref{prop:models} that the barcode associated with $\utot{1}H^*(Z_j;{\mathbb Q})$ is given by $\{ [1,5), [1,5) \}$ for $j=0$ and $1$.
On the other hand, the barcode associated with $\utot{1}H^*({\mathbb C}P^n;{\mathbb Q})$ is the empty set since the cohomology of ${\mathbb C}P^n$ is concentrated in even degrees.
Therefore, we have 
\[
d^1_{\text{\rm CohI},\Q}((Z_j , v_j) , ({\mathbb C} P^n , f_{n,1})) = d_{{\rm I}}(\chi_{[1,5)}, \chi_{\emptyset}) = 2
\]
This follows from Lemma \ref{lem:finiteIntervals_d} (2).
\end{proof}

\begin{rem}\label{rem:CP}
Let $X$ and $Y$ be spaces over  $BS^1$. Then, as seen in the proof of Proposition \ref{prop:cup_2},
the triangle inequality of the interleaving distance allows us to deduce an inequality
\[
\big| \ d_{\text{CohI}, \K}^k(X, \mathbb{C}P^n) - d_{\text{CohI}, \K}^k(Y, \mathbb{C}P^n) \ \big| \leq d_{\text{CohI}, \K}^k(X, Y)
\]
for each $n\geq 1$, $k = 0, 1$ and an arbitrary field $\K$.
The computation of the distance
$d_{\text{CohI}, \K}^k(X, \mathbb{C}P^n)$ in the proof of Proposition \ref{prop:M_CP} is comprehensively easy with the bottleneck distance because
the barcode associated with  $\mathbb{C}P^n$ consists of one bar. This is an advantage of the lower bound of the interleaving distance mentioned above.
\end{rem}

\section{Some rational homotopy invariants and the CohID}\label{sect:RHT}

We begin by briefly reviewing a (relative) Sullivan algebra.
Let
$\mathscr{M}(X)=(\wedge {V},d)$\label{index:minimalModel}
be the Sullivan minimal model of a simply-connected space $X$; see \cite[Section 12]{FHT}.
It is a free commutative differential graded algebra (DGA) over $\Q$ with a locally finite $\Q$-graded vector space $V=\bigoplus_{i\geq 1}V^i$ and a decomposable differential;
that is,
\[
\dim V^i<\infty,  d(V^i) \subset (\wedge^+{V} \cdot \wedge^+{V})^{i+1} \  \text{and} \ d \circ d=0.
\]
Here  $\wedge^+{V}$ denotes the ideal of $\wedge{V}$ generated by elements of positive degree.
The degree of a homogeneous element $x$ of the graded algebra is denoted by $|{x}|$.
By definition, the commutativity of the model gives the formula $xy=(-1)^{|{x}||{y}|}yx$ and the differential $d$ satisfies the condition that $d(xy)=d(x)y+(-1)^{|{x}|}xd(y)$ for homogeneous elements $x$ and $y$ in $\wedge V$.
Note that $\mathscr{M}(X)$ determines the rational homotopy type of $X$.
In particular, we see that $V^*\cong {\rm Hom}(\pi_*(X),\Q) \mbox{\ \ and\ \ }H^*(\wedge {V},d)\cong H^*(X;\Q )$.

Let $f:X\to Y$ be a map between simply-connected spaces. Then,
the relative Sullivan model of $f$
is given by $$\mathscr{M}(Y)=(\wedge W,d_Y)\to (\wedge W\otimes \wedge V,D)\to (\wedge V,\overline{D}),$$
where
$D\mid_W=d_Y$ and $(\wedge W\otimes \wedge V,D)$ is quasi-isomorphic to $\mathscr{M}(X)$ \cite[Section 14]{FHT}.

We also recall a spectral sequence introduced in \cite[Section 32 (b)]{FHT}.  Let $(\wedge V, d)$ be a Sullivan algebra for which $V$ is finite dimensional.  We give the Sullivan algebra a bigrading $(\wedge V, d)^{*,*}$ defined by
$(\wedge V^{\text{even}}\otimes \wedge^kV^{\text{odd}})^n= (\wedge V)^{k+n, k}$. Then, a generator $x$ with odd degree and
a generator $y$ with even degree have the bidegrees $(\deg x +1, -1)$ and $(\deg y, 0)$, respectively.
The filtration $F^*(\wedge V)$ of $\wedge V$ defined by $F^p(\wedge V) = (\wedge V)^{\geq p, *}$ gives rise to the forth quadrant spectral sequence converging to $H(\wedge V, d)$, which is called the {\it odd spectral sequence} of the Sullivan algebra $(\wedge V, d)$. Observe that the $E_0$-term is a DGA of the form $(\wedge V, d_\sigma)$ with the differential of bidegree $(0, +1)$ characterized by
\[
d_\sigma(V^\text{even}) =0, d_\sigma : V^\text{odd} \to \wedge V^\text{even} \ \text{and} \ d- d_\sigma : V^\text{odd} \to \wedge V^\text{even}\otimes \wedge^+V^\text{odd}.
\]
\begin{prop}\label{prop:OddSS}\cite[Proposition 32.4]{FHT}
Let $(\wedge V, d)$ be a minimal Sullivan algebra in which $V$ is of finite dimension and $V = V^{\geq 2}$. Then the following three conditions are equivalent. \text{\rm (i)} $\dim E_1 = \dim H(\wedge V, d_\sigma) < \infty$, \text{\rm (ii)} $\dim H(\wedge V, d) < \infty$  and \text{\rm (iii)} the LS category $\text{\em cat}(\wedge V, d)$  is finite; see \cite[Section 29]{FHT} for the definition of $\text{\em cat}(\wedge V, d)$.
\end{prop}

Let $ r_0(X)$ be the {\it rational toral rank}
of a simply-connected CW complex $X$
of $\dim H^*(X;\Q )<\infty $; that is, the largest integer $r$ such that an $r$-torus
$T^r=S^1 \times\dots\times S^1$ ($r$-factors) can act continuously
on a CW-complex $Y$ having the rational homotopy type of $X$
with all its isotropy subgroups finite (almost free action); see
\cite[7.3]{FOT} and \cite{Ha}. If an $r$-torus $T^r$ acts on $X$
by $\rho :T^r\times X\to X$, then the Borel fibration 
\[
X \to ET^r \times_{T^r}^{\rho} X \to BT^r
\]
is constructed.
Thus, we have a relative Sullivan model
\[
(\Q[u_1,u_2,\dots,u_r],0)
 \to (\Q[u_1,u_2,\dots,u_r] \otimes \wedge {V},D)
 \to (\wedge {V},d)\ \ \ \ (*)_r
 \]
for the fibration, where $\deg {u_i}=2$ for $i=1,2,\dots,r$, $Du_i=0$ and
$Dv \equiv dv$ modulo the ideal $(u_1,u_2,\dots,u_r)$ for $v\in V$. 
According to \cite[Proposition 4.2]{Ha}, $r_0(X) \geq r$ if and only if there exists a relative Sullivan algebra of the form $(*)_r$
such that $(\wedge {V},d)$ is the minimal model for $X$ and $\dim H(\Q[u_1,u_2,\ldots, u_r] \otimes  \wedge V,D) < \infty$.

We recall the spaces $Z_0$ and $Z_1$ over $BS^1$ in Proposition \ref{prop:models}.

\begin{prop}\label{prop:ranks}
 $r_0(Z_0)=2$ and $r_0(Z_1)=0$.
\end{prop}

\begin{proof}  It follows that $r_0(Z_0)\geq 2$.  In fact, we define $Dx=u_1^2$, $Dy=u_2^2$ and
$Dz=dz= u^4$ in $(*)_2$. Then, we have $\dim H(\Q[u_1, u_2] \otimes  \wedge V,D)< \infty$.  
If $r_0(M_0)\geq 3$, then there is a relative Sullivan model
\[
(\Q[u_1,u_2,u_3],0)
 \to (\Q[u_1,u_2,u_3] \otimes \wedge {(u, x,y,z)},D)
 \to (\wedge {(u, x,y,z)},d)\ \ \ \ (*)_3
 \]
such that $\dim H^*(\Q[u_1,u_2,u_3] \otimes \wedge {(u, x,y,z)},D)<\infty$. We write
$(\wedge W, D)$ for $(\Q[u_1,u_2,u_3] \otimes \wedge {(u, x,y,z)},D)$. Then, the result \cite[Proposition 32.10]{FHT} implies that
$\dim W^\text{odd} - \dim W^\text{even} \geq 0$. However, it follows from $(*)_3$ that $\dim W^\text{odd} - \dim W^\text{even} = 3-4=-1$, which is a contradiction.

Suppose that $r_0(Z_1)\geq 1$. Then, the DGA $(\wedge W, D):=(\wedge (x,y)\otimes \wedge (u_1, u, z),D)$ in $(*)_1$ for $Z_1$ satisfies the condition that $Dx=Dy=0$.
Indeed, let $Dz=xyu+u^4+f+axy u_1$ for some $f=f(u,u_1)\in \Q[u,u_1]$ and $a\in\Q$. We have
\[
DDz=gyu-hxu+agyu_1-ahxu_1\neq 0
\]
if $Dx=g(u,u_1)\neq 0$ or $Dy=h(u,u_1)\neq 0$ in $\Q [u,u_1]$.  Thus the differential $D$ is trivial on $\wedge (x, y)$. We consider the odd spectral sequence converging to $H(\wedge W, D)$. The $E_0$-term is a DGA of the form
$(\wedge (x, y), 0) \otimes \wedge (u_1, u, z), D_\sigma)$ with $D_\sigma z =  u^4 +f$. Thus, by applying \cite[Proposition 32.10]{FHT} again, we see that the $E_1$-term is of infinite dimension. Proposition \ref{prop:OddSS} implies that
$\dim H(\wedge W, D) = \infty$, which is a contradiction. We have  $r_0(Z_1)=0$.
\end{proof}



We conclude this section with comments on upper and lower bounds of the cohomology interleaving distance.
The proof of Proposition \ref{prop:models} enables us to deduce the following result.

\begin{prop}\label{prop:B}
Let $v_j : Z_j \to BS^1$ be the space over $BS^1$ in Proposition \ref{prop:models} for each $j=0$ and $1$. Then,
$\text{\rm cup}(v_0)_{\Q}=3$ and
$\text{\rm cup}(v_1)_{\Q}=6$.
\end{prop}

It follows that
the equalities in the inequalities in Proposition \ref{prop:cup} and Remark \ref{rem:CP} do not hold in general.
In fact, we have
\begin{eqnarray*}
& & \big| \ d_{\text{CohI}, \Q}^k(Z_0, \mathbb{C}P^6) - d_{\text{CohI}, \Q}^k(Z_1, \mathbb{C}P^6 ) \big| \\
&=&  3- \frac{1}{2} <
 d^0_{\text{CohI}, \Q}(Z_0, Z_1) =3 < \frac{7}{2}= \frac{1}{2}\text{max}\{\text{cup}(v_0)_\Q+1, \text{cup}(v_1)_\Q+1\}.
\end{eqnarray*}
Observe that the first equality follows from Proposition \ref{prop:M_CP}.

We have $\text{cup}_{\Q} (\mathbb{C} P^3)=\text{cup}_{\Q} (f_{3,1})=3$. Then, Proposition \ref{prop:M_CP} (1) and Proposition \ref{prop:cup} allow us to deduce  that
\[
d^0_{\text{CohI}, \Q}((\mathbb{C} P^3,f_{3,1}), (Z_0 , v_0 ))  = 2=\frac{4}{2}= \frac{1}{2}\text{max}\{3+1, \text{cup}(v_0)_\Q+1\}.
\]
On the other hand,
the inclusion $\mathbb{C}P^3 \to (S^3\times S^3)\times \mathbb{C}P^3=Z_0$ defined with a base point in $S^3\times S^3$ satisfies the assumption in Proposition \ref{prop:shriek_map} (i); see also Remark \ref{rem:realization}.
Thus, the evaluation in the proposition gives the inequality  $d^0_{\text{CohI}, \Q}(\mathbb{C} P^3, Z_0) \leq 3$.  The equality in Proposition \ref{prop:shriek_map} does not hold either in general.


\bigskip
\noindent
{\bf List of Symbols}
\begin{figure}[h]
  \begin{center}
    \begin{tabular}{r p{9cm} lc}
    symbol & meaning & page \\
      \hline
     $\mathcal{B}_V$  &    the barcode  associated with a persistence vector space $V$
      & \pageref{barcode}\\
      $\eta^k$ & the functor $\eta^k :  \mathsf{grMod}_\K^{(\R, \leq)} \to  \mathsf{Mod}_\K^{(\R, \leq)}$ & \pageref{index:eta}\\
      $\utot{k}$ & the functor $\utot{k} : \K[u]\text{-}\mathsf{grMod} \to \K[t]\text{-}\mathsf{grMod}$ & \pageref{index:S} \\
      $\nu^k$ & the functor $\nu^k:  \K[u]\text{-}\mathsf{Ch} \to \mathsf{Mod}_\K^{(\R, \leq)}$ & \pageref{index:nu} \\
        \(\alpha \) & the functor $\alpha : \K[u]\text{-}\mathsf{Ch} \to \mathsf{Ch}_\K^{(\R, \leq)}$ & \pageref{index:alpha} \\
             \(\text{cup}\),  \(\text{cup}^k\) & the cup-length of a map &  \pageref{index:cup1}
                    \pageref{index:cup2}\\
             \(\text{cup}^k\) & the cup-length of a graded $\K[u]$-module & \pageref{index:cup_M}\\
      \( d_{\text{CohI}} \) & the cohomology interleaving distance of persistence dg modules& \pageref{index:CohI} \\
         \( d_{\text{CohI}}^k \) & the (even, odd) cohomology interleaving distance of dg $\K[u]$-modules  &
        \pageref{index:eo_CohI} \pageref{index:CohI_total} \\
            \( d_{\text{CohI}, \K}^k\) & the (even, odd) cohomology interleaving distance of spaces &\pageref{index:CohI_space}  \\
             \(d_{\text{IHC}} \) & the interleaving distance in the homotopy category & \pageref{index:homotopyID} \\
      \( \lfloor \ \rfloor\), \(\lceil \ \rceil \) & the floor function, the ceiling function  & \pageref{index:floor}
      \pageref{index:ceiling}\\
      $\text{D}(\K[u])$ & the derived category of dg $\K[u]$-modules & \pageref{index:derivedcat}\\
           \(\mathscr{M}(X)\) & the Sullivan minimal model for a space $X$ & \pageref{index:minimalModel} \\
    \end{tabular}
  \end{center}
\end{figure}

\end{document}